\documentclass[11pt,reqno]{amsart}
\usepackage[percent]{overpic}
\usepackage{amsfonts,amssymb,amsmath,color}

\usepackage{bbm}

\usepackage[pdftex,colorlinks,
citecolor=black,linkcolor=black]{hyperref}

\DeclareMathOperator{\bmin}{\mathbf{min}} \DeclareMathOperator{\bmax}{\mathbf{max}}

 \DeclareMathOperator{\supp}{\mathrm{supp}}
 \DeclareMathOperator{\Arg}{\mathrm{Arg}}

\DeclareMathOperator{\hwt}{\mathtt{hwt}}

\DeclareMathOperator{\compl}{\scriptscriptstyle\complement}

\DeclareMathOperator{\ro}{\textcolor[rgb]{0.00,0.00,0.55}{\mathtt{r\hspace{-0.05mm}o}}}
\DeclareMathOperator{\rn}{\textcolor[rgb]{0.00,0.00,0.55}{\mathtt{rn}}}

\newtheorem{theorem}{Theorem}

\newtheorem{corollary}[theorem]{Corollary}

\newtheorem{proposition}[theorem]{Proposition}

\newtheorem{remark}[theorem]{Remark}
\newtheorem{example}[theorem]{Example}
\newtheorem{algorithm}[theorem]{Algorithm}

\numberwithin{equation}{section} \numberwithin{theorem}{section}

\begin{document}

\title{Pattern Recognition on Oriented Matroids: Symmetric Cycles in the Hypercube Graphs. V}

\author{Andrey O. Matveev}
\email{andrey.o.matveev@gmail.com}

\begin{abstract}
We consider decompositions of topes of the oriented matroid realizable as the arrangement of coordinate hyperplanes in $\mathbb{R}^{2^t}$, with respect to a distinguished symmetric $2\cdot 2^t$-cycle in its hypercube graph of topes $\boldsymbol{H}(2^t,2)$. We seek interpretations of such decompositions in the context of subset families on the ground set $E_t:=\{1,\ldots,t\}$ and of the families of their blocking sets, in the context of clutters on $E_t$ and of their blockers.
\end{abstract}

\maketitle

\pagestyle{myheadings}

\markboth{PATTERN RECOGNITION ON ORIENTED MATROIDS}{A.O.~MATVEEV}

\thispagestyle{empty}

\setcounter{tocdepth}{3}
\tableofcontents

\section{Introduction}

Let $\mathcal{H}:=(E_t,\{1,-1\}^t)$ be the oriented matroid on its {\em ground set\/}~$E_t:=[t]$ $:=[1,t]:=\{1,\ldots,t\}$, where~$t\geq 3$, and with its set of {\em topes\/} $\{1,-1\}^t$. This oriented matroid is realizable as the {\em arrangement\/} of {\em coordinate hyperplanes\/} in the real Euclidean space~$\mathbb{R}^t\supset\{1,-1\}^t$ of row vectors, see~\cite[Example~4.1.4]{BLSWZ}.

See, e.g.,~\cite{BK,Bo,BS,DRS,K,S,Z} on {\em oriented matroids}.

Each of the $2^t$ {\em maximal covectors\/} $T:=(T(1),\ldots,T(t))\in\{1,-1\}^t$ of~$\mathcal{H}$ can be regarded as the {\em characteristic tope\/} of the {\em negative part\/} $T^-:=\{e\in E_t\colon T(e)=-1\}$. Conversely, given an arbitrary subset $A\subseteq E_t$, we define the {\em characteristic tope\/} of~$A$ to be the {\em reorientation\/}~${}_{-A}\mathrm{T}^{(+)}$ of the {\em positive tope\/}~$\mathrm{T}^{(+)}:=(1,\ldots,1)$ on the subset~$A$; recall that~$({}_{-A}\mathrm{T}^{(+)})^-:=A$.
Let~$\boldsymbol{H}(t,2)$ denote the {\em hypercube graph\/} of {\em topes\/} of the oriented matroid~$\mathcal{H}$, that is, the {\em vertex\/} set of the graph~$\boldsymbol{H}(t,2)$ is the set~$\{1.-1\}^t$, and the {\em edges\/} of~$\boldsymbol{H}(t,2)$ are the
pairs $\{T',T''\}\subset\{1,-1\}^t$, such that $|\{e\in E_t\colon T'(e)\neq T''(e)\}|=1$.

Let $\boldsymbol{R}:=(R^0,R^1,\ldots,R^{2t-1},R^0)$ be a distinguished {\em symmetric cycle\/} in the graph~$\boldsymbol{H}(t,2)$, where
\begin{equation}
\label{eq:12}
\begin{split}
R^0:\!&=\mathrm{T}^{(+)}\; ,\\
R^s:\!&={}_{-[s]}R^0\; ,\ \ \ 1\leq s\leq t-1\; ,\\
\end{split}
\end{equation}
and
\begin{equation}
\label{eq:13}
R^{t+k}:=-R^k\; ,\ \ \ 0\leq k\leq t-1\; .
\end{equation}

For any vertex $T\in\{1,-1\}^t$ of the
graph~$\boldsymbol{H}(t,2)$, there exists a {\em unique inclusion-minimal\/} subset
\begin{equation}
\label{eq:37}
\boldsymbol{Q}(T,\boldsymbol{R})\subset\mathrm{V}(\boldsymbol{R}):=(R^0,R^1,\ldots,R^{2t-1})
\end{equation}
of the vertex sequence $\mathrm{V}(\boldsymbol{R})$ of the cycle~$\boldsymbol{R}$,
such that
\begin{equation}
\label{eq:32}
T=\sum_{Q\in\boldsymbol{Q}(T,\boldsymbol{R})}Q\; ,
\end{equation}
see~\cite[\S{}11.1]{M-PROM}. This subset $\boldsymbol{Q}(T,\boldsymbol{R})\subset\mathbb{R}^t$ is {\em linearly independent}, and it contains an {\em odd\/} number~$\mathfrak{q}(T):=\mathfrak{q}(T,\boldsymbol{R}):=|\boldsymbol{Q}(T,\boldsymbol{R})|$ of topes.
In fact, the linear algebraic decomposition~(\ref{eq:32}) is just a way to describe a particular mechanism of {\em majority voting}.

Let~$\boldsymbol{\sigma}(e)$ denote the $e$th standard unit vector
of the space~$\mathbb{R}^t$, $e\in[t]$. The bijections
\begin{align}
\label{eq:8}
\{1,-1\}^t&\to\{0,1\}^t\colon & T&\mapsto\tfrac{1}{2}(\mathrm{T}^{(+)}-T)\; ,\\
\intertext{and}
\label{eq:9}
\{0,1\}^t&\to\{1,-1\}^t\colon & \widetilde{T}&\mapsto\mathrm{T}^{(+)}-2\widetilde{T}\; ,
\end{align}
between the vertex set $\{1,-1\}^t$ of the hypercube graph~$\boldsymbol{H}(t,2)$ and the vertex set $\{0,1\}^t$ of the hypercube graph~$\widetilde{\boldsymbol{H}}(t,2)$ allow us to associate with the symmetric cycle~$\boldsymbol{R}$ in the graph~$\boldsymbol{H}(t,2)$ a symmetric cycle~$\widetilde{\boldsymbol{R}}:=(\widetilde{R}^0,\widetilde{R}^1,\ldots,\widetilde{R}^{2t-1},\widetilde{R}^0)$ in the graph~$\widetilde{\boldsymbol{H}}(t,2)$, where
\begin{equation*}
\begin{split}
\widetilde{R}^0:\!&=(0,\ldots,0)\; ,\\
\widetilde{R}^s:\!&=\sum\nolimits_{e\in[s]}\boldsymbol{\sigma}(e)\; ,\ \ \ 1\leq s\leq t-1\; ,
\end{split}
\end{equation*}
and
\begin{equation*}
\widetilde{R}^{t+k}:=\mathrm{T}^{(+)}-\widetilde{R}^k\; ,\ \ \ 0\leq k\leq t-1\; .
\end{equation*}
For any vertex~$\widetilde{T}$ of the hypercube graph~$\widetilde{\boldsymbol{H}}(t,2)$, let us define a subset~$\widetilde{\boldsymbol{Q}}(\widetilde{T},\widetilde{\boldsymbol{R}})\subset\mathrm{V}(\widetilde{\boldsymbol{R}}):=
(\widetilde{R}^0,\widetilde{R}^1,\ldots,\widetilde{R}^{2t-1})$ indirectly, via the mapping
\begin{equation*}
\widetilde{T}\; \overset{(\ref{eq:9})}{\mapsto}\; T\; ,
\end{equation*}
and via the bijection
\begin{equation*}
\boldsymbol{Q}(T,\boldsymbol{R})\; \xrightarrow{(\ref{eq:8})}\;
\widetilde{\boldsymbol{Q}}(\widetilde{T},\widetilde{\boldsymbol{R}})\; .
\end{equation*}
Involving the quantity~$\mathfrak{q}(\widetilde{T}):=\mathfrak{q}(\widetilde{T},\widetilde{\boldsymbol{R}})
:=|\widetilde{\boldsymbol{Q}}(\widetilde{T},\widetilde{\boldsymbol{R}})|=\mathfrak{q}(T)$,
we can write down the decomposition
\begin{equation}
\label{eq:69}
\widetilde{T}=-\tfrac{1}{2}(\mathfrak{q}(\widetilde{T})-1)\cdot\mathrm{T}^{(+)}
+\sum_{\substack{\widetilde{Q}\in\widetilde{\boldsymbol{Q}}(\widetilde{T},\widetilde{\boldsymbol{R}})\colon\\
\widetilde{Q}\neq(0,\ldots,0)=:\widetilde{R}^0}}\widetilde{Q}\; ,
\end{equation}
that describes yet another mechanism of {\em majority voting},\footnote{$\quad$Let $\mathbf{2}^{[t]}$ denote the {\em power set\/} (i.e., the family of all subsets) of $E_t$. Recall that maps
\begin{align*}
f\colon\{0,1\}^t\;&\to\;\mathbb{R}\; , & -\tfrac{1}{2}(\mathfrak{q}(\widetilde{T})-1)\cdot\mathrm{T}^{(+)}
+\sum_{\substack{\widetilde{Q}\in\widetilde{\boldsymbol{Q}}(\widetilde{T},\widetilde{\boldsymbol{R}})\colon\\
\widetilde{Q}\neq(0,\ldots,0)=:\widetilde{R}^0}}\widetilde{Q}\;=\;\widetilde{T}\,\;&\mapsto\,\; f(\widetilde{T})\; ,\\
\intertext{and}
g\colon\{1,-1\}^t\;&\to\;\mathbb{R}\; , & \sum_{Q\in\boldsymbol{Q}(T,\boldsymbol{R})}Q\;=\;T\,\;&\mapsto\,\; g(T)\; ,
\end{align*}
are {\em pseudo-Boolean functions}. Maps $\mathfrak{f}\colon\mathbf{2}^{[t]}\to\mathbb{R}$ are {\em set functions}. Maps $f\colon\{0,1\}^t\to\{0,1\}$, and $g\colon\{1,-1\}^t\to\{1,-1\}$ are {\em Boolean functions}.
}
but this decomposition has no essential meaning from the linear algebraic viewpoint, since the set~$\widetilde{Q}\in\widetilde{\boldsymbol{Q}}(\widetilde{T},\widetilde{\boldsymbol{R}})$ can contain the origin~$(0,\ldots,0)=:\widetilde{R}^0$
of the space~$\mathbb{R}^t$, which should be omitted in calculations.

For the topes $T\in\{1,-1\}^t$ of the oriented matroid $\mathcal{H}$, we define topes $\ro(T)\in\{1,-1\}^t$ by\footnote{
$\quad$$\ro(T)$ means the {\em relabeled opposite\/} of~$T$.}
\begin{equation}
\label{eq:6}
\ro(T):=-T\,\overline{\mathbf{U}}(t)\; ,
\end{equation}
where $\overline{\mathbf{U}}(t)$ denotes\footnote{$\quad$In~\cite[Sect.~2.1]{M-PROM} the similar notation $\mathbf{U}(m)$ was used to denote the backward identity matrix of order $(m+1)$ whose rows and columns were indexed starting with zero.} the {\em backward identity matrix\/} (with the rows and columns indexed starting with $1$) of order $t$ whose $(i,j)$th
entry is the Kronecker delta~$\delta_{i+j,t+1}$.

For vertices $\widetilde{T}$ of the discrete hypercube $\{0,1\}^t$, the counterparts of topes $\ro(T)$ of the oriented matroid~$\mathcal{H}$ are  vertices $\rn(\widetilde{T})\in\{0,1\}^t$, defined by\footnote{$\quad$$\rn(\widetilde{T})$ means the {\em relabeled negation\/} of~$\widetilde{T}$.}
\begin{equation}
\label{eq:35}
\rn(\widetilde{T}):=\mathrm{T}^{(+)}-\widetilde{T}\,\overline{\mathbf{U}}(t)\; .
\end{equation}
For example, suppose{\small
\begin{align*}
T:\!&=(1,-1,\phantom{-}1, -1,-1)\in\{1,-1\}^5\; ,\\
\widetilde{T}:\!&=(0,\phantom{-}1,\phantom{-}0, \phantom{-}1,\phantom{-}1)\in\{0,1\}^5\; .\\
\intertext{\normalsize{Then we have}}
\ro(T)&=(1,\phantom{-}1, -1,\phantom{-}1,-1)\; ,\\
\rn(\widetilde{T})&=(0,\phantom{-}0,\phantom{-}1,\phantom{-}0,\phantom{-}1)\; .
\end{align*}
}

\noindent$\bullet$ In the first part of the paper we compare the decompositions~$\boldsymbol{Q}(T,\boldsymbol{R})$ and~$\boldsymbol{Q}(\ro(T),\boldsymbol{R})$ of topes $T$ and $\ro(T)$ with respect to the symmetric cycle~$\boldsymbol{R}$ in the graph~$\boldsymbol{H}(t,2)$, defined by~(\ref{eq:12})(\ref{eq:13}).

Our interest in considering {\em relabeled opposites\/}~$\ro(T)$ and
{\em relabeled negations\/}~$\rn(\widetilde{T})$ lies in their application to combined {\em blocking/voting-models\/} of increasing families of sets and of clutters. We study those impractical $2^t$-dimensional vector models in order to gain a better understanding of the structure of families.

Recall that a family $\mathcal{A}:=\{A_1,\ldots,A_{\alpha}\}\subset\mathbf{2}^{[t]}$ of subsets\footnote{$\quad$We denote by $\hat{0}$ the {\em empty subset\/} of the ground set~$E_t$, and we let $\emptyset$ denote the {\em empty family} containing no sets.

Given a family $\mathcal{F}\subseteq\mathbf{2}^{[t]}$, such that $\emptyset\neq\mathcal{F}\not\ni\hat{0}$, the set $E_t:=[t]$ is the {\em ground set\/} of $\mathcal{F}$, while the set $\mathrm{V}(\mathcal{F}):=\bigcup_{F\in\mathcal{F}}F\subseteq E_t$ is
the {\em vertex set\/} of $\mathcal{F}$.

The families $\emptyset$ and $\{\hat{0}\}$ are the two {\em trivial clutters\/} on the ground set~$E_t$. The other clutters on $E_t$ are {\em nontrivial}.
} of the ground set~$E_t$ is called a {\em clutter}\footnote{$\quad$Or {\em Sperner family}, {\em antichain}, {\em simple hypergraph}.} if {\em no set\/} $A_i$ from $\mathcal{A}$ contains another set~$A_j$.

Given a family $\mathcal{F}\subseteq\mathbf{2}^{[t]}$, we let $\bmin\mathcal{F}$ denote the clutter composed of the {\em inclusion-minimal\/} sets in $\mathcal{F}$.

We say that a family of subsets $\mathcal{F}\subseteq\mathbf{2}^{[t]}$
is an {\em increasing family}\footnote{$\quad$Or {\em up-set}, {\em upward-closed family of sets}, {\em filter of sets}.} if the following implications hold:
\begin{equation*}
A\in\mathcal{F},\ \ \mathbf{2}^{[t]}\ni B\supset A\ \ \ \Longrightarrow\ \ \ B\in\mathcal{F}\; .
\end{equation*}

If $C\subseteq E_t$, then the family $\{C\}^{\triangledown}:=\{D\subseteq E_t\colon D\supseteq C\}$ is called the {\em principal\/} increasing family generated by the {\em one-member\/} clutter $\{C\}$. Conversely, an increasing family~$\mathcal{F}\subseteq\mathbf{2}^{[t]}$ is said to be {\em principal\/} if\footnote{$\quad$We denote by $|\!\cdot\!|$ the cardinality of a set, and we denote by $\#\,\cdot$ the number of sets in a family.} $\#\bmin\mathcal{F}=1$.

Given an arbitrary nonempty family $\mathcal{C}\subseteq\mathbf{2}^{[t]}$, we
denote by~$\mathcal{C}^{\triangledown}$ the {\em increasing family\/} on $E_t$, generated by~$\mathcal{C}$:
\begin{equation*}
\mathcal{C}^{\triangledown}:=\bigcup\nolimits_{C\in\mathcal{C}}\{C\}^{\triangledown}
=\bigcup\nolimits_{C\in\bmin\mathcal{C}}\{C\}^{\triangledown}\; .
\end{equation*}
``Decreasing'' constructs are defined in the obvious similar way.\footnote{
$\quad$For a family $\mathcal{F}\subseteq\mathbf{2}^{[t]}$, we use the notation $\bmax\mathcal{F}$ to denote the clutter composed of the {\em inclusion-maximal\/} sets in $\mathcal{F}$.

A family $\mathcal{F}\subseteq\mathbf{2}^{[t]}$ is said to be
a {\em decreasing family\/} (or
{\em down-set}, {\em downward-closed family of sets}, {\em ideal of sets}) if the following implications hold:
\begin{equation*}
B\in\mathcal{F},\ \ A\subset B\ \ \ \Longrightarrow\ \ \ A\in\mathcal{F}\; .
\end{equation*}
If $\emptyset\neq\mathcal{F}\neq\{\hat{0}\}$, then this decreasing family is the {\em abstract simplicial complex\/} on its {\em vertex set\/} $\bigcup_{M\in\bmax\mathcal{F}}M$, with the {\em facet\/} family~$\bmax\mathcal{F}$.

If $D\subseteq E_t$, then the family $\{D\}^{\vartriangle}:=\{C\colon C\subseteq D\}$ is called the {\em principal\/} decreasing family generated by the {\em one-member\/} clutter $\{D\}$. Conversely, a decreasing family~$\mathcal{F}\subseteq\mathbf{2}^{[t]}$ is said to be {\em principal\/} if $\#\bmax\mathcal{F}=1$.

Given an arbitrary nonempty family $\mathcal{D}\subseteq\mathbf{2}^{[t]}$, we
denote by~$\mathcal{D}^{\vartriangle}$ the {\em decreasing family\/} on $E_t$, generated by~$\mathcal{D}$:
\begin{equation*}
\mathcal{D}^{\vartriangle}:=\bigcup\nolimits_{D\in\mathcal{D}}\{D\}^{\vartriangle}=\bigcup\nolimits_{D\in\bmax\mathcal{D}}\{D\}^{\vartriangle}\; .
\end{equation*}
}

The duality philosophy behind clutters and increasing families is that any clutter is the {\em blocker}\footnote{$\quad$
\vspace{-5.5mm}
\begin{quote}
{\sl I enjoyed working with Ray} [{\sl{}Fulkerson}] {\sl and I coined the terms ``clutter'' and~``blocker''.}
\hfill Jack~Edmonds~\cite[p.~201]{50Y}
\end{quote}
} of a unique clutter, and any increasing family is the family of {\em blocking sets\/} of a unique clutter.

We often meet in the literature the {\em free distributive lattice\/} of antichains in the Boolean lattice of subsets of a finite nonempty set, ordered by containment of the corresponding generated {\em order ideals}, but an intrinsically related construct, the {\em free distributive lattice\/} of those antichains ordered by containment of the corresponding generated {\em order filters\/} has greater discrete mathematical expressiveness, because the latter lattice can be interpreted as the {\em lattice\/} of {\em blockers}, for which the {\em blocker map\/} is its anti-automorphism.\footnote{$\quad$For this paper, we chose the language of {\em power sets}, {\em clutters}, and {\em increasing\/} and {\em decreasing families}. A parallel exposition could be presented in poset-theoretic terms of {\em Boolean lattices}, {\em antichains}, and {\em order filters\/} and {\em ideals}.}

Recall that a subset $B\subseteq E_t$ is called a {\em blocking set\footnote{$\quad$Or
{\em transversal}, {\em hitting set}, {\em vertex cover} (or {\em node cover}), {\em system of representatives}.}\/} of a subset family~\mbox{$\mathcal{F}\subset\mathbf{2}^{[t]}$,} where~$\emptyset\neq\mathcal{F}\not\ni\hat{0}$, if we have
\begin{equation*}
|B\cap F|>0\; ,
\end{equation*}
for each set $F\in\mathcal{F}$. The {\em blocker}\footnote{$\quad$Or {\em blocking hypergraph} (or {\em transversal hypergraph}), {\em blocking clutter}, {\em dual clutter}, {\em Alexander dual clutter}.} $\mathfrak{B}(\mathcal{F})$ of the family~$\mathcal{F}$ is
the family of all {\em inclusion-minimal blocking sets\/} of  $\mathcal{F}$; note that we have
$\mathfrak{B}(\mathcal{F})=\mathfrak{B}(\bmin\mathcal{F})$. The notation~$\mathfrak{B}(\mathcal{F})^{\triangledown}$ just means the increasing family of all blocking sets of the family~$\mathcal{F}$.

For a nonempty family of subsets $\mathcal{F}\subseteq\mathbf{2}^{[t]}$,
we define a family\footnote{$\quad$Given a nonempty family of subsets $\mathcal{F}\subseteq\mathbf{2}^{[t]}$,
we define a family of complements~$\mathcal{F}^{\bot}$ by $\mathcal{F}^{\bot}:=\{F^{\bot}\colon F\in\mathcal{F}\}$, where $F^{\bot}:=
\mathrm{V}(\mathcal{F})-F$.
}
of complements $\mathcal{F}^{\compl}$ by $\mathcal{F}^{\compl}:=\{F^{\compl}\colon F\in\mathcal{F}\}$, where $F^{\compl}:=E_t-F$.

Given a nontrivial clutter $\mathcal{A}\subset\mathbf{2}^{[t]}$, one associates with $\mathcal{A}$ the four
extensively studied partitions of the power set of the ground set $E_t$:
\begin{align}
\label{eq:22}
\mathbf{2}^{[t]}&=\mathcal{A}^{\triangledown}\ \dot{\cup}\ (\mathfrak{B}(\mathcal{A})^{\compl})^{\vartriangle}\; ,\\
\nonumber
\mathbf{2}^{[t]}&=\mathcal{A}^{\vartriangle}\ \dot{\cup}\ \mathfrak{B}(\mathcal{A}^{\compl})^{\triangledown}\; ,\\
\nonumber
\mathbf{2}^{[t]}&=\mathfrak{B}(\mathcal{A})^{\triangledown}\ \dot{\cup}\ (\mathcal{A}^{\compl})^{\vartriangle}\; ,\\
\intertext{and}
\nonumber
\mathbf{2}^{[t]}&=\mathfrak{B}(\mathcal{A})^{\vartriangle}\ \dot{\cup}\ \mathfrak{B}(\mathfrak{B}(\mathcal{A})^{\compl})^{\triangledown}\; .
\end{align}

\noindent$\bullet$ In the second part of the paper we arrange the subsets of the ground set~$E_t$ in linear order. We then turn to the so-called {\em characteristic vectors\/} $\boldsymbol{\gamma}(\mathcal{F})\in\{0,1\}^{2^t}$ of subset families~$\mathcal{F}\subseteq\mathbf{2}^{[t]}$. If~$\mathcal{A}\subset\mathbf{2}^{[t]}$ is a nontrivial clutter on~$E_t$, then relation~(\ref{eq:22}) reformulated in the form (cf.~(\ref{eq:35}))
\begin{equation*}
\boldsymbol{\gamma}(\mathfrak{B}(\mathcal{A})^{\triangledown})
=\mathrm{T}_{2^t}^{(+)}-\boldsymbol{\gamma}(\mathcal{A}^{\triangledown})\cdot\overline{\mathbf{U}}(2^t)\; ,
\end{equation*}
where $\mathrm{T}_{2^t}^{(+)}:=(1,\ldots,1)$ is the $2^t$-dimensional row vector of all $1$'s, provides us with the characteristic vector
\begin{equation*}
\boldsymbol{\gamma}(\mathfrak{B}(\mathcal{A})^{\triangledown})=\rn(\boldsymbol{\gamma}(\mathcal{A}^{\triangledown}))
\end{equation*}
of the increasing family of blocking sets~$\mathfrak{B}(\mathcal{A})^{\triangledown}$ of the clutter~$\mathcal{A}$.

\noindent$\bullet$ In the third part of the paper we mention a blocking/voting-connection of the characteristic vectors $\boldsymbol{\gamma}(\mathcal{A}^{\triangledown})$ and $\boldsymbol{\gamma}(\mathfrak{B}(\mathcal{A})^{\triangledown})$
with the decompositions of the corresponding characteristic topes of the increasing families~$\mathcal{A}^{\triangledown}$ and~$\mathfrak{B}(\mathcal{A})^{\triangledown}$ with respect to a distinguished symmetric cycle in the hypercube graph~$\boldsymbol{H}(2^t,2)$, which is analogous to the cycle~(\ref{eq:12})(\ref{eq:13}) in the graph~$\boldsymbol{H}(t,2)$.

\part*{Decomposing}

\section{Topes, their relabeled opposites, and decompositions}

In this section we consider vertices $T$ of the discrete hypercube~$\{1,-1\}^t$, their relabeled opposites~$\ro(T)$ defined by~(\ref{eq:6}), and we discuss basic properties of the decompositions $\boldsymbol{Q}(T,\boldsymbol{R})$ and $\boldsymbol{Q}(\ro(T),\boldsymbol{R})$ of the topes~$T$ and~$\ro(T)$ with respect to the distinguished symmetric cycle~$\boldsymbol{R}:=(R^0,R^1,\ldots,R^{2t-1},R^0)$ in the graph~$\boldsymbol{H}(t,2)$, defined by~(\ref{eq:12})(\ref{eq:13}).

\noindent$\bullet$ Definitions~(\ref{eq:6}) and~(\ref{eq:35}) determine the maps
\begin{align}
\label{eq:10}
\{1,-1\}^t &\to \{1,-1\}^t: & T &\mapsto \ro(T):=\phantom{\mathrm{T}^{(+)}}-T\,\overline{\mathbf{U}}(t)\; ,\\
\label{eq:11}
\{0,1\}^t &\to \{0,1\}^t: & \widetilde{T} &\mapsto \rn(\widetilde{T})
:=\mathrm{T}^{(+)}-\widetilde{T}\,\overline{\mathbf{U}}(t)\; ,
\end{align}
and since we deal with the standard one-to-one correspondences between the vertex sets of the discrete hypercubes $\{1,-1\}^t$ and $\{0,1\}^t$, established by means of the maps~(\ref{eq:8}) and~(\ref{eq:9}), we mention the mappings
\begin{align}
\nonumber
 \{1,-1\}^t\ni \ro(T)\ &\overset{(\ref{eq:8})}{\mapsto}\  \rn(\widetilde{T})=\tfrac{1}{2}(\mathrm{T}^{(+)}+T\,\overline{\mathbf{U}}(t))\in\{0,1\}^t\; ,\\
\intertext{and}
\label{eq:36}
 \{0,1\}^t\ni\rn(\widetilde{T})\  &\overset{(\ref{eq:9})}{\mapsto}\    \ro(T)=-\mathrm{T}^{(+)}+2\widetilde{T}\,\overline{\mathbf{U}}(t)\in\{1,-1\}^t\; .
\end{align}

\noindent$\bullet$ Of course, the maps~(\ref{eq:10}) and~(\ref{eq:11}) are both {\em involutions}:
\begin{equation*}
\{1,-1\}^t\ni\ro(\ro(T))=T\; ,\ \ \ \text{and}\ \ \ \ \{0,1\}^t\ni\rn(\rn(\widetilde{T}))=\widetilde{T}\; .
\end{equation*}

\noindent$\bullet$ Given a vector $\boldsymbol{z}:=(z_1,\ldots,z_t)\in\mathbb{R}^t$, we denote its {\em support} $\{e\in E_t\colon z_e\neq 0\}$
by~$\supp(\boldsymbol{z})$. For a vertex~$\widetilde{T}$ of the discrete hypercube~$\{0,1\}^t$, we let~$\hwt(\widetilde{T})$ denote its {\em Hamming weight}: $\hwt(\widetilde{T}):=|\supp(\widetilde{T})|$.

Note that we have
\begin{align*}
\{1,-1\}^t\ni \ro(T)=T\ \ \ &\Longleftrightarrow\ \ \ -T=T\,\overline{\mathbf{U}}(t)\; ;\\
\ro(T)=T\ \ \ &\,\Longrightarrow\ \ \ |T^-|=\tfrac{t}{2}\; .\\
\intertext{We also have}
\{0,1\}^t\ni\rn(\widetilde{T})=\widetilde{T}\ \ \ &\Longleftrightarrow\ \ \ \mathrm{T}^{(+)}-\widetilde{T}=\widetilde{T}\,\overline{\mathbf{U}}(t)\; ;\\
\rn(\widetilde{T})=\widetilde{T}\ \ \ &\,\Longrightarrow\ \ \
\hwt(\widetilde{T})
=\tfrac{t}{2}\; .
\end{align*}
Thus, if $t$ is {\em odd}, then we always have $\ro(T)\neq T$, and $\rn(\widetilde{T})\neq\widetilde{T}$.
\newpage
\noindent$\bullet$ Let $\langle\cdot,\cdot\rangle$ denote the standard scalar product on the space~$\mathbb{R}^t$.

For vertices~$\widetilde{T}\in\{0,1\}^t$ and $T:={}_{-\supp(\widetilde{T})}\mathrm{T}^{(+)}\in\{1,-1\}^t$,
we have
\begin{multline*}
\langle T,\ro(T)\rangle=-T\,\overline{\mathbf{U}}(t)T^{\top}
=-\sum_{e\in[t]}T(e)T(t-e+1)
\\
=\begin{cases}
-1-2\sum_{e\in[(t-1)/2]}T(e)T(t-e+1)\; , & \text{if $t$ is odd},\\
\quad\vspace{-3mm}\\
\phantom{-1}-2\sum_{e\in[t/2]}T(e)T(t-e+1)\; , & \text{if $t$ is even};
\end{cases}
\end{multline*}
\begin{multline*}
\langle\widetilde{T},\rn(\widetilde{T})\rangle:=\langle\widetilde{T},\mathrm{T}^{(+)}-\widetilde{T}\,\overline{\mathbf{U}}(t)\rangle
=\hwt(\widetilde{T})-\sum_{e\in[t]}\widetilde{T}(e)\widetilde{T}(t-e+1)\\
=\hwt(\widetilde{T})-\begin{cases}
\widetilde{T}((t+1)/2)+2\sum_{e\in[(t-1)/2]}\widetilde{T}(e)\widetilde{T}(t-e+1)\; , & \text{if $t$ is odd},\\
\quad\vspace{-3mm}\\
\phantom{T((t+1)/2)+}\ \, 2\sum_{e\in[t/2]}\widetilde{T}(e)\widetilde{T}(t-e+1)\; , & \text{if $t$ is even}.
\end{cases}
\end{multline*}

\noindent$\bullet$ Given two {\em words\/} $X,Y\in\{-1,0,1\}^t$, we let $d(X,Y):=|\{e\in E_t\colon X(e)\neq Y(e)\}|$ denote the {\em Hamming distance\/} between them.\footnote{$\quad$If $X$ and $Y$ are {\em topes}, then one says that $d(X,Y)$ is the {\em graph distance}.}

Since the equal distances $d(T,\ro(T))=d(\widetilde{T},\rn(\widetilde{T}))$ can be calculated with the help of the formulas (see~(\ref{eq:9}) and~(\ref{eq:36}))
\begin{align*}
d(T,\ro(T))&=\tfrac{1}{2}\bigl(t-\langle T,\ro(T)\rangle\bigr)\; ,\\
d(\widetilde{T},\rn(\widetilde{T}))&=\tfrac{1}{2}\bigl(t-\langle \mathrm{T}^{(+)}-2\widetilde{T},-\mathrm{T}^{(+)}+2\widetilde{T}\,\overline{\mathbf{U}}(t)\rangle\bigr)\; ,
\end{align*}
we see that
\begin{align*}
d(T,\ro(T))&=\tfrac{1}{2}\bigl(t+T\,\overline{\mathbf{U}}(t)T^{\top}\bigr)
=\tfrac{1}{2}\bigl(t+\sum_{e\in[t]}T(e)T(t-e+1)\bigr)\\
&=\frac{t}{2}+
\begin{cases}
\frac{1}{2}+\sum_{e\in[(t-1)/2]}T(e)T(t-e+1)\; , & \text{if $t$ is odd},\\
\quad\vspace{-3mm}\\
\phantom{\frac{1}{2}+}\ \sum_{e\in[t/2]}T(e)T(t-e+1)\; , & \text{if $t$ is even},
\end{cases}
\end{align*}
and
\begin{multline*}
d(\widetilde{T},\rn(\widetilde{T}))
=t-2\cdot\hwt(\widetilde{T})+2\sum_{e\in[t]}\widetilde{T}(e)\widetilde{T}(t-e+1)\\
=t-2\cdot\hwt(\widetilde{T})
\\+2\cdot\begin{cases}
\widetilde{T}((t+1)/2)+2\sum_{e\in[(t-1)/2]}\widetilde{T}(e)\widetilde{T}(t-e+1)\; , & \text{if $t$ is odd},\\
\quad\vspace{-3mm}\\
\phantom{T((t+1)/2)+}\ \, 2\sum_{e\in[t/2]}\widetilde{T}(e)\widetilde{T}(t-e+1)\; , & \text{if $t$ is even}.
\end{cases}
\end{multline*}

\noindent$\bullet$ Suppose that $4|t$ (i.e., $t$ is divisible by $4$). Note that
\begin{align*}
\langle T,\ro(T)\rangle=0\ \ \ &\Longleftrightarrow\ \ \
\sum_{e\in[t/2]}T(e)T(t-e+1)=0\; ;\\
\langle T,\ro(T)\rangle=0\ \ \ &\Longleftrightarrow\ \ \
\sum_{e\in[t/2]}\widetilde{T}(e)\widetilde{T}(t-e+1)=\frac{4\cdot\hwt(\widetilde{T})-t}{8}\; .
\end{align*}

\noindent$\bullet$ Considering the restriction of the map~(\ref{eq:10})
to the vertex set $\mathrm{V}(\boldsymbol{R})$ of the symmetric cycle $\boldsymbol{R}$ in the hypercube graph~$\boldsymbol{H}(t,2)$, defined by~(\ref{eq:12})(\ref{eq:13}), we have the mappings
\begin{equation*}
R^i\ \overset{(\ref{eq:10})}{\mapsto}\ \ro(R^i)=R^{(3t-i)\hspace{-1.5mm}\mod{2t}}
=\begin{cases}
R^{t-i}, & \text{if $0\leq i\leq t$}\; ,\\
\quad\vspace{-3mm}\\
R^{3t-i}, & \text{if $t+1\leq i\leq 2t-1$}\; .
\end{cases}
\end{equation*}
If $t$ is {\em even}, then the following implication holds:
\begin{equation*}
R^i\in\mathrm{V}(\boldsymbol{R}),\ \ \ro(R^i)=R^i\ \ \ \Longrightarrow\ \ \ i\in\{\tfrac{t}{2},\tfrac{3t}{2}\}\; .
\end{equation*}

\begin{remark} Let $\boldsymbol{R}$ be the symmetric cycle in the hypercube graph~$\boldsymbol{H}(t,2)$, defined by~{\rm(\ref{eq:12})(\ref{eq:13})}.
Given a vertex $T\in\{1,-1\}^t$ of $\boldsymbol{H}(t,2)$,
suppose that
\begin{equation*}
(R^0,R^1,\ldots,R^{2t-1})=:\mathrm{V}(\boldsymbol{R})\supset\boldsymbol{Q}(T,\boldsymbol{R})=(R^{i_0},R^{i_1},\ldots,R^{i_{\mathfrak{q}(T)-1}})\; ,
\end{equation*}
for some indices $i_0<i_1<\cdots<\mathfrak{q}(T)-1$.
\begin{itemize}
\item[\rm(i)]
We have
\begin{equation*}
\boldsymbol{Q}(\ro(T),\boldsymbol{R})=(R^{(3t-i_0)\hspace{-1.5mm}\mod{2t}},R^{(3t-i_1)\hspace{-1.5mm}\mod{2t}},\ldots,
R^{(3t-i_{\mathfrak{q}(T)-1})\hspace{-1.5mm}\mod{2t}})\; ,
\end{equation*}
or, in other words,
\begin{equation*}
\boldsymbol{Q}(\ro(T),\boldsymbol{R})=\{R^{(3t-i)\hspace{-1.5mm}\mod{2t}}\colon R^i\in\boldsymbol{Q}(T,\boldsymbol{R})\}\; .
\end{equation*}
\item[\rm(ii)] If $t$ is {\em even}, then
\begin{equation*}
\ro(T)=T\ \ \ \Longleftrightarrow\ \ \
\Bigl(Q\in\boldsymbol{Q}(T,\boldsymbol{R}) \Longrightarrow \ro(Q)\in\boldsymbol{Q}(T,\boldsymbol{R})\Bigr)\; .
\end{equation*}
Note that the following implication holds:
\begin{equation*}
\ro(T)=T\ \ \ \Longrightarrow\ \ \
|\{R^{t/2},R^{3t/2}\}\cap\boldsymbol{Q}(T,\boldsymbol{R})|=1\; .
\end{equation*}
\end{itemize}
\end{remark}

\noindent$\bullet$ Recall that for any vertex $T\in\{1,-1\}^t$ of the hypercube graph~$\boldsymbol{H}(t,2)$ with its distinguished symmetric cycle~$\boldsymbol{R}$ defined by~{\rm(\ref{eq:12})(\ref{eq:13})}, there exists a unique row vector $\boldsymbol{x}:=\boldsymbol{x}(T):=\boldsymbol{x}(T,\boldsymbol{R}):=(x_1,\ldots,x_t)\in\{-1,0,1\}^t$ such that
\begin{equation*}
T=\sum_{i\in[t]}x_i\cdot R^{i-1}=\boldsymbol{x}\mathbf{M}\; ,
\end{equation*}
where
\begin{equation*}
\mathbf{M}:=\mathbf{M}(\boldsymbol{R}):=\left(\begin{smallmatrix}
R^0\\
R^1\vspace{-2mm}\\
\vdots\\
R^{t-1}
\end{smallmatrix}
\right)\; .
\end{equation*}
In other words, the inclusion-minimal linearly independent set $\boldsymbol{Q}(T,\boldsymbol{R})$ of odd cardinality, given in~(\ref{eq:37})(\ref{eq:32}), is described as
\begin{equation*}
\boldsymbol{Q}(T,\boldsymbol{R})=\{x_i\cdot R^{i-1}\colon x_i\neq 0\}\; .
\end{equation*}
Recall that if $x_e\neq 0$ for some $e\in E_t$, then $x_e=T(e)$.

We will now give an explicit description of decompositions $\boldsymbol{Q}(T,\boldsymbol{R})$ and $\boldsymbol{Q}(\ro(T),\boldsymbol{R})$ via the corresponding ``$\boldsymbol{x}$-vectors''.

Let $\overline{\mathbf{T}}(t)$ denote\footnote{$\quad$
In~\cite[Sect.~2.1]{M-PROM} the similar notation $\mathbf{T}(m)$ was used to denote the forward shift matrix of order $(m+1)$ whose rows and columns were indexed starting with zero.
}
the {\em forward shift matrix\/} of order $t$ whose $(i,j)$th entry is $\delta_{j-i,1}$.

\begin{proposition}\cite[Prop.~2.4, extended]{M-SC-II} \label{th:10}
Let $\boldsymbol{R}$ be the symmetric cycle in the hypercube graph~$\boldsymbol{H}(t,2)$, defined by~{\rm(\ref{eq:12})(\ref{eq:13})}.

Let $A$ be a nonempty subset of the ground set $E_t$, regarded as a disjoint union
\begin{equation*}
A=[i_1,j_1]\;\dot\cup\;[i_2,j_2]\;\dot\cup\;\cdots\;\dot\cup\;[i_{\varrho-1},j_{\varrho-1}]\;\dot\cup\;[i_{\varrho},j_{\varrho}]
\end{equation*}
of intervals such that
\begin{equation*}
j_1+2\leq i_2,\ \ j_2+2\leq i_3,\ \ \ldots,\ \
j_{\varrho-2}+2\leq i_{\varrho-1},\ \ j_{\varrho-1}+2\leq i_{\varrho}\; ,
\end{equation*}
for some $\varrho:=\varrho(A)$.

\begin{itemize}
\item[\rm(i)]
{\rm(a)} If $\{1,t\}\cap A=\{\underset{\overset{\uparrow}{i_1}}{1}\}$, then
we have
\begin{align*}
|\boldsymbol{Q}({}_{-A}\mathrm{T}^{(+)},\boldsymbol{R})|&=
2\varrho-1\; ,\\
\boldsymbol{x}({}_{-A}\mathrm{T}^{(+)},\boldsymbol{R})&=
\sum_{1\leq k\leq\varrho}\boldsymbol{\sigma}(j_k+1)-\sum_{2\leq \ell\leq\varrho}\boldsymbol{\sigma}(i_{\ell})
\; .
\end{align*}
{\rm(b)} Since
\begin{multline*}
\{t-e+1\colon e\in E_t-A\}=[\underset{\overset{\uparrow}{i_1}}{1},t-j_{\varrho}]\;\dot\cup\;[t-i_{\varrho}+2,t-j_{\varrho-1}]\\
\dot\cup\;\cdots\;\dot\cup\;[t-i_3+2,t-j_2]
\;\dot\cup\;[t-i_2+2,t-j_1]\; ,
\end{multline*}
and $\{1,t\}\cap\{t-e+1\colon e\in E_t-A\}=\{1\}$, we see that
\begin{align*}
|\boldsymbol{Q}(\ro({}_{-A}\mathrm{T}^{(+)}),\boldsymbol{R})|&=
2\varrho-1\; ,\\
\boldsymbol{x}(\ro({}_{-A}\mathrm{T}^{(+)}),\boldsymbol{R})&=
\sum_{1\leq k\leq\varrho}\boldsymbol{\sigma}(t-j_k+1)-\sum_{2\leq \ell\leq\varrho}\boldsymbol{\sigma}(t-i_{\ell}+2)\; .
\end{align*}
{\rm(c)} Note that
\begin{equation*}
\boldsymbol{x}(\ro({}_{-A}\mathrm{T}^{(+)}),\boldsymbol{R})=
\boldsymbol{x}({}_{-A}\mathrm{T}^{(+)},\boldsymbol{R})\cdot\overline{\mathbf{U}}(t)\cdot
\overline{\mathbf{T}}(t)\; .
\end{equation*}

\item[\rm(ii)]
{\rm(a)} If $\{1,t\}\cap A=\{\underset{\overset{\uparrow}{i_1}}{1},\underset{\overset{\uparrow}{j_{\varrho}}}{t}\}$, then
\begin{align*}
|\boldsymbol{Q}({}_{-A}\mathrm{T}^{(+)},\boldsymbol{R})|&=2\varrho-1\; ,\\
\boldsymbol{x}({}_{-A}\mathrm{T}^{(+)},\boldsymbol{R})&=-\boldsymbol{\sigma}(1)+\sum_{1\leq k\leq\varrho-1}\boldsymbol{\sigma}(j_k+1)-
\sum_{2\leq\ell \leq\varrho}\boldsymbol{\sigma}(i_{\ell})
\; .
\end{align*}
{\rm(b)} Since
\begin{multline*}
\{t-e+1\colon e\in E_t-A\}=[t-i_{\varrho}+2,t-j_{\varrho-1}]\;\dot\cup\;[t-i_{\varrho-1}+2,t-j_{\varrho-2}]\\
\dot\cup\;\cdots\;\dot\cup\;[t-i_3+2,t-j_2]
\;\dot\cup\;[t-i_2+2,t-j_1]\; ,
\end{multline*}
and $|\{1,t\}\cap\{t-e+1\colon e\in E_t-A\}|=0$, we have
\begin{align*}
|\boldsymbol{Q}(\ro({}_{-A}\mathrm{T}^{(+)}),\boldsymbol{R})|&=
2\varrho-1\; ,\\
\boldsymbol{x}(\ro({}_{-A}\mathrm{T}^{(+)}),\boldsymbol{R})&=\boldsymbol{\sigma}(1)+
\sum_{1\leq k\leq\varrho-1}\boldsymbol{\sigma}(t-j_k+1)-\sum_{2\leq \ell\leq\varrho}\boldsymbol{\sigma}(t-i_{\ell}+2)\; .
\end{align*}
{\rm(c)} Note that
\begin{equation*}
\boldsymbol{x}(\ro({}_{-A}\mathrm{T}^{(+)}),\boldsymbol{R})=
\boldsymbol{\sigma}(1)+\boldsymbol{x}({}_{-A}\mathrm{T}^{(+)},\boldsymbol{R})\cdot\overline{\mathbf{U}}(t)\cdot
\overline{\mathbf{T}}(t)\; .
\end{equation*}

\item[\rm(iii)]
{\rm(a)} If $|\{1,t\}\cap A|=0$, then
\begin{align*}
|\boldsymbol{Q}({}_{-A}\mathrm{T}^{(+)},\boldsymbol{R})|&=2\varrho\;+1\; ,\\
\boldsymbol{x}({}_{-A}\mathrm{T}^{(+)},\boldsymbol{R})&=\boldsymbol{\sigma}(1)\;+\sum_{1\leq k\leq\varrho}\boldsymbol{\sigma}(j_k+1)\;-
\sum_{1\leq\ell \leq\varrho}\boldsymbol{\sigma}(i_{\ell})
\; .
\end{align*}
{\rm(b)} Since
\begin{multline*}
\{t-e+1\colon e\in E_t-A\}=[1,t-j_{\varrho}]\;\dot\cup\;[t-i_{\varrho}+2,t-j_{\varrho-1}]\\
\dot\cup\;\cdots\;\dot\cup\;[t-i_2+2,t-j_1]
\;\dot\cup\;[t-i_1+2,t]\; ,
\end{multline*}
and $\{1,t\}\cap\{t-e+1\colon e\in E_t-A\}=\{1,t\}$, we have
\begin{align*}
|\boldsymbol{Q}(\ro({}_{-A}\mathrm{T}^{(+)}),\boldsymbol{R})|&=
2\varrho+1\; ,\\
\boldsymbol{x}(\ro({}_{-A}\mathrm{T}^{(+)}),\boldsymbol{R})&=
-\boldsymbol{\sigma}(1)+\sum_{1\leq k\leq\varrho}\boldsymbol{\sigma}(t-j_k+1)-\sum_{1\leq \ell\leq\varrho}\boldsymbol{\sigma}(t-i_{\ell}+2)\; .
\end{align*}
{\rm(c)} Note that
\begin{equation*}
\boldsymbol{x}(-({}_{-A}\mathrm{T}^{(+)})\cdot\overline{\mathbf{U}}(t),\boldsymbol{R})=
-\boldsymbol{\sigma}(1)+
\boldsymbol{x}({}_{-A}\mathrm{T}^{(+)},\boldsymbol{R})\cdot\overline{\mathbf{U}}(t)\cdot
\overline{\mathbf{T}}(t)\; .
\end{equation*}

\item[\rm(iv)]
{\rm(a)} If $\{1,t\}\cap A=\{\underset{\overset{\uparrow}{j_{\varrho}}}{t}\}$, then
\begin{align*}
|\boldsymbol{Q}({}_{-A}\mathrm{T}^{(+)},\boldsymbol{R})|&=2\varrho-1\; ,\\
\boldsymbol{x}({}_{-A}\mathrm{T}^{(+)},\boldsymbol{R})&=
\sum_{1\leq k\leq\varrho-1}\boldsymbol{\sigma}(j_k+1)
-\sum_{1\leq \ell\leq\varrho}\boldsymbol{\sigma}(i_{\ell})
\; .
\end{align*}
{\rm(b)} Since
\begin{multline*}
\{t-e+1\colon e\in E_t-A\}=[t-i_{\varrho}+2,t-j_{\varrho-1}]\;\dot\cup\;[t-i_{\varrho-1}+2,t-j_{\varrho-2}]\\
\dot\cup\;\cdots\;\dot\cup\;[t-i_2+2,t-j_1]
\;\dot\cup\;[t-i_1+2,\underset{\overset{\uparrow}{j_{\varrho}}}{t}]\; ,
\end{multline*}
and $\{1,t\}\cap\{t-e+1\colon e\in E_t-A\}=\{t\}$, we see that
\begin{align*}
|\boldsymbol{Q}(\ro({}_{-A}\mathrm{T}^{(+)}),\boldsymbol{R})|&=
2\varrho-1\; ,\\
\boldsymbol{x}(\ro({}_{-A}\mathrm{T}^{(+)}),\boldsymbol{R})&=
\sum_{1\leq k\leq\varrho-1}\boldsymbol{\sigma}(t-j_k+1)-\sum_{1\leq \ell\leq\varrho}\boldsymbol{\sigma}(t-i_{\ell}+2)\; .
\end{align*}
{\rm(c)} Note that
\begin{equation*}
\boldsymbol{x}(\ro({}_{-A}\mathrm{T}^{(+)}),\boldsymbol{R})=
\boldsymbol{x}({}_{-A}\mathrm{T}^{(+)},\boldsymbol{R})\cdot\overline{\mathbf{U}}(t)\cdot
\overline{\mathbf{T}}(t)\; .
\end{equation*}
\end{itemize}
\end{proposition}

\begin{corollary}\label{th:11}
Let $\boldsymbol{R}$ be the symmetric cycle in the hypercube graph~$\boldsymbol{H}(t,2)$, defined by~{\rm(\ref{eq:12})(\ref{eq:13})}.

For any vertex $T\in\{1,-1\}^t$ of the graph~$\boldsymbol{H}(t,2)$ we have
\begin{equation*}
\mathfrak{q}(\ro(T)):=|\supp(\boldsymbol{x}(\ro(T),\boldsymbol{R}))|=
|\supp(\boldsymbol{x}(T,\boldsymbol{R}))|=:\mathfrak{q}(T)\; .
\end{equation*}

\begin{itemize}
\item[\rm(i)]
If\/ $|T^-\cap\{1,t\}|=1$, then
\begin{equation*}
\boldsymbol{x}(\ro(T))=
\boldsymbol{x}(T)\cdot\overline{\mathbf{U}}(t)\cdot
\overline{\mathbf{T}}(t)\; .
\end{equation*}

\item[\rm(ii)]
If\/ $|T^-\cap\{1,t\}|=2$, then
\begin{equation*}
\boldsymbol{x}(\ro(T))=
\boldsymbol{\sigma}(1)+\boldsymbol{x}(T)\cdot\overline{\mathbf{U}}(t)\cdot
\overline{\mathbf{T}}(t)\; .
\end{equation*}

\item[\rm(iii)]
If\/ $|T^-\cap\{1,t\}|=0$, then
\begin{equation*}
\boldsymbol{x}(\ro(T))=
-\boldsymbol{\sigma}(1)+
\boldsymbol{x}(T)\cdot\overline{\mathbf{U}}(t)\cdot
\overline{\mathbf{T}}(t)\; .
\end{equation*}
\end{itemize}
\end{corollary}

\part*{Blocking}

Blocking sets and the blockers of set families
(families are often regarded as the {\em hyperedge\/} families of {\em hypergraphs}) are discussed, e.g., in the monographs~\cite{Berge89,Bretto,CCZ,Co,CH,CFKLMPPS,FLSZ,G,GJ,GLS,HHH,HY-Total,HY,Jukna-E2,M-PROM,Mirsky,NI,NW,SchU,Sch-C,V}
and in the works~\cite{AN-2017,AH-2017,Alon-1990,AKB-2021,ABCS-2021,BS-2005,Barat-2021,vBS-2020,BH-2004,BFLMSCh-2019,BFSch-2021,BM-2009,Bou-2013,
BKMN-2016,BHT-2012,BRT-2021,ChHuZ-2013,ChMcD-1992,CCGK-2016,CFM-1991,D-2011,DNU-2020,EG-1995,EGT-2003,EMG-2008,Elbassioni-2008,FHNS-2020,
FKh-1996,F-1988,G-DV-L-2017,HBC-2007,KLMN-2014,Karp-1972,KS-2005,KhBEG-2006,LS-2008,LJ-2003,MU-2014,NP-Kao, PQ-2012,Reiter-1987,Y-2018}.

\noindent$\bullet$ Let
\begin{equation*}
\mho([t]):=\{\mathcal{A}\subset\mathbf{2}^{[t]}\colon\ \mathcal{A}\,=\,\bmin\mathcal{A}\,=\,\bmax\mathcal{A}\}
\end{equation*}
denote the family of clutters on the ground set~$E_t$.
The map
\begin{equation}
\label{eq:38}
\mho([t])\to\mho([t])\; ,\ \ \ \mathcal{A}\mapsto\mathfrak{B}(\mathcal{A})\; ,
\end{equation}
is called the {\em blocker map\/} on clutters~\cite{CFM-1991}.

\noindent$\bullet$ If the (\!{\em{}abstract simplicial}) {\em complex\/}~$\Delta:=(\mathfrak{B}(\mathcal{A})^{\compl})^{\vartriangle}$ appearing  in~(\ref{eq:27}) and the {\em complex\/}~$\Delta^{\!\!\vee}:=\{F^{\compl}\colon F\in\mathcal{A}^{\triangledown}\}$ both have the same vertex set~$E_t$, then the complex $\Delta^{\!\!\vee}$ is called the {\em Alexander dual\/} of the complex~$\Delta$; see, e.g.,~\cite{V} and~\cite{BT-2009} on combinatorial {\em Alexander duality}.

\noindent$\bullet$ Given a clutter~$\mathcal{A}$, the quantity
\begin{equation*}
\tau(\mathcal{A}):=\min\{|B|\colon B\in\mathfrak{B}(\mathcal{A})\}
\end{equation*}
is called the {\em transversal number}\footnote{$\quad$Or {\em covering number}, {\em vertex cover number}, {\em blocking number}.} of $\mathcal{A}$.

\noindent$\bullet$ Recall a classical result in combinatorial optimization: For any clutter $\mathcal{A}$
we have
\begin{equation*}
\mathfrak{B}(\mathfrak{B}(\mathcal{A}))=\mathcal{A}\; ,
\end{equation*}
see~\cite{EF-1970,I-1958,La-1966,Leh-164}.

\noindent$\bullet$ For a nontrivial clutter~$\mathcal{A}\subset\mathbf{2}^{[t]}$ on the ground set~$E_t$, we have
\begin{equation}
\label{eq:27}
\#\mathcal{A}^{\triangledown}\; +\; \#\mathfrak{B}(\mathcal{A})^{\triangledown}=2^t\; .
\end{equation}
More precisely, for any $s$, where $0\leq s\leq t$, we have
\begin{equation}
\label{eq:28}
\#(\mathfrak{B}(\mathcal{A})^{\triangledown}\cap\tbinom{E_t}{s})\; +\; \#(\mathcal{A}^{\triangledown}\cap\tbinom{E_t}{t-s})
=\tbinom{t}{s}\; ,
\end{equation}
where
\begin{equation*}
\tbinom{E_t}{s}:=\{F\subseteq E_t\colon |F|=s\}
\end{equation*}
is the {\em complete $s$-uniform clutter\/} on the vertex set $E_t$.

The increasing families~$\mathcal{A}^{\triangledown}$ and~$\mathfrak{B}(\mathcal{A})^{\triangledown}$ are {\em comparable by containment}: either we have
\begin{equation*}
\mathfrak{B}(\mathcal{A})^{\triangledown}\subseteq\mathcal{A}^{\triangledown}\; ,
\ \ \ \text{or\/}\ \ \ \
\mathfrak{B}(\mathcal{A})^{\triangledown}\supseteq\mathcal{A}^{\triangledown}\; .
\end{equation*}
The following implications hold:
\begin{gather*}
\mathfrak{B}(\mathcal{A})^{\triangledown}\subsetneqq \mathcal{A}^{\triangledown}\ \ \Longleftrightarrow\ \ \#\mathcal{A}^{\triangledown}>2^{t-1}\; ;
\\
\mathfrak{B}(\mathcal{A})^{\triangledown}\supsetneqq \mathcal{A}^{\triangledown}\ \ \Longleftrightarrow\ \ \#\mathcal{A}^{\triangledown}<2^{t-1}\; .
\end{gather*}
Note also that the following implications hold:
\begin{gather*}
\#\mathcal{A}^{\triangledown}>2^{t-1}\ \ \ \Longrightarrow\ \ \ \min\{|A|\colon A\in\mathcal{A}\}\leq\min\{|B|\colon B\in\mathfrak{B}(\mathcal{A})\}\; ;
\\
\#\mathcal{A}^{\triangledown}<2^{t-1}\ \ \ \Longrightarrow\ \ \ \min\{|A|\colon A\in\mathcal{A}\}\geq\min\{|B|\colon B\in\mathfrak{B}(\mathcal{A})\}\; .
\end{gather*}

\noindent$\bullet$ A clutter $\mathcal{A}$
is called {\em self-dual\/}~\cite[Ch.~9]{Jukna-E2}\cite[\S{}5.7]{M-PROM} or {\em identically self-blocking\/}~\cite{ACL-2019,AP-2018} if
\begin{equation*}
\mathfrak{B}(\mathcal{A})=\mathcal{A}\; ;
\end{equation*}
see also the early reference~\cite[\S{}2.1]{Berge89}.
In other words, the self-dual clutters~$\mathcal{A}\subset\mathbf{2}^{[t]}$ on the ground set~$E_t$ are the {\em fixed points\/} of the {\em blocker map\/}~(\ref{eq:38}); for each of them we also have
\begin{equation*}
\mathfrak{B}(\mathcal{A})^{\triangledown}=\mathcal{A}^{\triangledown}\; .
\end{equation*}

As noted in~\cite[Cor.~5.28(i)]{M-PROM}, one criterion for a {\em clutter\/} $\mathcal{A}\subset\mathbf{2}^{[t]}$ on the ground set~$E_t$ to be {\em self-dual\/} is as follows:
\begin{equation*}
\mathfrak{B}(\mathcal{A})=\mathcal{A}\ \ \Longleftrightarrow\ \ \#\mathcal{A}^{\triangledown}=2^{t-1}\; .
\end{equation*}

\noindent$\bullet$ Let~$X$ be a subset of the ground set~$E_t$. Given a nontrivial clutter~$\mathcal{A}$ on~$E_t$, its {\em deletion\/}
$\mathcal{A}\,\backslash\,X$ is defined to be the clutter
\begin{equation*}
\mathcal{A}\,\backslash\,X:=\{A\in\mathcal{A}\colon |A\cap X|=0\}\; .
\end{equation*}
The {\em contraction\/} $\mathcal{A}\,/X$ is defined to be the clutter
\begin{equation*}
\mathcal{A}\,/X:=\bmin\{A-X\colon A\in\mathcal{A}\}\; .
\end{equation*}

A classical result in combinatorial optimization is as follows:
\begin{equation*}
\mathfrak{B}(\mathcal{A})\,\backslash\,X=\mathfrak{B}(\mathcal{A}\,/X)\; ,\ \ \ \ \ \text{and}\ \ \ \ \ \
\mathfrak{B}(\mathcal{A})\,/X=\mathfrak{B}(\mathcal{A}\,\backslash\,X)\; ,
\end{equation*}
see~\cite{S-1975}.

We also have
\begin{equation*}
(\mathfrak{B}(\mathcal{A})\,\backslash\,X)^{\triangledown}=\mathfrak{B}(\mathcal{A}\,/X)^{\triangledown}\subseteq
\mathfrak{B}(\mathcal{A})^{\triangledown}
\subseteq
(\mathfrak{B}(\mathcal{A})\,/X)^{\triangledown}=\mathfrak{B}(\mathcal{A}\,\backslash\,X)^{\triangledown}\; ,
\end{equation*}
cf.~\cite[Eq.\,(5.4)]{M-PROM}. Further,
\begin{align*}
\#(\mathcal{A}\,\backslash\,X)^{\triangledown}\ +\ \#(\mathfrak{B}(\mathcal{A})\,/X)^{\triangledown}\ &=\ 2^t\; ,\\
\#(\mathcal{A}\,/X)^{\triangledown}\ +\ \#(\mathfrak{B}(\mathcal{A})\,\backslash\,X)^{\triangledown}\ &=\ 2^t\; ,
\end{align*}
see~\cite[Cor.~5.28(ii)]{M-PROM}. More precisely, for any $s$, where $0\leq s\leq t$, we have
\begin{align*}
\#\bigl((\mathfrak{B}(\mathcal{A})\,/X)^{\triangledown}\cap\tbinom{E_t}{s}\bigr)\ +\
\#\bigl((\mathcal{A}\,\backslash\,X)^{\triangledown}\cap\tbinom{E_t}{t-s}\bigr)\ &=\ \tbinom{t}{s}\; ,\\
\#\bigl((\mathfrak{B}(\mathcal{A})\,\backslash\,X)^{\triangledown}\cap\tbinom{E_t}{s}\bigr)\ +\
\#\bigl((\mathcal{A}\,/X)^{\triangledown}\cap\tbinom{E_t}{t-s}\bigr)\ &=\ \tbinom{t}{s}\; .
\end{align*}

\noindent$\bullet$ Let $p$ be a rational number such that $0\leq p<1$. Given a nontrivial clutter $\mathcal{A}:=\{A_1,\ldots,A_{\alpha}\}\subset\mathbf{2}^{[t]}$ on the ground set~$E_t$, a subset $B\subseteq E_t$ is called a {\em $p$-committee}\footnote{$\quad$By convention, a $\tfrac{1}{2}$-committee of a clutter $\mathcal{A}$ is called its {\em committee}.} of the clutter $\mathcal{A}$, if we have
\begin{equation*}
|B\cap A_i|>p\cdot|B|\; ,
\end{equation*}
for each $i\in[\alpha]$. The {\em $0$-committees\/} of the clutter~$\mathcal{A}$ are its {\em blocking sets}.

\noindent$\bullet$ For a nontrivial clutter $\mathcal{A}:=\{A_1,\ldots,A_{\alpha}\}\subset\mathbf{2}^{[t]}$ on the ground set~$E_t$, we have
\begin{equation*}
\#\bigl(\mathfrak{B}(\mathcal{A})^{\triangledown}\cap\tbinom{E_t}{k}\bigr)
=\tbinom{t}{k}+\sum_{S\subseteq[\alpha]\colon |S|>0}(-1)^{|S|}\cdot\tbinom{t-|\bigcup_{s\in S}A_s|}{k}
\; ,\ \ \ \ \ 1\leq k\leq t\; .
\end{equation*}
Several ways to count the blocking $k$-sets of clutters are mentioned in~\cite{M-PROM}.

\section{Increasing families of blocking sets, and blockers: Set covering problems}

In this section we recall
the set covering problem(s); see, e.g.,~\cite[Sect.~2.4]{CCZ} and~\cite[Ch.~1]{Co}.

Let ${\boldsymbol{\chi}}(A):=(\chi_1(A),\ldots,\chi_t(A))\in\{0,1\}^t$ denote the familiar row {\em characteristic vector\/} of a subset $A$ of the ground set $E_t$, defined for each element~$j\in E_t$ by
\begin{equation*}
\chi_j(A):={\small\begin{cases}
1\; , & \text{if $j\in A$},\\
0\; , & \text{if $j\not\in A$}.
\end{cases}
}
\end{equation*}
If $\mathcal{A}:=\{A_1,\ldots,A_{\alpha}\}\subset\mathbf{2}^{[t]}$ is a nontrivial clutter on
$E_t$, then
\begin{equation}
\label{eq:2}
\boldsymbol{A}:=\boldsymbol{A}(\mathcal{A}):=\left(\begin{smallmatrix}
{\boldsymbol{\chi}}(A_1)\vspace{-2mm}\\ \vdots\\ {\boldsymbol{\chi}}(A_{\alpha})
\end{smallmatrix}
\right)
\end{equation}
is its {\em incidence matrix}.

Consider\footnote{$\quad$We will denote by~$\mathbbm{1}$ and~$\mathbbm{2}$ the $\alpha$-dimensional column vectors $\left(\begin{smallmatrix}1\vspace{-2mm}\\ \vdots\\ 1\end{smallmatrix}\right)$ and $\left(\begin{smallmatrix}2\vspace{-2mm}\\ \vdots\\ 2\end{smallmatrix}\right)$, respectively.}
the {\em set covering collection}
\begin{equation}
\label{eq:24}
\widetilde{\boldsymbol{\mathcal{S}}}:=\widetilde{\boldsymbol{\mathcal{S}}}{}^{\mathtt{C}}(\boldsymbol{A}):=\bigl\{\tilde{\mathbf{z}}\in\{0,1\}^t\colon \boldsymbol{A}\tilde{\mathbf{z}}^{\top}\geq\mathbbm{1}\bigr\}\; ,
\end{equation}
which is the collection of {\em characteristic vectors\/} of the {\em blocking sets\/} of the clutter~$\mathcal{A}$, that is,
\begin{equation*}
\widetilde{\boldsymbol{\mathcal{S}}}=
\bigl\{\boldsymbol{\chi}(B)\colon B\in\mathfrak{B}(\mathcal{A})^{\triangledown}\bigr\}\; , \ \ \text{and}\ \ \ \
\mathfrak{B}(\mathcal{A})^{\triangledown}=
\bigl\{\supp(\tilde{\boldsymbol{z}})\colon \tilde{\boldsymbol{z}}\in \widetilde{\boldsymbol{\mathcal{S}}}\bigr\}
\; .
\end{equation*}
The latter expression just rephrases the convention according to which the {\em supports\/} of the vectors in the collection~$\widetilde{\boldsymbol{\mathcal{S}}}\subset\{0,1\}^t$ are the {\em blocking sets\/} of the clutter~$\mathcal{A}$.

Let us redefine the collection
\begin{equation*}
\widetilde{\boldsymbol{\mathcal{S}}}:=
\Bigl\{\tilde{\mathbf{z}}\in\{0,1\}^t\colon \left(\begin{smallmatrix}
\boldsymbol{\chi}(A_1)\vspace{-2mm}\\ \vdots\\ \boldsymbol{\chi}(A_{\alpha})
\end{smallmatrix}
\right)\tilde{\mathbf{z}}^{\top}\geq\mathbbm{1}\Bigr\}
\end{equation*}
as
\begin{equation}
\label{eq:23}
\widetilde{\boldsymbol{\mathcal{S}}}:=\biggl\{\tfrac{1}{2}(\mathrm{T}^{(+)}-\mathbf{z})\in\{0,1\}^t\colon \left(\begin{smallmatrix}
\frac{1}{2}(\mathrm{T}^{(+)}-T^1)\vspace{-2mm}\\ \vdots\\ \frac{1}{2}(\mathrm{T}^{(+)}-T^{\alpha})
\end{smallmatrix}
\right)\cdot\tfrac{1}{2}(\mathrm{T}^{(+)}-\mathbf{z})^{\top}\geq\mathbbm{1}\biggr\}\; ,
\end{equation}
where the vertices $T^i$ of the discrete hypercube~$\{1,-1\}^t$ and the vector of unknowns $\mathbf{z}\in\{1,-1\}^t$ are given by
\begin{align*}
T^i:={}_{-A_i}\mathrm{T}^{(+)}&=\mathrm{T}^{(+)}-2\boldsymbol{\chi}(A_i)
\; ,\ \ \ i\in[\alpha]\; ,\\
\intertext{and}
\mathbf{z}:\!&=\mathrm{T}^{(+)}-2\tilde{\mathbf{z}}\; .
\end{align*}

Let us now associate with the collection~$\widetilde{\boldsymbol{\mathcal{S}}}\subset\{0,1\}^t$, described in~(\ref{eq:23}), a collection~\mbox{$\boldsymbol{\mathcal{S}}:=\boldsymbol{\mathcal{S}}{}^{\mathtt{C}}(\boldsymbol{A})\subset\{1,-1\}^t$,} defined by
\begin{equation*}
\boldsymbol{\mathcal{S}}:=
\biggl\{\mathbf{z}\in\{1,-1\}^t\colon \boldsymbol{A}\mathbf{z}^{\top}\leq
\left(\begin{smallmatrix}
|(T^1)^-|\vspace{-2mm}\\ \vdots\\ |(T^{\alpha})^-|
\end{smallmatrix}
\right)
-2\cdot\mathbbm{1}
\biggr\}\; ,
\end{equation*}
that is, the collection
\begin{equation}
\label{eq:25}
\boldsymbol{\mathcal{S}}:=
\biggl\{\mathbf{z}\in\{1,-1\}^t\colon \boldsymbol{A}\mathbf{z}^{\top}\leq
\left(\begin{smallmatrix}
|A_1|\vspace{-2mm}\\ \vdots\\ |A_{\alpha}|
\end{smallmatrix}
\right)
-\mathbbm{2}
\biggr\}\; .
\end{equation}

We have defined the twin collections $\widetilde{\boldsymbol{\mathcal{S}}}\subset\{0,1\}^t$ and $\boldsymbol{\mathcal{S}}\subset\{1,-1\}^t$, given in~(\ref{eq:24}) and~(\ref{eq:25}), respectively, that are equipped with the {\em bijections\/}
$\widetilde{\boldsymbol{\mathcal{S}}}\to\boldsymbol{\mathcal{S}}\colon$ \mbox{$\widetilde{T}\mapsto\mathrm{T}^{(+)}-2\widetilde{T}$,} and
$\boldsymbol{\mathcal{S}}\to\widetilde{\boldsymbol{\mathcal{S}}}\colon$ $T\mapsto\tfrac{1}{2}(\mathrm{T}^{(+)}-T)$; see Example~\ref{th:3}.

\begin{example}
\label{th:3}
Consider the clutter $\mathcal{A}:=\{A_1,A_{\alpha:=2}\}:=\bigl\{\{1,2\},\{2,3\}\bigr\}$, on the ground set $E_{t:=3}:=\{1,2,3\}$, with its incidence matrix
\begin{equation*}
\boldsymbol{A}:=\boldsymbol{A}(\mathcal{A})=\left(
\begin{smallmatrix}
1&1&0\\
0&1&1
\end{smallmatrix}
\right)\; .
\end{equation*}
The {\em set covering $\{0,1\}$-collection}
\begin{equation*}
\begin{split}
\widetilde{\boldsymbol{\mathcal{S}}}:\!&=
\bigl\{\tilde{\mathbf{z}}\in\{0,1\}^t\colon \boldsymbol{A}\tilde{\mathbf{z}}^{\top}\geq\mathbbm{1}\bigr\}\\
&=\bigl\{\tilde{\mathbf{z}}\in\{0,1\}^3\colon \left(
\begin{smallmatrix}
1&1&0\\
0&1&1
\end{smallmatrix}
\right)\tilde{\mathbf{z}}^{\top}\geq\left(
\begin{smallmatrix}
1\\
1
\end{smallmatrix}
\right)\bigr\}
\end{split}
\end{equation*}
is the collection
\begin{equation*}
\begin{split}
\widetilde{\boldsymbol{\mathcal{S}}}&=
\bigl\{
\left(\begin{smallmatrix}0&1&0\end{smallmatrix}\right),
\left(\begin{smallmatrix}1&1&0\end{smallmatrix}\right),
\left(\begin{smallmatrix}1&0&1\end{smallmatrix}\right),
\left(\begin{smallmatrix}0&1&1\end{smallmatrix}\right),
\left(\begin{smallmatrix}1&1&1\end{smallmatrix}\right)\bigr\}\\
&=
\bigl\{\boldsymbol{\chi}(\{2\}),\boldsymbol{\chi}(\{1,2\}),\boldsymbol{\chi}(\{1,3\}),
\boldsymbol{\chi}(\{2,3\}),\boldsymbol{\chi}(\{1,2,3\})\bigr\}\; .
\end{split}
\end{equation*}
The {\em set covering $\{1,-1\}$-collection}
\begin{equation*}
\begin{split}
\boldsymbol{\mathcal{S}}:\!&=
\biggl\{\mathbf{z}\in\{1,-1\}^t\colon
\boldsymbol{A}\mathbf{z}^{\top}
\leq
\left(\begin{smallmatrix}
|A_1|\vspace{-2mm}\\ \vdots\\ |A_{\alpha}|
\end{smallmatrix}
\right)
-\mathbbm{2}
\biggr\}
\\
&=
\biggl\{\mathbf{z}\in\{1,-1\}^3\colon
\left(\begin{smallmatrix}
1&1&0\\ 0&1&1
\end{smallmatrix}
\right)\mathbf{z}^{\top}
\leq
\left(\begin{smallmatrix}
2\\ 2
\end{smallmatrix}
\right)
-\mathbbm{2}
=\left(\begin{smallmatrix}
0\\ 0
\end{smallmatrix}
\right)
\biggr\}
\end{split}
\end{equation*}
is the collection
\begin{equation*}
\begin{split}
\boldsymbol{\mathcal{S}}&=
\bigl\{
\left(\begin{smallmatrix}1& -1& 1\end{smallmatrix}\right),
\left(\begin{smallmatrix}-1& -1& 1\end{smallmatrix}\right),
\left(\begin{smallmatrix}-1& 1& -1\end{smallmatrix}\right),
\left(\begin{smallmatrix}1& -1& -1\end{smallmatrix}\right),
\left(\begin{smallmatrix}-1& -1& -1\end{smallmatrix}\right)
\bigr\}
\\
&=\bigl\{{}_{-\{2\}}\mathrm{T}^{(+)}, {}_{-\{1,2\}}\mathrm{T}^{(+)}, {}_{-\{1,3\}}\mathrm{T}^{(+)},
{}_{-\{2,3\}}\mathrm{T}^{(+)}, {}_{-\{1,2,3\}}\mathrm{T}^{(+)}\bigr\}\; .
\end{split}
\end{equation*}
\end{example}

\noindent{}$\bullet$ Let $\boldsymbol{w}\in\mathbb{R}^t$ be a row vector of {\em nonnegative\/} weights. The {\em set covering\/} problems are
\begin{equation*}
\min\{\boldsymbol{w}\tilde{\mathbf{z}}^{\top}\colon \tilde{\mathbf{z}}\in\widetilde{\boldsymbol{\mathcal{S}}}\}
=\min\bigl\{\boldsymbol{w}\cdot\tfrac{1}{2}(\mathrm{T}^{(+)}-\mathbf{z})^{\top}\colon \mathbf{z}\in\boldsymbol{\mathcal{S}}\bigr\}
\; .
\end{equation*}

\noindent$\bullet$ Suppose $\boldsymbol{w}:=\mathrm{T}^{(+)}$, and consider the (\!{\em unweighted}) set covering problems
\begin{multline*}
\tau(\mathcal{A}):=\min\{\hwt(\tilde{\mathbf{z}})\colon \tilde{\mathbf{z}}\in\widetilde{\boldsymbol{\mathcal{S}}}\}=\min\{\mathrm{T}^{(+)}\tilde{\mathbf{z}}^{\top}\colon \tilde{\mathbf{z}}\in\widetilde{\boldsymbol{\mathcal{S}}}\}
\\=
\min\bigl\{\mathrm{T}^{(+)}\cdot\tfrac{1}{2}(\mathrm{T}^{(+)}-\mathbf{z})^{\top}\colon \mathbf{z}\in\boldsymbol{\mathcal{S}}\bigr\}
=\min\bigl\{\tfrac{1}{2}(t-\underbrace{\mathrm{T}^{(+)}\mathbf{z}^{\top}}_{t-2|\mathbf{z}^-|})
\colon \mathbf{z}\in\boldsymbol{\mathcal{S}}\bigr\}
\\=
\min\bigl\{|\mathbf{z}^-|\colon \mathbf{z}\in\boldsymbol{\mathcal{S}}\bigr\}=:\tau(\mathcal{A})\; ,
\end{multline*}
that is, the problem
\begin{equation}
\label{eq:29}
\underbrace{
\min\{\mathrm{T}^{(+)}\tilde{\mathbf{z}}^{\top}\colon \tilde{\mathbf{z}}\in\widetilde{\boldsymbol{\mathcal{S}}}\}
}_{
\tau(\mathcal{A}):=\min\{\hwt(\tilde{\mathbf{z}})\colon \tilde{\mathbf{z}}\in\widetilde{\boldsymbol{\mathcal{S}}}\}
}
:=\min\{\mathrm{T}^{(+)}\tilde{\mathbf{z}}^{\top}\colon \tilde{\mathbf{z}}
\in\{0,1\}^t,\
\boldsymbol{A}\tilde{\mathbf{z}}\geq\boldsymbol{1}\}\; ,
\end{equation}
and the problem
\begin{gather}
\nonumber
\underbrace{
\tfrac{1}{2}\min\bigl\{t-\mathrm{T}^{(+)}\mathbf{z}^{\top}
\colon \mathbf{z}\in\boldsymbol{\mathcal{S}}\bigr\}
}_{
\tau(\mathcal{A}):=\min\bigl\{|\mathbf{z}^-|\colon \mathbf{z}\in\boldsymbol{\mathcal{S}}\bigr\}
}
\\
\nonumber
:=
\frac{1}{2}\min\biggl\{t-\mathrm{T}^{(+)}\mathbf{z}^{\top}
\colon \mathbf{z}\in\{1,-1\}^t,\
\boldsymbol{A}
\mathbf{z}^{\top}
\leq
\left(\begin{smallmatrix}
|A_1|\vspace{-2mm}\\ \vdots\\ |A_{\alpha}|
\end{smallmatrix}
\right)
-\mathbbm{2}
\biggr\}
\\
\nonumber
=\frac{1}{2}\cdot\!\left(t-\max\biggl\{\mathrm{T}^{(+)}\mathbf{z}^{\top}
\colon \mathbf{z}\in\{1,-1\}^t,\
\boldsymbol{A}
\mathbf{z}^{\top}
\leq
\left(\begin{smallmatrix}
|A_1|\vspace{-2mm}\\ \vdots\\ |A_{\alpha}|
\end{smallmatrix}
\right)
-\mathbbm{2}
\biggr\}\right)
\\
\label{eq:30}
=:\underbrace{
\tfrac{1}{2}\cdot\left(t-\max\bigl\{\mathrm{T}^{(+)}\mathbf{z}^{\top}
\colon \mathbf{z}\in\boldsymbol{\mathcal{S}}\bigr\}\right)
}_{
\tau(\mathcal{A}):=\min\bigl\{|\mathbf{z}^-|\colon \mathbf{z}\in\boldsymbol{\mathcal{S}}\bigr\}
}
\; .
\end{gather}

For vectors $\tilde{\boldsymbol{z}}\in\widetilde{\boldsymbol{\mathcal{S}}}$ and $\boldsymbol{z}\in\boldsymbol{\mathcal{S}}$, where $\tilde{\boldsymbol{z}}:=\tfrac{1}{2}(\mathrm{T}^{(+)}-\boldsymbol{z})$, we have the inclusions
\begin{align*}
\tilde{\boldsymbol{z}}&\in
\Arg\min\{\mathrm{T}^{(+)}\tilde{\mathbf{z}}^{\top}\colon \tilde{\mathbf{z}}
\in \widetilde{\boldsymbol{\mathcal{S}}}\}\; ,\\
\boldsymbol{z}&\in
\Arg\max\{\mathrm{T}^{(+)}\mathbf{z}^{\top}\colon \mathbf{z}
\in \boldsymbol{\mathcal{S}}\}\; ,
\end{align*}
that is, $\tilde{\boldsymbol{z}}$ and~$\boldsymbol{z}$ provide the solution
to the problems~(\ref{eq:29}) and~(\ref{eq:30}), respectively, if and only if the member
\begin{equation*}
B:=\supp(\tilde{\boldsymbol{z}})=\boldsymbol{z}^-\in\mathfrak{B}(\mathcal{A})
\end{equation*}
of the blocker of the clutter $\mathcal{A}$ has the {\em minimum\/} cardinality
\begin{equation*}
|B|=\tau(\mathcal{A})\; .
\end{equation*}

\noindent$\bullet$ We conclude this section by noting that the rows of incidence matrices~$\boldsymbol{A}$, as well as the vectors in the set covering collections $\widetilde{\boldsymbol{\mathcal{S}}}
\subset\{0,1\}^t$ and $\boldsymbol{\mathcal{S}}\subset\{1,-1\}^t$, admit their decompositions with respect to symmetric cycles in the corresponding hypercube graphs~$\widetilde{\boldsymbol{H}}(t,2)$ and~$\boldsymbol{H}(t,2)$.

\section{Families of subsets of the ground set $E_t$: Characteristic vectors and characteristic topes}

The {\em generation\/} of fundamental combinatorial objects is extensively treated in~\cite{Knuth-4A}.

\noindent$\bullet$ Consider the family $\tbinom{E_t}{s}$, for some $s$, where $0\leq s\leq t$. We denote this family of all $s$-subsets $L^s_j\subseteq E_t$, ordered {\em lexicographically}, by $\overrightarrow{\tbinom{E_t}{s}}=:(L^s_1,\ldots,L^s_{\binom{t}{s}})$.

For an $s$-uniform clutter $\mathcal{G}:=\{G_1,\ldots,G_k\}\subseteq \tbinom{E_t}{s}$, we define its row {\em characteristic vector\/} $\boldsymbol{\gamma}^{(s)}(\mathcal{G}):=(\gamma^{(s)}_1(\mathcal{G}),\ldots,\gamma^{(s)}_{\binom{t}{s}}(\mathcal{G}))\in\{0,1\}^{\binom{t}{s}}$ in the familiar way: for each $j$, where~$1\leq j\leq\tbinom{t}{s}$, we set
\begin{equation*}
\gamma^{(s)}_j(\mathcal{G}):=\begin{cases}
1\; , & \text{if $\overrightarrow{\tbinom{E_t}{s}}\ni L^s_j\in \mathcal{G}$},\\
\quad\vspace{-3mm}\\
0\; , & \text{if $\overrightarrow{\tbinom{E_t}{s}}\ni L^s_j\not\in \mathcal{G}$};
\end{cases}
\end{equation*}
see~(\ref{eq:14})--(\ref{eq:15}) in Example~\ref{th:1}.

Now, given an arbitrary family $\mathcal{F}\subseteq\mathbf{2}^{[t]}$, we set
\begin{equation*}
\boldsymbol{\gamma}^{(s)}(\mathcal{F}):=\boldsymbol{\gamma}^{(s)}(\mathcal{F}\cap\tbinom{E_t}{s})\; ,\ \ \ 0\leq s\leq t\; ,
\end{equation*}
and in a natural way\footnote{$\quad$Indeed, this is the
most popular ordering of all subsets of a nonempty finite set, the {\em lexicographic ordering subordinated to cardinality}, that is,
considering smaller subsets first, and in case of equal cardinality, the lexicographic
ordering, see~\cite[p.~77]{Gr}.} we define the {\em characteristic vector\/}
$\boldsymbol{\gamma}(\mathcal{F})
:=(\gamma_1(\mathcal{F}),\ldots,$ $\gamma_{2^t}(\mathcal{F}))\in\{0,1\}^{2^t}$ of the family~$\mathcal{F}$ to be the {\em concatenation\/}
\begin{equation*}
\boldsymbol{\gamma}(\mathcal{F})
:=\boldsymbol{\gamma}^{(0)}(\mathcal{F})\;\centerdot\;
\boldsymbol{\gamma}^{(1)}(\mathcal{F})
\;\centerdot\;\cdots
\;\centerdot\;
\boldsymbol{\gamma}^{(t-1)}(\mathcal{F})
\;\centerdot\;
\boldsymbol{\gamma}^{(t)}(\mathcal{F})\; ;
\end{equation*}
see~(\ref{eq:16})--(\ref{eq:17}).

\noindent$\bullet$ The characteristic vector $\boldsymbol{\gamma}(\mathbf{2}^{[t]})=\mathrm{T}_{2^t}^{(+)}$, whose components are all $1$'s, describes the linearly ordered {\em power set\/} $\mathbf{2}^{[t]}$
of the ground set $E_t$; see~(\ref{eq:20}).

\noindent$\bullet$ The Hamming weights $\hwt(\boldsymbol{\gamma}^{(s)}(\mathcal{F}))$ of the vectors $\boldsymbol{\gamma}^{(s)}(\mathcal{F})$, $0\leq s\leq t$, are the components $f_s(\mathcal{F};t)$ of the so-called {\em long $f$-vectors\/} $\boldsymbol{f}(\mathcal{F};t)$ associated with families $\mathcal{F}\subseteq\mathbf{2}^{[t]}$, see~\cite[Sect.~2.1]{M-PROM}.

\noindent$\bullet$ If $\mathcal{F}'\subseteq\mathbf{2}^{[t]}$ and~$\mathcal{F}''\subseteq\mathbf{2}^{[t]}$ are families of subsets of the ground set~$E_t$, then we will use the {\em componentwise product\/} of their characteristic vectors
\begin{equation*}
\boldsymbol{\gamma}(\mathcal{F}')\ast\boldsymbol{\gamma}(\mathcal{F}'')
:=\bigl(\gamma_1(\mathcal{F}')\cdot\gamma_1(\mathcal{F}''),\ldots,\gamma_{2^t}(\mathcal{F}')\cdot\gamma_{2^t}(\mathcal{F}'')\bigr)
\in\{0,1\}^{2^t}
\end{equation*}
to describe\footnote{$\quad$ The notation {\tiny$\sideset{}{^\ast}\prod$} will be used to denote the {\em componentwise product\/} of several vectors.} the {\em intersection\/} of these families:
\begin{equation*}
\boldsymbol{\gamma}(\mathcal{F}'\cap\mathcal{F}'')=\boldsymbol{\gamma}(\mathcal{F}')\ast\boldsymbol{\gamma}(\mathcal{F}'')\; .
\end{equation*}

\noindent$\bullet$ Let $\varGamma(k)$ denote the subset $A\subseteq E_t$, for which the characteristic vector
of the corresponding
one-member clutter $\{A\}$ on $E_t$
by convention is the $k$th standard unit vector~$\boldsymbol{\sigma}(k)$ of the space~$\mathbb{R}^{2^t}$; we thus use the map
\begin{equation*}
\varGamma: [2^t]\to\mathbf{2}^{[t]}\; ,\ \ \ \ \ k\mapsto A\colon\ \boldsymbol{\gamma}(\{A\})=\boldsymbol{\sigma}(k)\in\{0,1\}^{2^t}\; ;
\end{equation*}
see~(\ref{eq:18})--(\ref{eq:19}).
Conversely, we denote by $\varGamma^{-1}(A)$, where $A\subseteq E_t$, the position number $k$ such that the vector~$\boldsymbol{\sigma}(k)$ is the characteristic vector of the
one-member clutter $\{A\}$ on $E_t$:
\begin{equation*}
\varGamma^{-1}: \mathbf{2}^{[t]}\to [2^t]\; ,\ \ \ \ \ A\mapsto k\colon\ \boldsymbol{\sigma}(k)=\boldsymbol{\gamma}(\{A\})\in\{0,1\}^{2^t}\; ;
\end{equation*}
see~(\ref{eq:18})--(\ref{eq:19}).

By construction, we have the implications
\begin{gather}
\label{eq:21}
\ell',\ell''\in [2^t]\; ,\ \ \ell'<\ell''\ \ \ \Longrightarrow\ \ \ |\varGamma(\ell')|\leq|\varGamma(\ell'')|\; ;\\
\nonumber
A,B\in\mathbf{2}^{[t]}\; ,\ \ |A|<|B|\ \ \ \Longrightarrow\ \ \ \varGamma^{-1}(A)<\varGamma^{-1}(B)\; ,
\end{gather}
and, in particular,
\begin{equation*}
A,B\in\mathbf{2}^{[t]}\; ,\ \ A\subsetneqq B\ \ \ \Longrightarrow\ \ \ \varGamma^{-1}(A)<\varGamma^{-1}(B)\; .
\end{equation*}
Note also that for any index $\ell\in[2^t]$, the disjoint union
\begin{equation*}
\varGamma(\ell)\;\,\dot\cup\;\,\varGamma(2^t-\ell+1)=E_t
\end{equation*}
is a {\em partition\/} of the ground set.

\begin{example}
\label{th:1}
Suppose $t:=3$, and $E_t=\{1,2,3\}$. We have
{\small
\begin{align}
\label{eq:14}
\boldsymbol{\gamma}^{(0)}(\tbinom{E_t}{0}):=\boldsymbol{\gamma}^{(0)}(\{\hat{0}\})&=(1)\in\{0,1\}^{\binom{t}{0}}\; ,\\
\boldsymbol{\gamma}^{(1)}(\tbinom{E_t}{1}):=\boldsymbol{\gamma}^{(1)}(\{\{1\},\{2\},\{3\}\})&=(1,1,1)\in\{0,1\}^{\binom{t}{1}}\; ,\\
\boldsymbol{\gamma}^{(2)}(\tbinom{E_t}{2}):=\boldsymbol{\gamma}^{(2)}(\{\{1,2\},\{1,3\},\{2,3\}\})&=(1,1,1)\in\{0,1\}^{\binom{t}{2}}\; ,\\
\boldsymbol{\gamma}^{(t)}(\tbinom{E_t}{t}):=\boldsymbol{\gamma}^{(t)}(\{\{1,2,3\}\})&=(1)\in\{0,1\}^{\binom{t}{t}}\; ,\\
\boldsymbol{\gamma}^{(1)}(\{\{2\}\})
&=(0,1,0)\in\{0,1\}^{\binom{t}{1}}\; ,\\
\boldsymbol{\gamma}^{(2)}(\{\{1,2\},\{2,3\}\})
&=(1,0,1)\in\{0,1\}^{\binom{t}{2}}\; ,\\
\label{eq:15}
\boldsymbol{\gamma}^{(2)}(\{\{1,3\}\})
&=(0,1,0)\in\{0,1\}^{\binom{t}{2}}\; , \\
\intertext{\normalsize{and}}
\label{eq:16}
\boldsymbol{\gamma}(\emptyset)&=(0,0,0,0,0,0,0,0)\in\{0,1\}^{2^t}\; ,\\
\boldsymbol{\gamma}(\tbinom{E_t}{0})
&=(1,0,0,0,0,0,0,0)
\; ,\\
\boldsymbol{\gamma}(\tbinom{E_t}{1})
&=(0,1,1,1,0,0,0,0)
\; ,\\
\boldsymbol{\gamma}(\tbinom{E_t}{2})
&=(0,0,0,0,1,1,1,0)
\; ,\\
\boldsymbol{\gamma}(\tbinom{E_t}{t})
&=(0,0,0,0,0,0,0,1)
\; ,\\
\label{eq:20}
\boldsymbol{\gamma}(\mathbf{2}^{[t]})=\mathrm{T}_{2^t}^{(+)}:\!&=(1,1,1,1,1,1,1,1)
\; .
\intertext{\label{page:1}\normalsize{If $\mathcal{A}:=\{A_1,A_2\}$ and $\mathcal{B}:=\{B_1,B_2\}$ are clutters on $E_t$, where $A_1:=\{1,2\}$, \mbox{$A_2:=\{2,3\}$,}
$B_1:=\{1,3\}$, $B_2:=\{2\}$, and $\mathcal{B}=\mathfrak{B}(\mathcal{A})$, then we have
}
}
\boldsymbol{\gamma}(\mathcal{A}):=\boldsymbol{\gamma}(\{A_1,A_2\}):=\boldsymbol{\gamma}(\{\{1,2\},\{2,3\}\})&=(0,0,0,0,1,0,1,0)\in\{0,1\}^{2^t}\; ,\\
\boldsymbol{\gamma}(\mathcal{B}):=\boldsymbol{\gamma}(\{B_1,B_2\}):=\boldsymbol{\gamma}(\{\{1,3\},\{2\}\})&=(0,0,1,0,0,1,0,0)
\; ,\\
\boldsymbol{\gamma}(\mathcal{A}^{\triangledown}):=\boldsymbol{\gamma}(\{\{1,2\},\{2,3\}\}^{\triangledown})&=(0,0,0,0,1,0,1,1)
\; ,\\
\boldsymbol{\gamma}(\mathcal{B}^{\triangledown}):=\boldsymbol{\gamma}(\{\{1,3\},\{2\}\}^{\triangledown})&=(0,0,1,0,1,1,1,1)
\; ,\\
\boldsymbol{\gamma}(\{A_1\}^{\triangledown}):=\boldsymbol{\gamma}(\{\{1,2\}\}^{\triangledown})&=(0,0,0,0,1,0,0,1)
\; ,\\
\boldsymbol{\gamma}(\{A_2\}^{\triangledown}):=\boldsymbol{\gamma}(\{\{2,3\}\}^{\triangledown})&=(0,0,0,0,0,0,1,1)
\; ,\\
\boldsymbol{\gamma}(\{B_1\}^{\triangledown}):=\boldsymbol{\gamma}(\{\{1,3\}\}^{\triangledown})&=(0,0,0,0,0,1,0,1)
\; ,\\
\label{eq:17}
\boldsymbol{\gamma}(\{B_2\}^{\triangledown}):=\boldsymbol{\gamma}(\{\{2\}\}^{\triangledown})&=(0,0,1,0,1,0,1,1)
\; .
\end{align}
}
We have
{\small
\begin{align}
\label{eq:18}
\varGamma(3)&=\{2\}\; , & \boldsymbol{\gamma}(\{\{2\}\})&=\boldsymbol{\sigma}(3)\in\{0,1\}^{2^t}\; , & \varGamma^{-1}(\{2\})&=3\; ,\\
\varGamma(6)&=\{1,3\}\; , & \boldsymbol{\gamma}(\{\{1,3\}\})&=\boldsymbol{\sigma}(6)\in\{0,1\}^{2^t}\; , & \varGamma^{-1}(\{1,3\})&=6\; , \\
\label{eq:19}
\varGamma(2^t)&=E_t\; , & \boldsymbol{\gamma}(\{E_t\})&=\boldsymbol{\sigma}(2^t)\in\{0,1\}^{2^t}\; , & \varGamma^{-1}(E_t)&=2^t\; .
\end{align}
}
\end{example}

\noindent$\bullet$ Given a family $\mathcal{F}\subseteq\mathbf{2}^{[t]}$ of subsets of the ground set $E_t$, we call the tope
\begin{equation}
\label{eq:68}
T_{\mathcal{F}}:={}_{-\supp(\boldsymbol{\gamma}(\mathcal{F}))}\mathrm{T}_{2^t}^{(+)}=\mathrm{T}_{2^t}^{(+)}-2\boldsymbol{\gamma}(\mathcal{F})
\end{equation}
of the oriented matroid $\mathcal{H}_{2^t}:=(E_{2^t},\{1,-1\}^{2^t})$ the {\em characteristic tope\/} of the family $\mathcal{F}$;
see Example~\ref{th:2}.

If we let $\mathrm{T}^{(+)}_{\binom{t}{s}}:=(1,\ldots,1)\in\mathbb{R}^{\binom{t}{s}}$ denote the $\tbinom{t}{s}$-dimensional row vector
of all~$1$'s, $0\leq s\leq t$, then
\begin{equation*}
T^{(s)}_{\mathcal{F}}:={}_{-\supp(\boldsymbol{\gamma}^{(s)}(\mathcal{F}))}\mathrm{T}_{\binom{t}{s}}^{(+)}
=\mathrm{T}_{\binom{t}{s}}^{(+)}-2\boldsymbol{\gamma}^{(s)}(\mathcal{F})\; , \ \ \ 0\leq s\leq t\; .
\end{equation*}

\begin{example}
\label{th:2}
Suppose $t:=3$ and $E_t=\{1,2,3\}$. If
$\mathcal{A}$ and
$\mathcal{B}=\mathfrak{B}(\mathcal{A})$ are clutters on the ground set~$E_t$, mentioned in~Example~{\rm\ref{th:1}} on page~{\rm\pageref{page:1}}, then we have
{\small
\begin{align*}
\boldsymbol{\gamma}(\mathcal{A})
&=
(0,\phantom{-}0,\phantom{-}0,\phantom{-}0,\phantom{-}1,\phantom{-}0,\phantom{-}1,\phantom{-}0)
\in\{0,1\}^{2^t}\; ,\\
T_{\mathcal{A}}:\!&=
(1,\phantom{-}1,\phantom{-}1,\phantom{-}1,-1,\phantom{-}1,-1,\phantom{-}1)\in\{1,-1\}^{2^t}\; ;\\
\boldsymbol{\gamma}(\mathcal{A}^{\triangledown})
&=(0,\phantom{-}0,\phantom{-}0,\phantom{-}0,\phantom{-}1,\phantom{-}0,\phantom{-}1,\phantom{-}1)\; ,\\
T_{\mathcal{A}^{\triangledown}}:\!&=(1,\phantom{-}1,\phantom{-}1,\phantom{-}1,-1,\phantom{-}1,-1,-1)\; ;\\
\boldsymbol{\gamma}(\mathcal{B})
&=(0,\phantom{-}0,\phantom{-}1,\phantom{-}0,\phantom{-}0,\phantom{-}1,\phantom{-}0,\phantom{-}0)\; ,\\
T_{\mathcal{B}}
:\!&=(1,\phantom{-}1,-1,\phantom{-}1,\phantom{-}1,-1,\phantom{-}1,\phantom{-}1)\; ;\\
\boldsymbol{\gamma}(\mathcal{B}^{\triangledown})
&=(0,\phantom{-}0,\phantom{-}1,\phantom{-}0,\phantom{-}1,\phantom{-}1,\phantom{-}1,\phantom{-}1)\; ,\\
T_{\mathcal{B}^{\triangledown}}
:\!&=(1,\phantom{-}1,-1,\phantom{-}1,-1,-1,-1,-1)
\; .
\end{align*}
}
\end{example}

\section{Increasing families of blocking sets, and blockers: Characteristic vectors and characteristic topes}

In this section, we begin with somewhat sophisticated restatements of the simple basic observations~(\ref{eq:22}), (\ref{eq:27}) and~(\ref{eq:28}) on set families in terms of their characteristic vectors.

\noindent$\bullet$
For a nontrivial clutter $\mathcal{A}\subset\mathbf{2}^{[t]}$ on the ground set~$E_t$, we have
\begin{gather*}
\boldsymbol{\gamma}(\mathfrak{B}(\mathcal{A})^{\triangledown})
\ast\boldsymbol{\gamma}(\mathcal{A}^{\triangledown})
\\
=\rn(\boldsymbol{\gamma}(\mathcal{A}^{\triangledown}))
\ast\boldsymbol{\gamma}(\mathcal{A}^{\triangledown}):=
\bigl(\mathrm{T}_{2^t}^{(+)}-\boldsymbol{\gamma}(\mathcal{A}^{\triangledown})\cdot\overline{\mathbf{U}}(2^t)\bigr)
\ast\boldsymbol{\gamma}(\mathcal{A}^{\triangledown})
\\
=\boldsymbol{\gamma}(\mathfrak{B}(\mathcal{A})^{\triangledown})
\ast\rn(\boldsymbol{\gamma}(\mathfrak{B}(\mathcal{A})^{\triangledown}))
:=
\boldsymbol{\gamma}(\mathfrak{B}(\mathcal{A})^{\triangledown})
\ast\bigl(\mathrm{T}_{2^t}^{(+)}-\boldsymbol{\gamma}(\mathfrak{B}(\mathcal{A})^{\triangledown})\cdot\overline{\mathbf{U}}(2^t)\bigr)
\\=
\begin{cases}
\boldsymbol{\gamma}(\mathcal{A}^{\triangledown})\; , & \text{if $\#\mathcal{A}^{\triangledown}<2^{t-1}$,}\\
\boldsymbol{\gamma}(\mathfrak{B}(\mathcal{A})^{\triangledown})\; , & \text{if $\#\mathcal{A}^{\triangledown}>2^{t-1}$,}\\
\boldsymbol{\gamma}(\mathcal{A}^{\triangledown})
=\boldsymbol{\gamma}(\mathfrak{B}(\mathcal{A})^{\triangledown})\; , & \text{if $\#\mathcal{A}^{\triangledown}=2^{t-1}$.}
\end{cases}
\end{gather*}

\begin{remark}
\label{th:4}
Let $\mathcal{A}\subset\mathbf{2}^{[t]}$ be a nontrivial clutter on the ground set~$E_t$. We have
\begin{gather}
\nonumber
\gamma_1(\mathfrak{B}(\mathcal{A})^{\triangledown})=0\; ,\ \ \ \text{and}\ \ \ \ \gamma_{2^t}(\mathfrak{B}(\mathcal{A})^{\triangledown})=1\; ;\\
\label{eq:3}
\boldsymbol{\gamma}(\mathfrak{B}(\mathcal{A})^{\triangledown})=\rn(\boldsymbol{\gamma}(\mathcal{A}^{\triangledown}))
:=\mathrm{T}^{(+)}_{2^t}-\boldsymbol{\gamma}(\mathcal{A}^{\triangledown})\cdot\overline{\mathbf{U}}(2^t)\; ;\\
\label{eq:86}
T_{\mathfrak{B}(\mathcal{A})^{\triangledown}}=\ro(T_{\mathcal{A}^{\triangledown}})
:=-T_{\mathcal{A}^{\triangledown}}\cdot\overline{\mathbf{U}}(2^t)\; ;\\
\nonumber
\underbrace{
\hwt(\boldsymbol{\gamma}(\mathfrak{B}(\mathcal{A})^{\triangledown}))
}_{\#\mathfrak{B}(\mathcal{A})^{\triangledown}}
\; +\;
\underbrace{
\hwt(\boldsymbol{\gamma}(\mathcal{A}^{\triangledown}))
}_{\#\mathcal{A}^{\triangledown}}
=2^t\; ;\\
\nonumber
\underbrace{
|(T_{\mathfrak{B}(\mathcal{A})^{\triangledown}})^-|
}_{\#\mathfrak{B}(\mathcal{A})^{\triangledown}}
\; +\;
\underbrace{
|(T_{\mathcal{A}^{\triangledown}})^-|
}_{\#\mathcal{A}^{\triangledown}}
=2^t\; ;\\
\label{eq:26}
\boldsymbol{\gamma}^{(s)}(\mathfrak{B}(\mathcal{A})^{\triangledown})=\rn(\boldsymbol{\gamma}^{(t-s)}(\mathcal{A}^{\triangledown}))
:=\mathrm{T}^{(+)}_{\binom{t}{s}}-\boldsymbol{\gamma}^{(t-s)}(\mathcal{A}^{\triangledown})\cdot\overline{\mathbf{U}}(\tbinom{t}{s})\; ,
\ \ \ 0\leq s\leq t\; ;\\
\label{eq:87}
T^{(s)}_{\mathfrak{B}(\mathcal{A})^{\triangledown}}=\ro(T^{(t-s)}_{\mathcal{A}^{\triangledown}})
:=-T^{(t-s)}_{\mathcal{A}^{\triangledown}}\cdot\overline{\mathbf{U}}(\tbinom{t}{s})\; ,\ \ \ 0\leq s\leq t\; ;\\
\nonumber
\underbrace{
\hwt(\boldsymbol{\gamma}^{(s)}(\mathfrak{B}(\mathcal{A})^{\triangledown}))
}_{\#(\mathfrak{B}(\mathcal{A})^{\triangledown}\cap\binom{E_t}{s})}
\; +\;
\underbrace{
\hwt(\boldsymbol{\gamma}^{(t-s)}(\mathcal{A}^{\triangledown}))
}_{\#(\mathcal{A}^{\triangledown}\cap\binom{E_t}{t-s})}
=\tbinom{t}{s}\; , \ \ \ 0\leq s\leq t\; ;\\
\nonumber
\underbrace{
|(T^{(s)}_{\mathfrak{B}(\mathcal{A})^{\triangledown}})^-|
}_{\#(\mathfrak{B}(\mathcal{A})^{\triangledown}\cap\binom{E_t}{s})}
+\;
\underbrace{
|(T^{(t-s)}_{\mathcal{A}^{\triangledown}})^-|
}_{\#(\mathcal{A}^{\triangledown}\cap\binom{E_t}{t-s})}
=\tbinom{t}{s}\; , \ \ \ 0\leq s\leq t\; .
\end{gather}
\end{remark}
In addition to~(\ref{eq:3}) and~(\ref{eq:86}), relations~(\ref{eq:26}) and~(\ref{eq:87}) imply that
\begin{gather*}
\boldsymbol{\gamma}(\mathfrak{B}(\mathcal{A})^{\triangledown})=
\underbrace{
\rn(\boldsymbol{\gamma}^{(t)}(\mathcal{A}^{\triangledown}))
}_{(0)}
\;\centerdot\;
\rn(\boldsymbol{\gamma}^{(t-1)}(\mathcal{A}^{\triangledown}))
\;\centerdot\; \cdots
\;\centerdot\;
\rn(\boldsymbol{\gamma}^{(1)}(\mathcal{A}^{\triangledown}))
\;\centerdot\;
\underbrace{
\rn(\boldsymbol{\gamma}^{(0)}(\mathcal{A}^{\triangledown}))
}_{(1)}
\; ,\\
T_{\mathfrak{B}(\mathcal{A})}^{\triangledown}=
\underbrace{
\ro(T^{(t)}_{\mathcal{A}^{\triangledown}})
}_{(1)}
\;\centerdot\;
\ro(T^{(t-1)}_{\mathcal{A}^{\triangledown}})
\;\centerdot\; \cdots
\;\centerdot\;
\ro(T^{(1)}_{\mathcal{A}^{\triangledown}})
\;\centerdot\;
\underbrace{
\ro(T^{(0)}_{\mathcal{A}^{\triangledown}})
}_{(-1)}\; .
\end{gather*}

\noindent$\bullet$ In view of~(\ref{eq:21}), if
\begin{equation*}
\ell^{\star}:=\min\supp(\boldsymbol{\gamma}(\mathfrak{B}(\mathcal{A})^{\triangledown}))
=\min\,(T_{\mathfrak{B}(\mathcal{A})^{\triangledown}})^-
\; ,
\end{equation*}
then the member $\varGamma(\ell^{\star})$ of the blocker $\mathfrak{B}(\mathcal{A})$ of a nontrivial clutter $\mathcal{A}\subset\mathbf{2}^{[t]}$ is a blocking set of {\em minimum\/} cardinality for~$\mathcal{A}$, that is,
the vectors~$\boldsymbol{\chi}(\varGamma(\ell^{\star}))$ and ${}_{-\varGamma(\ell^{\star})}\mathrm{T}^{(+)}$ provide
the solution
$|\varGamma(\ell^{\star})|=\tau(\mathcal{A})$
to the set covering problems~(\ref{eq:29}) and~(\ref{eq:30}), respectively:
\begin{align*}
\boldsymbol{\chi}(\varGamma(\ell^{\star}))&\in
\Arg\min\{\mathrm{T}^{(+)}\tilde{\mathbf{z}}^{\top}\colon \tilde{\mathbf{z}}
\in \widetilde{\boldsymbol{\mathcal{S}}}\}\; ,\\
{}_{-\varGamma(\ell^{\star})}\mathrm{T}^{(+)}&\in
\Arg\max\{\mathrm{T}^{(+)}\mathbf{z}^{\top}\colon \mathbf{z}
\in \boldsymbol{\mathcal{S}}\}\; .
\end{align*}

\subsection{A clutter $\{\{a\}\}$} $\quad$ \label{section:1}

Let $\{\{a\}\}$ be a (nontrivial) clutter on the ground set $E_t$, whose only member is a {\em one-element\/} subset $\{a\}\subset E_t$.

\subsubsection{The principal increasing family of blocking sets $\mathfrak{B}(\{\{a\}\})^{\triangledown}=\{\{a\}\}^{\triangledown}$} $\quad$

\noindent$\bullet$ The {\em increasing family\/} of {\em blocking sets\/} $\mathfrak{B}(\{\{a\}\})^{\triangledown}$ of the {\em self-dual\/} clutter~$\{\{a\}\}$ coincides with the principal increasing family~$\{\{a\}\}^{\triangledown}$.

\noindent$\bullet$ We will use the notation $\widetilde{\boldsymbol{\mathfrak{a}}}(a):=\widetilde{\boldsymbol{\mathfrak{a}}}(a;2^t)$ and $\boldsymbol{\mathfrak{a}}(a):=\boldsymbol{\mathfrak{a}}(a;2^t)$ to denote the characteristic vector and the characteristic tope, respectively, that are associated with the principal increasing family~$\{\{a\}\}^{\triangledown}=\mathfrak{B}(\{\{a\}\})^{\triangledown}$\,:
\begin{gather*}
\widetilde{\boldsymbol{\mathfrak{a}}}(a):=\boldsymbol{\gamma}(\{\{a\}\}^{\triangledown})=
\boldsymbol{\gamma}(\mathfrak{B}(\{\{a\}\})^{\triangledown})
\in\{0,1\}^{2^t}\; ,
\\
\boldsymbol{\mathfrak{a}}(a):=T_{\{\{a\}\}^{\triangledown}}=T_{\mathfrak{B}(\{\{a\}\})^{\triangledown}}
\in\{1,-1\}^{2^t}\; .
\end{gather*}
We have
\begin{equation*}
\widetilde{\boldsymbol{\mathfrak{a}}}(a)=
\underbrace{
(0)
}_{\widetilde{\boldsymbol{\mathfrak{a}}}^{(0)}(a)}
\;\centerdot\;
\underbrace{
\boldsymbol{\chi}(\{a\})
}_{\widetilde{\boldsymbol{\mathfrak{a}}}^{(1)}(a)}
\;\centerdot\;
\underbrace{
\boldsymbol{\gamma}^{(2)}(\{\{a\}\}^{\triangledown})
}_{\widetilde{\boldsymbol{\mathfrak{a}}}^{(2)}(a)}
\;\centerdot\;
\cdots
\;\centerdot\;
\underbrace{
\boldsymbol{\gamma}^{(t-1)}(\{\{a\}\}^{\triangledown})
}_{\widetilde{\boldsymbol{\mathfrak{a}}}^{(t-1)}(a)}
\;\centerdot\;
\underbrace{
(1)
}_{\widetilde{\boldsymbol{\mathfrak{a}}}^{(t)}(a)}\; ,
\end{equation*}
see~(\ref{eq:39}), (\ref{eq:40}) and~(\ref{eq:41}) in~Example~\ref{th:16};
\begin{equation*}
\boldsymbol{\mathfrak{a}}(a)=
\underbrace{
(1)
}_{\boldsymbol{\mathfrak{a}}^{(0)}(a)}
\;\centerdot\;
\underbrace{
{}_{-\{a\}}\mathrm{T}^{(+)}
}_{\boldsymbol{\mathfrak{a}}^{(1)}(a)}
\;\centerdot\;
\underbrace{
T_{\{\{a\}\}^{\triangledown}}^{(2)}
}_{\boldsymbol{\mathfrak{a}}^{(2)}(a)}
\;\centerdot\;
\cdots
\;\centerdot\;
\underbrace{
T_{\{\{a\}\}^{\triangledown}}^{(t-1)}
}_{\boldsymbol{\mathfrak{a}}^{(t-1)}(a)}
\;\centerdot\;
\underbrace{
(-1)
}_{\boldsymbol{\mathfrak{a}}^{(t)}(a)}
\; ,
\end{equation*}
see~(\ref{eq:42}), (\ref{eq:43}) and~(\ref{eq:44}).

\begin{remark}[see~Remark~{\rm\ref{th:4}}, and cf.~Remark~{\rm\ref{th:13}}]
\label{th:12}
Note that
\begin{gather}
\label{eq:1}
\widetilde{\boldsymbol{\mathfrak{a}}}(a)=\rn(\widetilde{\boldsymbol{\mathfrak{a}}}(a))\; ,\ \ \ \text{and}\ \ \ \
\boldsymbol{\mathfrak{a}}(a)=\ro(\boldsymbol{\mathfrak{a}}(a))\; ;\\
\nonumber
\hwt(\widetilde{\boldsymbol{\mathfrak{a}}}(a))=|\boldsymbol{\mathfrak{a}}(a)^-|=
\#\{\{a\}\}^{\triangledown}=\#\mathfrak{B}(\{\{a\}\})^{\triangledown}=2^{t-1}\; ;\\
\nonumber
\widetilde{\boldsymbol{\mathfrak{a}}}^{(s)}(a)=\rn(\widetilde{\boldsymbol{\mathfrak{a}}}^{(t-s)}(a))\; ,\ \ \ \text{and}\ \ \ \
\boldsymbol{\mathfrak{a}}^{(s)}(a)=\ro(\boldsymbol{\mathfrak{a}}^{(t-s)}(a))\; ,\ \ \ 0\leq s\leq t\; ;\\
\nonumber
\hwt(\widetilde{\boldsymbol{\mathfrak{a}}}^{(s)}(a))=
|\boldsymbol{\mathfrak{a}}^{(s)}(a)^-|=\tbinom{t-1}{s-1}\; ,\ \ \ 0\leq s\leq t\; .
\end{gather}
\end{remark}

\subsubsection{The blocker $\mathfrak{B}(\{\{a\}\})=\{\{a\}\}$} $\quad$

\noindent$\bullet$ The {\em blocker\/} $\mathfrak{B}(\{\{a\}\})$ coincides with the self-dual clutter $\{\{a\}\}$.

\noindent$\bullet$ We associate with the clutter $\mathfrak{B}(\{\{a\}\})=\{\{a\}\}$ its characteristic vector~$\boldsymbol{\gamma}(\{\{a\}\})=\boldsymbol{\gamma}(\mathfrak{B}(\{\{a\}\}))\in\{0,1\}^{2^t}$ and its characteristic tope~\mbox{$T_{\{\{a\}\}}$} $=T_{\mathfrak{B}(\{\{a\}\})}\in\{1,-1\}^{2^t}$,
where
\begin{equation*}
\boldsymbol{\gamma}(\{\{a\}\})=\boldsymbol{\gamma}(\mathfrak{B}(\{\{a\}\}))=
(0)\;\centerdot\;\boldsymbol{\chi}(\{a\})
\;\centerdot\;(0,\ldots,0)\;\centerdot\;\cdots\;\centerdot\;(0)\; ,
\end{equation*}
see~(\ref{eq:45}), (\ref{eq:46}) and~(\ref{eq:47}) in~Example~\ref{th:16}.

\subsubsection{More on the principal increasing family~$\mathfrak{B}(\{\{a\}\})^{\triangledown}=\{\{a\}\}^{\triangledown}$} $\quad$

In view of~(\ref{eq:1}), we can make the following observation:
\begin{remark}[cf.~Remark~\ref{th:15}] \label{th:14}
For any element $a\in E_t$, we have
\begin{equation*}
\begin{split}
\min\{i\in E_{2^t}\colon T_{\{\{a\}\}^{\triangledown}}(i)=-1\}:\!&=\min\{i\in E_{2^t}\colon \gamma_i(\{\{a\}\}^{\triangledown})=1\}\\
:\!&=\varGamma^{-1}(\{a\})=1+a\; ,\\
\max\{i\in E_{2^t}\colon T_{\{\{a\}\}^{\triangledown}}(i)=1\}:\!&=\max\{i\in E_{2^t}\colon \gamma_i(\{\{a\}\}^{\triangledown})=0\}\\
&=2^t-\min\{a\}=2^t-a\; ;\\
\underbrace{\min\{j\in E_{2^t}\colon T_{\mathfrak{B}(\{\{a\}\})^{\triangledown}}(j)=-1\}}_{
\min\{j\in E_{2^t}\colon T_{\mathfrak{B}(\{\{a\}\})}(j)=-1\}
}
:\!&=
\underbrace{
\min\{j\in E_{2^t}\colon \gamma_j(\mathfrak{B}(\{\{a\}\})^{\triangledown})=1\}
}_{
\min\{j\in E_{2^t}\colon \gamma_j(\mathfrak{B}(\{\{a\}\}))=1\}
}
\\
&=1+\min\{a\}=1+a\; ,\\
\max\{j\in E_{2^t}\colon T_{\mathfrak{B}(\{\{a\}\})^{\triangledown}}(j)=1\}:\!&=\max\{j\in E_{2^t}\colon \gamma_j(\mathfrak{B}(\{\{a\}\})^{\triangledown})=0\}\\
&=1+2^t-\varGamma^{-1}(\{a\})=2^t-a\; .
\end{split}
\end{equation*}
\end{remark}

\noindent$\bullet$
We have
\begin{equation*}
\{\{a\}\}^{\triangledown}\ \dot{\cup}\ \{E_t-\{a\}\}^{\vartriangle}=\mathbf{2}^{[t]}\; .
\end{equation*}
Let us denote by $\widetilde{\boldsymbol{\mathfrak{c}}}(a):=\widetilde{\boldsymbol{\mathfrak{c}}}(a;2^t)$
and~$\boldsymbol{\mathfrak{c}}(a):=\boldsymbol{\mathfrak{c}}(a;2^t)$ the characteristic vector and the characteristic tope, respectively, of the principal decreasing family~$\{E_t-\{a\}\}^{\vartriangle}$:
\begin{equation*}
\widetilde{\boldsymbol{\mathfrak{c}}}(a):=\boldsymbol{\gamma}(\{E_t-\{a\}\}^{\vartriangle})\in\{0,1\}^{2^t}\; ,
\end{equation*}
see~(\ref{eq:48}), (\ref{eq:49}) and~(\ref{eq:50});
\begin{equation*}
\boldsymbol{\mathfrak{c}}(a):=T_{\{E_t-\{a\}\}^{\vartriangle}}\in\{1,-1\}^{2^t}\; ,
\end{equation*}
see~(\ref{eq:51}), (\ref{eq:52}) and~(\ref{eq:53}).
We have
\begin{gather*}
\widetilde{\boldsymbol{\mathfrak{c}}}(a)=\mathrm{T}_{2^t}^{(+)}-\widetilde{\boldsymbol{\mathfrak{a}}}(a)=
\widetilde{\boldsymbol{\mathfrak{a}}}(a)\cdot\overline{\mathbf{U}}(2^t)\; ,\\
\boldsymbol{\mathfrak{c}}(a)=-\boldsymbol{\mathfrak{a}}(a)=\boldsymbol{\mathfrak{a}}(a)\cdot\overline{\mathbf{U}}(2^t)\; .
\end{gather*}

\noindent$\bullet$ For any two-element subset $\{i,j\}\subset E_t$ of the ground set, we have
\begin{equation*}
\#(\{\{i\}\}^{\triangledown}\cap\{\{j\}\}^{\triangledown})
=\#\{\{i,j\}\}^{\triangledown}=2^{t-2}\; ,
\end{equation*}
and
\begin{equation*}
\#(\;\underbrace{
(\mathbf{2}^{[t]}-\{\{i\}\}^{\triangledown})
}_{
\{E_t-\{i\}\}^{\vartriangle}
}
\,\cap\,
\underbrace{
(\mathbf{2}^{[t]}-\{\{j\}\}^{\triangledown})
}_{
\{E_t-\{j\}\}^{\vartriangle}
}
\;)
=\#\{E_t-\{i,j\}\}^{\vartriangle}=2^{t-2}\; .
\end{equation*}
Thus, if $i$ and $j$ are elements of the ground set $E_t$, and $i\neq j$, then we have
\begin{equation*}
\begin{split}
d\bigl(\widetilde{\boldsymbol{\mathfrak{a}}}(i),\widetilde{\boldsymbol{\mathfrak{a}}}(j)\bigr)
=d\bigl(\boldsymbol{\mathfrak{a}}(i),\boldsymbol{\mathfrak{a}}(j)\bigr)
&=d\bigl(\widetilde{\boldsymbol{\mathfrak{c}}}(i),\widetilde{\boldsymbol{\mathfrak{c}}}(j)\bigr)
=d\bigl(\boldsymbol{\mathfrak{c}}(i),\boldsymbol{\mathfrak{c}}(j)\bigr)
\\
&=2^{t-1}\; .
\end{split}
\end{equation*}

\begin{remark}
For any two elements $i$ and $j$ of the ground set~$E_t$ we have
\begin{equation*}
\bigl\langle \boldsymbol{\mathfrak{a}}(i),\boldsymbol{\mathfrak{a}}(j)\bigr\rangle=
\bigl\langle \boldsymbol{\mathfrak{c}}(i),\boldsymbol{\mathfrak{c}}(j)\bigr\rangle=
\delta_{i,j}\cdot 2^t\; .
\end{equation*}
In other words, the sequences of\/ $t$ row vectors
\begin{equation*}
\bigl(\,\tfrac{1}{\sqrt{2^{t}}}\cdot
\boldsymbol{\mathfrak{a}}(1),\; \tfrac{1}{\sqrt{2^{t}}}\cdot
\boldsymbol{\mathfrak{a}}(2),\;\ldots,\;\tfrac{1}{\sqrt{2^{t}}}\cdot
\boldsymbol{\mathfrak{a}}(t)\bigr)\subset\mathbb{R}^{2^t}
\end{equation*}
and
\begin{equation*}
\bigl(\,\tfrac{1}{\sqrt{2^{t}}}\cdot
\boldsymbol{\mathfrak{c}}(t),\;\tfrac{1}{\sqrt{2^{t}}}\cdot \boldsymbol{\mathfrak{c}}(t-1),\;\ldots,\;\tfrac{1}{\sqrt{2^{t}}}\cdot
\boldsymbol{\mathfrak{c}}(1)\bigr)\subset\mathbb{R}^{2^t}
\end{equation*}
are both {\em orthonormal}.
\end{remark}

\begin{example}\label{th:16}
Suppose $t:=3$, and $E_t=\{1,2,3\}$. We have
{\small
\begin{align}
\label{eq:39}
\widetilde{\boldsymbol{\mathfrak{a}}}(1):=\boldsymbol{\gamma}(\{\{1\}\}^{\triangledown})=\boldsymbol{\gamma}(\mathfrak{B}(\{\{1\}\})^{\triangledown})
&=(\phantom{-}0,\phantom{-}1,\phantom{-}0,\phantom{-}0,\phantom{-}1,\phantom{-}1,\phantom{-}0,\phantom{-}1)\in\{0,1\}^{2^t}\; ,\\
\label{eq:42}
\boldsymbol{\mathfrak{a}}(1):=T_{\{\{1\}\}^{\triangledown}}=T_{\mathfrak{B}(\{\{1\}\})^{\triangledown}}
&=(\phantom{-}1,-1,\phantom{-}1,\phantom{-}1,-1,-1,\phantom{-}1,-1)\in\{1,-1\}^{2^t}\; ,
\\
\label{eq:48}
\widetilde{\boldsymbol{\mathfrak{c}}}(1):=\boldsymbol{\gamma}(\{E_t-\{1\}\}^{\vartriangle})
&=(\phantom{-}1,\phantom{-}0,\phantom{-}1,\phantom{-}1,\phantom{-}0,\phantom{-}0,\phantom{-}1,\phantom{-}0)
\; ,\\
\label{eq:51}
\boldsymbol{\mathfrak{c}}(1):=T_{\{E_t-\{1\}\}^{\vartriangle}}
&=(-1,\phantom{-}1,-1,-1,\phantom{-}1,\phantom{-}1,-1,\phantom{-}1)
\; ,
\\
\label{eq:45}
\boldsymbol{\gamma}(\{\{1\}\})=\boldsymbol{\gamma}(\mathfrak{B}(\{\{1\}\}))
&=(\phantom{-}0,\phantom{-}1,\phantom{-}0,\phantom{-}0,\phantom{-}0,\phantom{-}0,\phantom{-}0,\phantom{-}0)
\; ,\\
\label{eq:40}
\widetilde{\boldsymbol{\mathfrak{a}}}(2):=\boldsymbol{\gamma}(\{\{2\}\}^{\triangledown})=\boldsymbol{\gamma}(\mathfrak{B}(\{\{2\}\})^{\triangledown})
&=(\phantom{-}0,\phantom{-}0,\phantom{-}1,\phantom{-}0,\phantom{-}1,\phantom{-}0,\phantom{-}1,\phantom{-}1)
\; ,\\
\label{eq:43}
\boldsymbol{\mathfrak{a}}(2):=T_{\{\{2\}\}^{\triangledown}}=T_{\mathfrak{B}(\{\{2\}\})^{\triangledown}}
&=(\phantom{-}1,\phantom{-}1,-1,\phantom{-}1,-1,\phantom{-}1,-1,-1)
\; ,\\
\label{eq:49}
\widetilde{\boldsymbol{\mathfrak{c}}}(2):=\boldsymbol{\gamma}(\{E_t-\{2\}\}^{\vartriangle})
&=(\phantom{-}1,\phantom{-}1,\phantom{-}0,\phantom{-}1,\phantom{-}0,\phantom{-}1,\phantom{-}0,\phantom{-}0)
\; ,\\
\label{eq:52}
\boldsymbol{\mathfrak{c}}(2):=T_{\{E_t-\{2\}\}^{\vartriangle}}
&=(-1,-1,\phantom{-}1,-1,\phantom{-}1,-1,\phantom{-}1,\phantom{-}1)
\; ,\\
\label{eq:46}
\boldsymbol{\gamma}(\{\{2\}\})=\boldsymbol{\gamma}(\mathfrak{B}(\{\{2\}\}))
&=(\phantom{-}0,\phantom{-}0,\phantom{-}1,\phantom{-}0,\phantom{-}0,\phantom{-}0,\phantom{-}0,\phantom{-}0)
\; ,\\
\label{eq:41}
\widetilde{\boldsymbol{\mathfrak{a}}}(t):=\boldsymbol{\gamma}(\{\{t\}\}^{\triangledown})=\boldsymbol{\gamma}(\mathfrak{B}(\{\{t\}\})^{\triangledown})
&=(\phantom{-}0,\phantom{-}0,\phantom{-}0,\phantom{-}1,\phantom{-}0,\phantom{-}1,\phantom{-}1,\phantom{-}1)
\; ,\\
\label{eq:44}
\boldsymbol{\mathfrak{a}}(t):=T_{\{\{t\}\}^{\triangledown}}=T_{\mathfrak{B}(\{\{t\}\})^{\triangledown}}
&=(\phantom{-}1,\phantom{-}1,\phantom{-}1,-1,\phantom{-}1,-1,-1,-1)
\; ,\\
\label{eq:50}
\widetilde{\boldsymbol{\mathfrak{c}}}(t):=\boldsymbol{\gamma}(\{E_t-\{t\}\}^{\vartriangle})
&=(\phantom{-}1,\phantom{-}1,\phantom{-}1,\phantom{-}0,\phantom{-}1,\phantom{-}0,\phantom{-}0,\phantom{-}0)
\; ,\\
\label{eq:53}
\boldsymbol{\mathfrak{c}}(t):=T_{\{E_t-\{t\}\}^{\vartriangle}}
&=(-1,-1,-1,\phantom{-}1,-1,\phantom{-}1,\phantom{-}1,\phantom{-}1)
\; ,\\
\label{eq:47}
\boldsymbol{\gamma}(\{\{t\}\})=\boldsymbol{\gamma}(\mathfrak{B}(\{\{t\}\}))
&=(\phantom{-}0,\phantom{-}0,\phantom{-}0,\phantom{-}1,\phantom{-}0,\phantom{-}0,\phantom{-}0,\phantom{-}0)
\; ,\\
\label{eq:60}
\boldsymbol{\gamma}(\{\{1,2\}\})&=(\phantom{-}0,\phantom{-}0,\phantom{-}0,\phantom{-}0,\phantom{-}1,\phantom{-}0,\phantom{-}0,\phantom{-}0)
\; ,\\
\label{eq:54}
\boldsymbol{\gamma}(\{\{1,2\}\}^{\triangledown})&=(\phantom{-}0,\phantom{-}0,\phantom{-}0,\phantom{-}0,\phantom{-}1,\phantom{-}0,\phantom{-}0,\phantom{-}1)
\; ,\\
\label{eq:63}
\boldsymbol{\gamma}(\mathfrak{B}(\{\{1,2\}\}))&=(\phantom{-}0,\phantom{-}1,\phantom{-}1,\phantom{-}0,\phantom{-}0,\phantom{-}0,\phantom{-}0,\phantom{-}0)
\; ,\\
\label{eq:57}
\boldsymbol{\gamma}(\mathfrak{B}(\{\{1,2\}\})^{\triangledown})&=(\phantom{-}0,\phantom{-}1,\phantom{-}1,\phantom{-}0,\phantom{-}1,\phantom{-}1,\phantom{-}1,\phantom{-}1)
\; ,\\
\label{eq:61}
\boldsymbol{\gamma}(\{\{1,t\}\})&=(\phantom{-}0,\phantom{-}0,\phantom{-}0,\phantom{-}0,\phantom{-}0,\phantom{-}1,\phantom{-}0,\phantom{-}0)
\; ,\\
\label{eq:55}
\boldsymbol{\gamma}(\{\{1,t\}\}^{\triangledown})&=(\phantom{-}0,\phantom{-}0,\phantom{-}0,\phantom{-}0,\phantom{-}0,\phantom{-}1,\phantom{-}0,\phantom{-}1)
\; ,\\
\label{eq:64}
\boldsymbol{\gamma}(\mathfrak{B}(\{\{1,t\}\}))&=(\phantom{-}0,\phantom{-}1,\phantom{-}0,\phantom{-}1,\phantom{-}0,\phantom{-}0,\phantom{-}0,\phantom{-}0)
\; ,\\
\label{eq:58}
\boldsymbol{\gamma}(\mathfrak{B}(\{\{1,t\}\})^{\triangledown})&=(\phantom{-}0,\phantom{-}1,\phantom{-}0,\phantom{-}1,\phantom{-}1,\phantom{-}1,\phantom{-}1,\phantom{-}1)
\; ,\\
\label{eq:62}
\boldsymbol{\gamma}(\{\{2,t\}\})&=(\phantom{-}0,\phantom{-}0,\phantom{-}0,\phantom{-}0,\phantom{-}0,\phantom{-}0,\phantom{-}1,\phantom{-}0)
\; ,\\
\label{eq:56}
\boldsymbol{\gamma}(\{\{2,t\}\}^{\triangledown})&=(\phantom{-}0,\phantom{-}0,\phantom{-}0,\phantom{-}0,\phantom{-}0,\phantom{-}0,\phantom{-}1,\phantom{-}1)
\; ,\\
\label{eq:65}
\boldsymbol{\gamma}(\mathfrak{B}(\{\{2,t\}\}))&=(\phantom{-}0,\phantom{-}0,\phantom{-}1,\phantom{-}1,\phantom{-}0,\phantom{-}0,\phantom{-}0,\phantom{-}0)
\; ,\\
\label{eq:59}
\boldsymbol{\gamma}(\mathfrak{B}(\{\{2,t\}\})^{\triangledown})&=(\phantom{-}0,\phantom{-}0,\phantom{-}1,\phantom{-}1,\phantom{-}1,\phantom{-}1,\phantom{-}1,\phantom{-}1)
\; .
\end{align}
}
\end{example}

\noindent$\bullet$ Given a nontrivial clutter~$\mathcal{A}:=\{A_1,\ldots,A_{\alpha}\}\subset\mathbf{2}^{[t]}$ on the ground set~$E_t$, such that $\mathcal{A}\neq\{E_t\}$, we have
\begin{equation*}
\begin{split}
\boldsymbol{\gamma}(\mathcal{A}):\!&=\sum_{i\in[\alpha]}
\boldsymbol{\sigma}(\varGamma^{-1}(A_i))
=\boldsymbol{\gamma}(\mathcal{A}^{\triangledown})\ast\boldsymbol{\gamma}(\mathcal{A}^{\vartriangle})
=\sum_{i\in[\alpha]}\bigl(\boldsymbol{\gamma}(\{A_i\}^{\triangledown})\ast\boldsymbol{\gamma}(\{A_i\}^{\vartriangle})\bigr)
\\
&=\sum_{i\in[\alpha]}\Bigl(\sideset{}{^\ast}\prod_{a^i\in A_i}\widetilde{\boldsymbol{\mathfrak{a}}}(a^i)
\;\ast\sideset{}{^\ast}\prod_{c^i\in E_t-A_i}\widetilde{\boldsymbol{\mathfrak{c}}}(c^i)\Bigr)
\\
&=\sum_{i\in[\alpha]}\Bigl(\sideset{}{^\ast}\prod_{a^i\in A_i}
\widetilde{\boldsymbol{\mathfrak{a}}}(a^i)
\;
\ast\Bigl(\Bigl(\sideset{}{^\ast}\prod_{c^i\in E_t-A_i}\widetilde{\boldsymbol{\mathfrak{a}}}(c^i)\Bigr)\cdot\overline{\mathbf{U}}(2^t)\Bigr)
\Bigr)
\\
&=\sum_{i\in[\alpha]}\Bigl(\Bigl(\Bigl(\sideset{}{^\ast}\prod_{a^i\in A_i}\widetilde{\boldsymbol{\mathfrak{c}}}(a^i)\Bigr)\cdot\overline{\mathbf{U}}(2^t)\Bigr)
\;\ast\sideset{}{^\ast}\prod_{c^i\in E_t-A_i}\widetilde{\boldsymbol{\mathfrak{c}}}(c^i)\Bigr)
\; .
\end{split}
\end{equation*}

\subsection{A clutter $\{A\}$} $\quad$ \label{section:2}

Let $\{A\}$ be a (nontrivial) clutter on the ground set~$E_t$, whose only member is a {\em nonempty subset} $A\subseteq E_t$.

\subsubsection{The increasing family of blocking sets $\mathfrak{B}(\{A\})^{\triangledown}=\{\{a\}\colon a\in A\}^{\triangledown}$} $\quad$

\noindent$\bullet$ The {\em family\/} of {\em blocking sets\/} $\mathfrak{B}(\{A\})^{\triangledown}$ of the clutter $\{A\}$ is the increasing family~$\{\{a\}\colon a\in A\}^{\triangledown}$.

We have
\begin{equation*}
\{A\}^{\triangledown}=\bigcap_{a\in A}\{\{a\}\}^{\triangledown}\; , \ \ \ \text{and}\ \ \ \
\mathfrak{B}(\{A\})^{\triangledown}=\bigcup_{a\in A}\{\{a\}\}^{\triangledown}\; .
\end{equation*}

Let us associate with the increasing families $\{A\}^{\triangledown}$ and $\mathfrak{B}(\{A\})^{\triangledown}$ their characteristic vectors
$\boldsymbol{\gamma}(\{A\}^{\triangledown})\in\{0,1\}^{2^t}$ and~$\boldsymbol{\gamma}(\mathfrak{B}(\{A\})^{\triangledown})\in\{0,1\}^{2^t}$,
and their characteristic topes~$T_{\{A\}^{\triangledown}}\in\{1,-1\}^{2^t}$ and~$T_{\mathfrak{B}(\{A\})^{\triangledown}}\in\{1,-1\}^{2^t}$,
where
\begin{multline*}
\boldsymbol{\gamma}(\{A\}^{\triangledown})=\boldsymbol{\gamma}\Bigl(\bigcap_{a\in A}\{\{a\}\}^{\triangledown}\Bigr)
\\
=
\underbrace{
(0)
}_{
\boldsymbol{\gamma}^{(0)}(\{A\}^{\triangledown})
}
\;\centerdot\;\cdots\;\centerdot\;
\underbrace{
(0,\ldots,0)
}_{
\boldsymbol{\gamma}^{(|A|-1)}(\{A\}^{\triangledown})
}
\;\centerdot\;
\underbrace{
\boldsymbol{\gamma}^{(|A|)}(\{A\})
}_{
\boldsymbol{\gamma}^{(|A|)}(\{A\}^{\triangledown})
}
\\
\;\centerdot\;
\underbrace{
\boldsymbol{\gamma}^{(|A|+1)}\Bigl(\bigcap_{a\in A}\{\{a\}\}^{\triangledown}\Bigr)
}_{
\boldsymbol{\gamma}^{(|A|+1)}(\{A\}^{\triangledown})
}
\;\centerdot\;
\cdots
\;\centerdot\;
\underbrace{
\boldsymbol{\gamma}^{(t-1)}\Bigl(\bigcap_{a\in A}\{\{a\}\}^{\triangledown}\Bigr)
}_{
\boldsymbol{\gamma}^{(t-1)}(\{A\}^{\triangledown})
}
\;\centerdot\;
\underbrace{
(1)
}_{
\boldsymbol{\gamma}^{(t)}(\{A\}^{\triangledown})
}
\; ,
\end{multline*}
see~(\ref{eq:54}), (\ref{eq:55}) and~(\ref{eq:56}) in Example~\ref{th:16},
and
\begin{multline*}
\boldsymbol{\gamma}(\mathfrak{B}(\{A\})^{\triangledown})=\boldsymbol{\gamma}\Bigl(\bigcup_{a\in A}\{\{a\}\}^{\triangledown}\Bigr)
\\=
\underbrace{
(0)
}_{
\boldsymbol{\gamma}^{(0)}(\mathfrak{B}(\{A\})^{\triangledown})
}
\;\centerdot\;
\underbrace{
\boldsymbol{\chi}(A)
}_{\boldsymbol{\gamma}^{(1)}(\mathfrak{B}(\{A\})^{\triangledown})
}
\\
\;\centerdot\;
\underbrace{
\boldsymbol{\gamma}^{(2)}(\bigcup_{a\in A}\{\{a\}\}^{\triangledown})
}_{
\boldsymbol{\gamma}^{(2)}(\mathfrak{B}(\{A\})^{\triangledown})
}
\;\centerdot\;
\cdots
\;\centerdot\;
\underbrace{
\boldsymbol{\gamma}^{(t-1)}(\bigcup_{a\in A}\{\{a\}\}^{\triangledown})
}_{
\boldsymbol{\gamma}^{(t-1)}(\mathfrak{B}(\{A\})^{\triangledown})
}
\;\centerdot\;
\underbrace{
(1)
}_{
\boldsymbol{\gamma}^{(t)}(\mathfrak{B}(\{A\})^{\triangledown})
}
\; ,
\end{multline*}
see~(\ref{eq:57}), (\ref{eq:58}) and~(\ref{eq:59}).

\begin{remark}[see~Remark~{\rm\ref{th:4}}, and cf.~Remark~{\rm\ref{th:12}}]
\label{th:13}
Note that
\begin{gather*}
\nonumber
\boldsymbol{\gamma}(\mathfrak{B}(\{A\})^{\triangledown})=\rn(\boldsymbol{\gamma}(\{A\}^{\triangledown}))
:=\mathrm{T}^{(+)}_{2^t}-\boldsymbol{\gamma}(\{A\}^{\triangledown})\cdot\overline{\mathbf{U}}(2^t)\; ;\\
\nonumber
T_{\mathfrak{B}(\{A\})^{\triangledown}}=\ro(T_{\{A\}^{\triangledown}})
:=-T_{\{A\}^{\triangledown}}\cdot\overline{\mathbf{U}}(2^t)\; ;\\
\nonumber
\underbrace{
\hwt(\boldsymbol{\gamma}(\mathfrak{B}(\{A\})^{\triangledown}))
}_{\#\mathfrak{B}(\{A\})^{\triangledown}}=\underbrace{
|(T_{\mathfrak{B}(\{A\})^{\triangledown}})^-|
}_{\#\mathfrak{B}(\{A\})^{\triangledown}}=2^t-2^{t-|A|}\; ;\\
\underbrace{
\hwt(\boldsymbol{\gamma}(\{A\}^{\triangledown}))
}_{\#\{A\}^{\triangledown}}=
\underbrace{
|(T_{\{A\}^{\triangledown}})^-|
}_{\#\{A\}^{\triangledown}}
=2^{t-|A|}\; ;\\
\boldsymbol{\gamma}^{(s)}(\mathfrak{B}(\{A\})^{\triangledown})=\rn(\boldsymbol{\gamma}^{(t-s)}(\{A\}^{\triangledown}))
:=\mathrm{T}^{(+)}_{\binom{t}{s}}-\boldsymbol{\gamma}^{(t-s)}(\{A\}^{\triangledown})\cdot\overline{\mathbf{U}}(\tbinom{t}{s})\; ,
\ \ \ 0\leq s\leq t\; ;\\
\nonumber
T^{(s)}_{\mathfrak{B}(\{A\})^{\triangledown}}=\ro(T^{(t-s)}_{\{A\}^{\triangledown}})
:=-T^{(t-s)}_{\{A\}^{\triangledown}}\cdot\overline{\mathbf{U}}(\tbinom{t}{s})\; ,\ \ \ 0\leq s\leq t\; ;\\
\nonumber
\underbrace{
\hwt(\boldsymbol{\gamma}^{(s)}(\mathfrak{B}(\{A\})^{\triangledown}))
}_{\#(\mathfrak{B}(\{A\})^{\triangledown}\cap\binom{E_t}{s})}\; =
\underbrace{
|(T^{(s)}_{\mathfrak{B}(\{A\})^{\triangledown}})^-|
}_{\#(\mathfrak{B}(\{A\})^{\triangledown}\cap\binom{E_t}{s})}
=\tbinom{t}{s}-\tbinom{t-|A|}{s}\; ,\ \ \ 0\leq s\leq t;\\
\underbrace{
\hwt(\boldsymbol{\gamma}^{(t-s)}(\{A\}^{\triangledown}))
}_{\#(\{A\}^{\triangledown}\cap\binom{E_t}{t-s})}\;
=\underbrace{
|(T^{(t-s)}_{\{A\}^{\triangledown}})^-|
}_{\#(\{A\}^{\triangledown}\cap\binom{E_t}{t-s})}
=\tbinom{t-|A|}{s}\; , \ \ \ 0\leq s\leq t\; .
\end{gather*}
\end{remark}

\subsubsection{The blocker $\mathfrak{B}(\{A\})=\{\{a\}\colon a\in A\}$} $\quad$

\noindent$\bullet$ The {\em blocker\/}
of the clutter $\{A\}$ is the clutter
\begin{equation*}
\mathfrak{B}(\{A\})=\{\{a\}\colon a\in A\}\; .
\end{equation*}
Thus, $\#\mathfrak{B}(\{A\})=|A|$, and the members of the blocker $\mathfrak{B}(\{A\})$ are the one-element subsets of the set $A$.

\noindent$\bullet$ We associate with the clutters $\{A\}$ and $\mathfrak{B}(\{A\})$ their characteristic vectors~$\boldsymbol{\gamma}(\{A\})\in\{0,1\}^{2^t}$ and $\boldsymbol{\gamma}(\mathfrak{B}(\{A\}))\in\{0,1\}^{2^t}$, and their characteristic topes~$T_{\{A\}}\in\{1,-1\}^{2^t}$ and $T_{\mathfrak{B}(\{A\})}\in\{1,-1\}^{2^t}$, where
\begin{equation*}
\boldsymbol{\gamma}(\{A\})=
\underbrace{
(0)
}_{
\boldsymbol{\gamma}^{(0)}(\{A\})
}
\;\centerdot\;\cdots\;\centerdot\;
\underbrace{
(0,\ldots,0)
}_{
\boldsymbol{\gamma}^{(|A|-1)}(\{A\})
}
\;\centerdot\;
\boldsymbol{\gamma}^{(|A|)}(\{A\})
\;\centerdot\;
\underbrace{
(0,\ldots,0)
}_{
\boldsymbol{\gamma}^{(|A|+1)}(\{A\})
}
\centerdot\;\cdots\;\centerdot\;
\underbrace{
(0)
}_{
\boldsymbol{\gamma}^{(t)}(\{A\})
}
\; ,
\end{equation*}
see~(\ref{eq:60}), (\ref{eq:61}) and~(\ref{eq:62}) in Example~\ref{th:16}, and
\begin{equation*}
\boldsymbol{\gamma}(\mathfrak{B}(\{A\}))
=
\underbrace{
(0)
}_{
\boldsymbol{\gamma}^{(0)}(\mathfrak{B}(\{A\}))
}
\;\centerdot\;
\underbrace{
\boldsymbol{\chi}(A)
}_{
\boldsymbol{\gamma}^{(1)}(\mathfrak{B}(\{A\}))
}
\;\centerdot\;
\underbrace{
(0,\ldots,0)
}_{
\boldsymbol{\gamma}^{(2)}(\mathfrak{B}(\{A\}))
}
\;\centerdot\;\cdots\;\centerdot\;
\underbrace{
(0)
}_{
\boldsymbol{\gamma}^{(t)}(\mathfrak{B}(\{A\}))
}
\; ,
\end{equation*}
see~(\ref{eq:63}), (\ref{eq:64}) and~(\ref{eq:65}).

\subsubsection{More on the increasing families $\{A\}^{\triangledown}$ and\/ $\mathfrak{B}(\{A\})^{\triangledown}$} $\quad$

We can make the following observation:
\begin{remark}[cf.~Remark~\ref{th:14}] \label{th:15}
For a nonempty subset $A\subseteq E_t$, we have
\begin{equation*}
\begin{split}
\min\{i\in E_{2^t}\colon T_{\{A\}^{\triangledown}}(i)=-1\}:\!&=\min\{i\in E_{2^t}\colon \gamma_i(\{A\}^{\triangledown})=1\}\\
:\!&=\varGamma^{-1}(A)\; ,\\
\max\{i\in E_{2^t}\colon T_{\{A\}^{\triangledown}}(i)=1\}:\!&=\max\{i\in E_{2^t}\colon \gamma_i(\{A\}^{\triangledown})=0\}\\
&=2^t-\min A\; ;\\
\underbrace{\min\{j\in E_{2^t}\colon T_{\mathfrak{B}(\{A\})^{\triangledown}}(j)=-1\}}_{
\min\{j\in E_{2^t}\colon T_{\mathfrak{B}(\{A\})}(j)=-1\}
}
:\!&=
\underbrace{
\min\{j\in E_{2^t}\colon \gamma_j(\mathfrak{B}(\{A\})^{\triangledown})=1\}
}_{
\min\{j\in E_{2^t}\colon \gamma_j(\mathfrak{B}(\{A\}))=1\}
}
\\
&=1+\min A\; ,\\
\max\{j\in E_{2^t}\colon T_{\mathfrak{B}(\{A\})^{\triangledown}}(j)=1\}:\!&=\max\{j\in E_{2^t}\colon \gamma_j(\mathfrak{B}(\{A\})^{\triangledown})=0\}\\
&=1+2^t-\varGamma^{-1}(A)\; .
\end{split}
\end{equation*}
\end{remark}

\noindent$\bullet$ Recall that the partition
\begin{equation*}
\{A\}^{\triangledown}\;\dot{\cup}\; (\mathfrak{B}(\{A\})^{\compl})^{\vartriangle}=\mathbf{2}^{[t]}
\end{equation*}
implies that
\begin{equation*}
\mathfrak{B}(\{A\})^{\triangledown}=\{D^{\compl}\colon D\in\mathbf{2}^{[t]}-\{A\}^{\triangledown}\}\; .
\end{equation*}

\noindent$\bullet$
Note that
\begin{gather*}
\boldsymbol{\gamma}(\{A\}^{\triangledown})=\sideset{}{^\ast}\prod_{a\in A}\boldsymbol{\gamma}(\{\{a\}\}^{\triangledown})
=:\sideset{}{^\ast}\prod_{a\in A}\widetilde{\boldsymbol{\mathfrak{a}}}(a)
\\
=\sideset{}{^\ast}\prod_{a\in A}\bigl(\mathrm{T}_{2^t}^{(+)}-\widetilde{\boldsymbol{\mathfrak{c}}}(a)\bigr)
=\Bigl(\sideset{}{^\ast}\prod_{a\in A}\widetilde{\boldsymbol{\mathfrak{c}}}(a)\Bigr)\cdot\overline{\mathbf{U}}(2^t)\; ,
\end{gather*}
and recall that
\begin{equation*}
\boldsymbol{\gamma}(\mathfrak{B}(\{A\})^{\triangledown})=\rn(\boldsymbol{\gamma}(\{A\}^{\triangledown}))\; .
\end{equation*}

\begin{remark}
\label{th:8}
For a nonempty subset $A\subseteq E_t$, we have{\rm:}
\begin{itemize}
\item[\rm(i)]
\begin{equation*}
\boldsymbol{\gamma}(\{A\}^{\triangledown})=\sideset{}{^\ast}\prod_{a\in A}\widetilde{\boldsymbol{\mathfrak{a}}}(a)\; .
\end{equation*}

\item[\rm(ii)]
\begin{equation*}
\boldsymbol{\gamma}(\mathfrak{B}(\{A\})^{\triangledown})
=\mathrm{T}_{2^t}^{(+)}
-\Bigl(\sideset{}{^\ast}\prod_{a\in A}\widetilde{\boldsymbol{\mathfrak{a}}}(a)\Bigr)\cdot\overline{\mathbf{U}}(2^t)
\; .
\end{equation*}
\end{itemize}
\end{remark}

\subsection{A clutter $\mathcal{A}:=\{A_1,\ldots,A_{\alpha}\}$} $\quad$ \label{section:3}

Let $\mathcal{A}:=\{A_1,\ldots,A_{\alpha}\}$ be a nontrivial clutter on the ground set $E_t$.

\subsubsection{The increasing family of blocking sets $\mathfrak{B}(\mathcal{A})^{\triangledown}$} $\quad$

\noindent$\bullet$ See~Remark~\ref{th:4}, and note that
\begin{align*}
\mathcal{A}^{\triangledown}&=
\bigcup_{k\in[\alpha]}\{A_k\}^{\triangledown}=
\bigcup_{k\in[\alpha]}\; \bigcap_{a^k\in A_k}\{\{a^k\}\}^{\triangledown}\; ,\\
\intertext{and}
\mathfrak{B}(\mathcal{A})^{\triangledown}&=
\bigcap_{k\in[\alpha]}\mathfrak{B}(\{A_k\})^{\triangledown}=
\bigcap_{k\in[\alpha]}\;\bigcup_{a^k\in A_k}\{\{a^k\}\}^{\triangledown}\; .
\end{align*}

\noindent$\bullet$ We associate with the increasing families $\mathcal{A}^{\triangledown}$ and $\mathfrak{B}(\mathcal{A})^{\triangledown}$ their characteristic vectors~$\boldsymbol{\gamma}(\mathcal{A}^{\triangledown})\in\{0,1\}^{2^t}$ and~$\boldsymbol{\gamma}(\mathfrak{B}(\mathcal{A})^{\triangledown})\in\{0,1\}^{2^t}$, and their characteristic topes~$T_{\mathcal{A}^{\triangledown}}\in\{1,-1\}^{2^t}$ and~$T_{\mathfrak{B}(\mathcal{A})^{\triangledown}}
\in\{1,-1\}^{2^t}$,
where
\begin{multline*}
\boldsymbol{\gamma}(\mathcal{A}^{\triangledown})=\boldsymbol{\gamma}\Bigl(\bigcup_{k\in[\alpha]}\; \bigcap_{a^k\in A_k}\{\{a^k\}\}^{\triangledown}\Bigr)
\\
=
\underbrace{
(0)
}_{
\boldsymbol{\gamma}^{(0)}(\mathcal{A}^{\triangledown})
}
\;\centerdot\;
\boldsymbol{\gamma}^{(1)}(\bigcup_{k\in[\alpha]}\; \bigcap_{a^k\in A_k}\{\{a^k\}\}^{\triangledown})
\;\centerdot\;\cdots\;\centerdot\;
\boldsymbol{\gamma}^{(t-1)}(\bigcup_{k\in[\alpha]}\; \bigcap_{a^k\in A_k}\{\{a^k\}\}^{\triangledown})
\;\centerdot\;
\underbrace{
(1)
}_{
\boldsymbol{\gamma}^{(t)}(\mathcal{A}^{\triangledown})
}
\; ,
\end{multline*}
and
\begin{multline*}
\boldsymbol{\gamma}(\mathfrak{B}(\mathcal{A})^{\triangledown})=\boldsymbol{\gamma}\Bigl(\bigcap_{k\in[\alpha]}\;\bigcup_{a^k\in A_k}\{\{a^k\}\}^{\triangledown}\Bigr)
\\=
\underbrace{
(0)
}_{
\boldsymbol{\gamma}^{(0)}(\mathfrak{B}(\mathcal{A})^{\triangledown})
}
\;\centerdot\;\boldsymbol{\gamma}^{(1)}(\bigcap_{k\in[\alpha]}\;\bigcup_{a^k\in A_k}\{\{a^k\}\}^{\triangledown})
\;\centerdot\;\cdots\;\centerdot\;
\boldsymbol{\gamma}^{(t-1)}(\bigcap_{k\in[\alpha]}\;\bigcup_{a^k\in A_k}\{\{a^k\}\}^{\triangledown})
\;\centerdot\;
\underbrace{
(1)
}_{
\boldsymbol{\gamma}^{(t)}(\mathfrak{B}(\mathcal{A})^{\triangledown})
}
\; .
\end{multline*}

\subsubsection{The blocker~$\mathfrak{B}(\mathcal{A})$} $\quad$

\noindent$\bullet$ The {\em blocker\/}
of the clutter $\mathcal{A}$ is the clutter
\begin{multline*}
\mathfrak{B}(\mathcal{A})=\bmin\bigcap_{k\in[\alpha]}\mathfrak{B}(\{A_k\})^{\triangledown}
\\=\bmin\bigcap_{k\in[\alpha]}\{\{a^k\}\colon a^k\in A_k\}^{\triangledown}=
\bmin\bigcap_{k\in[\alpha]}\;\bigcup_{a^k\in A_k}\{\{a^k\}\}^{\triangledown}
\; .
\end{multline*}

\noindent$\bullet$ We associate with the clutters $\mathcal{A}$ and $\mathfrak{B}(\mathcal{A})$ their characteristic vectors $\boldsymbol{\gamma}(\mathcal{A})\in\{0,1\}^{2^t}$ and $\boldsymbol{\gamma}(\mathfrak{B}(\mathcal{A}))\in\{0,1\}^{2^t}$, and their characteristic topes \mbox{$T_{\mathcal{A}}\in\{1,-1\}^{2^t}$} and $T_{\mathfrak{B}(\mathcal{A})}\in\{1,-1\}^{2^t}$,
where
\begin{equation*}
\boldsymbol{\gamma}(\mathcal{A}):=
\underbrace{
(0)
}_{
\boldsymbol{\gamma}^{(0)}(\mathcal{A})
}
\;\centerdot\;
\boldsymbol{\gamma}^{(1)}(\mathcal{A})
\;\centerdot\;\cdots\;\centerdot\;
\boldsymbol{\gamma}^{(t)}(\mathcal{A})\; ,
\end{equation*}
and
\begin{multline*}
\boldsymbol{\gamma}(\mathfrak{B}(\mathcal{A}))=
\underbrace{
(0)
}_{
\boldsymbol{\gamma}^{(0)}(\mathfrak{B}(\mathcal{A}))
}
\;\centerdot\;\boldsymbol{\gamma}^{(1)}(\bmin\bigcap_{k\in[\alpha]}\;\bigcup_{a^k\in A_k}\{\{a^k\}\}^{\triangledown})
\\
\;\centerdot\;\cdots\;\centerdot\;
\boldsymbol{\gamma}^{(t)}(\bmin\bigcap_{k\in[\alpha]}\;\bigcup_{a^k\in A_k}\{\{a^k\}\}^{\triangledown})\; .
\end{multline*}

\subsubsection{More on the increasing families $\mathcal{A}^{\triangledown}$ and\/ $\mathfrak{B}(\mathcal{A})^{\triangledown}$} $\quad$

\noindent$\bullet$ Recall that we have
\begin{equation*}
\mathcal{A}^{\triangledown}\;\dot{\cup}\; (\mathfrak{B}(\mathcal{A})^{\compl})^{\vartriangle}=\mathbf{2}^{[t]}\; ,
\end{equation*}
that is,
\begin{equation*}
\mathfrak{B}(\mathcal{A})^{\triangledown}=\{D^{\compl}\colon D\in\mathbf{2}^{[t]}-\mathcal{A}^{\triangledown}\}\; .
\end{equation*}

\noindent$\bullet$ According to~Remark~\ref{th:8}(ii), we have
\begin{equation*}
\begin{split}
\boldsymbol{\gamma}(\mathfrak{B}(\mathcal{A})^{\triangledown})&=\sideset{}{^\ast}\prod_{i\in[\alpha]}\boldsymbol{\gamma}(\mathfrak{B}(\{A_i\})^{\triangledown})
=\sideset{}{^\ast}\prod_{i\in[\alpha]}\Bigl(\mathrm{T}_{2^t}^{(+)}
-\Bigl(\sideset{}{^\ast}\prod_{a^i\in A_i}\widetilde{\boldsymbol{\mathfrak{a}}}(a^i)\Bigr)\cdot\overline{\mathbf{U}}(2^t)\Bigr)\\
&=\Bigl(\sideset{}{^\ast}\prod_{i\in[\alpha]}\Bigl(\mathrm{T}^{(+)}_{2^t}
-\sideset{}{^\ast}\prod_{a^i\in A_i}\widetilde{\boldsymbol{\mathfrak{a}}}(a^i)\Bigr)\Bigr)\cdot\overline{\mathbf{U}}(2^t)
\; .
\end{split}
\end{equation*}
Since
\begin{equation*}
\boldsymbol{\gamma}(\mathfrak{B}(\mathcal{A})^{\triangledown})=\rn(\boldsymbol{\gamma}(\mathcal{A}^{\triangledown}))
:=\mathrm{T}_{2^t}^{(+)}
-\boldsymbol{\gamma}(\mathcal{A}^{\triangledown})\cdot\overline{\mathbf{U}}(2^t)\; ,
\end{equation*}
by~(\ref{eq:3}), we have
\begin{equation*}
\Bigl(\sideset{}{^\ast}\prod_{i\in[\alpha]}\Bigl(\mathrm{T}_{2^t}^{(+)}
-\sideset{}{^\ast}\prod_{a^i\in A_i}\widetilde{\boldsymbol{\mathfrak{a}}}(a^i)\Bigr)\Bigr)\cdot\overline{\mathbf{U}}(2^t)
=\mathrm{T}_{2^t}^{(+)}
-\boldsymbol{\gamma}(\mathcal{A}^{\triangledown})\cdot\overline{\mathbf{U}}(2^t)\; ,
\end{equation*}
that is,
\begin{equation*}
\boldsymbol{\gamma}(\mathcal{A}^{\triangledown})\cdot\overline{\mathbf{U}}(2^t)=
\mathrm{T}_{2^t}^{(+)}-
\Bigl(\sideset{}{^\ast}\prod_{i\in[\alpha]}\Bigl(\mathrm{T}_{2^t}^{(+)}
-\sideset{}{^\ast}\prod_{a^i\in A_i}\widetilde{\boldsymbol{\mathfrak{a}}}(a^i)\Bigr)\Bigr)\cdot\overline{\mathbf{U}}(2^t)
\; .
\end{equation*}

\begin{theorem} \label{th:5} If $\mathcal{A}:=\{A_1,\ldots,A_{\alpha}\}$ is a nontrivial clutter on the ground set~$E_t$, then we have{\rm:}
\begin{itemize}
\item[\rm(i)]
\begin{equation}
\label{eq:5}
\boldsymbol{\gamma}(\mathcal{A}^{\triangledown})=
\mathrm{T}_{2^t}^{(+)}-
\Bigl(\sideset{}{^\ast}\prod_{i\in[\alpha]}\Bigl(\mathrm{T}^{(+)}_{2^t}
-\sideset{}{^\ast}\prod_{a^i\in A_i}\widetilde{\boldsymbol{\mathfrak{a}}}(a^i)\Bigr)\Bigr)\; .
\end{equation}
\item[\rm(ii)]
\begin{equation}
\label{eq:4}
\boldsymbol{\gamma}(\mathfrak{B}(\mathcal{A})^{\triangledown})
=\Bigl(\sideset{}{^\ast}\prod_{i\in[\alpha]}\Bigl(\mathrm{T}^{(+)}_{2^t}
-\sideset{}{^\ast}\prod_{a^i\in A_i}\widetilde{\boldsymbol{\mathfrak{a}}}(a^i)\Bigr)\Bigr)\cdot\overline{\mathbf{U}}(2^t)
\; .
\end{equation}
\end{itemize}
\end{theorem}

\begin{example}\label{page:3} \label{th:17}
Suppose $t:=3$, and $E_t=\{1,2,3\}$.
We have in our hands the characteristic vectors
{\small
\begin{align*}
\widetilde{\boldsymbol{\mathfrak{a}}}(1):=\widetilde{\boldsymbol{\mathfrak{a}}}(1;2^t):=\boldsymbol{\gamma}(\{\{1\}\}^{\triangledown})
&=(0,\phantom{-}1,\phantom{-}0,\phantom{-}0,\phantom{-}1,\phantom{-}1,\phantom{-}0,\phantom{-}1)\in\{0,1\}^{2^t}\; ,\\
\widetilde{\boldsymbol{\mathfrak{a}}}(2):=\widetilde{\boldsymbol{\mathfrak{a}}}(2;2^t):=\boldsymbol{\gamma}(\{\{2\}\}^{\triangledown})
&=(0,\phantom{-}0,\phantom{-}1,\phantom{-}0,\phantom{-}1,\phantom{-}0,\phantom{-}1,\phantom{-}1)\; ,\\
\widetilde{\boldsymbol{\mathfrak{a}}}(3):=\widetilde{\boldsymbol{\mathfrak{a}}}(3;2^t):=\boldsymbol{\gamma}(\{\{3\}\}^{\triangledown})
&=(0,\phantom{-}0,\phantom{-}0,\phantom{-}1,\phantom{-}0,\phantom{-}1,\phantom{-}1,\phantom{-}1)\; ,
\end{align*}
}
$\quad$\vspace{-5mm}\\
associated with the principal increasing families that are generated by the clutters~$\{\{a\}\}$, for the elements $a\in E_t$ of the ground set.

We are given the clutter~$\mathcal{A}:=\{A_1,A_2\}$ on the ground set~$E_t$, where $A_1:=\{1,2\}$ and $A_2:=\{2,3\}$,
and we want to know the characteristic vector~$\boldsymbol{\gamma}(\mathfrak{B}(\mathcal{A})^{\triangledown})$ of the increasing family of blocking sets~$\mathfrak{B}(\mathcal{A})^{\triangledown}$ of the clutter~$\mathcal{A}$.

Turning to Theorem~{\rm\ref{th:5}(ii)}, we see that
{\small
\begin{align*}
\sideset{}{^\ast}\prod_{a^1\in A_1}\widetilde{\boldsymbol{\mathfrak{a}}}(a^1):=\sideset{}{^\ast}\prod_{a^1\in \{1,2\}}\widetilde{\boldsymbol{\mathfrak{a}}}(a^1)
=\quad &(0,\phantom{-}1,\phantom{-}0,\phantom{-}0,\phantom{-}1,\phantom{-}1,\phantom{-}0,\phantom{-}1)\\
\ast\, &(0,\phantom{-}0,\phantom{-}1,\phantom{-}0,\phantom{-}1,\phantom{-}0,\phantom{-}1,\phantom{-}1)\\
=\quad &(0,\phantom{-}0,\phantom{-}0,\phantom{-}0,\phantom{-}1,\phantom{-}0,\phantom{-}0,\phantom{-}1)\; ,
\\
\sideset{}{^\ast}\prod_{a^2\in A_2}\widetilde{\boldsymbol{\mathfrak{a}}}(a^2):=
\sideset{}{^\ast}\prod_{a^2\in\{2,3\}}\widetilde{\boldsymbol{\mathfrak{a}}}(a^2)
=\quad &(0,\phantom{-}0,\phantom{-}1,\phantom{-}0,\phantom{-}1,\phantom{-}0,\phantom{-}1,\phantom{-}1)\\
\ast\, &(0,\phantom{-}0,\phantom{-}0,\phantom{-}1,\phantom{-}0,\phantom{-}1,\phantom{-}1,\phantom{-}1)\\
=\quad&(0,\phantom{-}0,\phantom{-}0,\phantom{-}0,\phantom{-}0,\phantom{-}0,\phantom{-}1,\phantom{-}1)\; ;\\
\mathrm{T}^{(+)}_{2^t}
-\sideset{}{^\ast}\prod_{a^1\in A_1}\widetilde{\boldsymbol{\mathfrak{a}}}(a^1)
=\quad &(1,\phantom{-}1,\phantom{-}1,\phantom{-}1,\phantom{-}0,\phantom{-}1,\phantom{-}1,\phantom{-}0)\; ,\\
\mathrm{T}^{(+)}_{2^t}
-\sideset{}{^\ast}\prod_{a^2\in A_2}\widetilde{\boldsymbol{\mathfrak{a}}}(a^2)
=\quad &(1,\phantom{-}1,\phantom{-}1,\phantom{-}1,\phantom{-}1,\phantom{-}1,\phantom{-}0,\phantom{-}0)\; ;\\
\sideset{}{^\ast}\prod_{i\in[2]}\Bigl(\mathrm{T}^{(+)}_{2^t}
-\sideset{}{^\ast}\prod_{a^i\in A_i}\widetilde{\boldsymbol{\mathfrak{a}}}(a^i)\Bigr)
=\quad &(1,\phantom{-}1,\phantom{-}1,\phantom{-}1,\phantom{-}0,\phantom{-}1,\phantom{-}1,\phantom{-}0)\\
\ast\, &(1,\phantom{-}1,\phantom{-}1,\phantom{-}1,\phantom{-}1,\phantom{-}1,\phantom{-}0,\phantom{-}0)\\
=\quad &(1,\phantom{-}1,\phantom{-}1,\phantom{-}1,\phantom{-}0,\phantom{-}1,\phantom{-}0,\phantom{-}0)\; ,\\
\intertext{\normalsize{and finally}}
\boldsymbol{\gamma}(\mathfrak{B}(\mathcal{A})^{\triangledown})
=\Bigl(\sideset{}{^\ast}\prod_{i\in[2]}\Bigl(\mathrm{T}^{(+)}_{2^t}
-\sideset{}{^\ast}\prod_{a^i\in A_i}\widetilde{\boldsymbol{\mathfrak{a}}}(a^i)\Bigr)\Bigr)\cdot\overline{\mathbf{U}}(2^t)
=\quad &(0,\phantom{-}0,\phantom{-}1,\phantom{-}0,\phantom{-}1,\phantom{-}1,\phantom{-}1,\phantom{-}1)\; .
\end{align*}
}
In Example~{\rm\ref{th:7}} on page~{\rm\pageref{page:2}}, we will attempt to extract from the above vector~$\boldsymbol{\gamma}(\mathfrak{B}(\mathcal{A})^{\triangledown})$
the characteristic vector~$\boldsymbol{\gamma}(\mathfrak{B}(\mathcal{A}))$ of the blocker~$\mathfrak{B}(\mathcal{A})$.
\end{example}

\subsubsection{The characteristic vector of the subfamily of inclusion-minimal sets $\bmin\mathcal{F}$ in a family~$\mathcal{F}$} $\quad$

Suppose we are given the characteristic vector $\boldsymbol{\gamma}(\mathcal{F})$ of a
nonempty
family~$\mathcal{F}\subset\mathbf{2}^{[t]}$ of subsets of the ground set $E_t$, such that $\mathcal{F}\not\ni\hat{0}$. We can read off the position numbers of all the inclusion-minimal sets in the family~$\mathcal{F}$ in the following straightforward way (see~Example~\ref{th:7} on page~\pageref{page:2}):

\begin{algorithm}
\label{th:6}
$\quad$ {\em
\begin{itemize}
\item[\tt Input:] The char.-vector $\boldsymbol{\gamma}(\mathcal{F})$ of a family~$\mathcal{F}\subset\mathbf{2}^{[t]}$, such that
$\emptyset\neq\mathcal{F}\not\ni\hat{0}$.

\item[\tt Output:] \noindent{}A set~$M$ is the set $\supp(\boldsymbol{\gamma}(\bmin\mathcal{F}))$ of position numbers of the members of the clutter~$\bmin\mathcal{F}$;

\noindent{}a vector $\boldsymbol{\beta}$ is the char.-vector $\boldsymbol{\gamma}(\bmin\mathcal{F})$ of the clutter~$\bmin\mathcal{F}$
(\!{\em this data is optional});

\noindent{}a family $\mathcal{B}$ is the clutter~$\bmin\mathcal{F}$ (\!{\em this data is optional}).
\end{itemize}
\begin{itemize}
\item[\rm(0).] \textsf{Define}
$\boldsymbol{\phi}\in\{0,1\}^{2^t}$, and \textsf{store}\ \ $\boldsymbol{\phi}\, \leftarrow\, \boldsymbol{\gamma}(\mathcal{F})$;

\noindent{}\textsf{define} $\boldsymbol{\beta}\in\{0,1\}^{2^t}$, and \textsf{store} $\boldsymbol{\beta}\, \leftarrow\, (0,\ldots,0)$;\hspace{0.2cm}\text{\% {\em this action is optional}}

\noindent{}\textsf{define} $\mathcal{B}\subset\mathbf{2}^{[t]}$, and \textsf{store} $\mathcal{B}\, \leftarrow\, \emptyset$;\hspace{2.25cm}\text{\% {\em this action is optional}}

\noindent{}\textsf{define} $M\subset[2^t]$, and \textsf{store} $M\, \leftarrow\, \hat{0}$;

\noindent{}\textsf{define} $m\in\mathbb{N}$, and \textsf{store} $m\, \leftarrow\, 0$;

\noindent{}\textsf{define} $B\in\mathbf{2}^{[t]}$,  and \textsf{store} $B\, \leftarrow\, \hat{0}$.

\item[\rm(1).] \textsf{If} $|\supp(\boldsymbol{\phi})|=0$, \textsf{then} \textsf{go to} Step~{\rm(3)},

\noindent{}\textsf{else} \textsf{go to} Step~{\rm(2)}.

\item[\rm(2).] \textsf{Store} $m\, \leftarrow\, \min\supp(\boldsymbol{\phi})$,

and \textsf{store} $M\, \leftarrow\, M\,\dot{\cup}\,\{m\}$,

and \textsf{store} $B\, \leftarrow\, \varGamma(m)$,

and \textsf{store} $\mathcal{B}\, \leftarrow\, \mathcal{B}\,\,\dot{\cup}\, \{B\}$;\hspace{3cm}\text{\% {\em this action is optional}}

\noindent{}\textsf{store} $\boldsymbol{\beta}\, \leftarrow\, \boldsymbol{\beta}+\boldsymbol{\sigma}(m)$;\hspace{4.05cm}\text{\% {\em this action is optional}}

\noindent{}\textsf{If} $|\supp(\boldsymbol{\phi})|=1$, \textsf{then go to} Step~(3),

\noindent{}\textsf{else store} $\boldsymbol{\phi}\, \leftarrow\,
\boldsymbol{\phi}\,-\,\boldsymbol{\phi}\,\ast\underbrace{
\sideset{}{^\ast}\prod_{e\in B}\widetilde{\boldsymbol{\mathfrak{a}}}(e)
}_{
\boldsymbol{\gamma}(\{B\}^{\triangledown})
}\
$.

\noindent{}\textsf{Go to} Step~{\rm(1)}.

\item[\rm(3).] \textsf{Stop}.
\end{itemize}
}
\end{algorithm}

\subsubsection{More on the blocker $\mathfrak{B}(\mathcal{A})$} $\quad$

If we know (see, e.g.,~Theorem~\ref{th:5}(ii)) the characteristic vector~$\boldsymbol{\gamma}(\mathcal{F})$ of the increasing family of  blocking sets~$\mathcal{F}:=\mathfrak{B}(\mathcal{A})^{\triangledown}$ of a clutter~\mbox{$\mathcal{A}:=\{A_1,\ldots,A_{\alpha}\}$} on the ground set $E_t$, then a description of the blocker $\bmin\mathcal{F}:=\mathfrak{B}(\mathcal{A})$ can be obtained by an application of Algorithm~\ref{th:6} to the vector~$\boldsymbol{\gamma}(\mathcal{F})$; see~Example~\ref{th:7}.

\begin{example}\label{page:2}
\label{th:7}
Suppose $t:=3$, and $E_t=\{1,2,3\}$. Note that
\begin{align*}
\widetilde{\boldsymbol{\mathfrak{a}}}(2):=\widetilde{\boldsymbol{\mathfrak{a}}}(2;2^t)
:=\boldsymbol{\gamma}(\{\{2\}\}^{\triangledown})
&=(0,0,1,0,1,0,1,1)\in\{0,1\}^{2^t}\; .\\
\intertext{We are given the characteristic vector}
\boldsymbol{\gamma}(\mathcal{F}):\!&=(0,0,1,0,1,1,1,1)\in\{0,1\}^{2^t}
\end{align*}
of the increasing family of blocking sets~$\mathcal{F}:=\mathfrak{B}(\mathcal{A})^{\triangledown}$ of the clutter~$\mathcal{A}:=\{\{1,2\},\{2,3\}\}$
on the ground set~$E_t$; see, e.g.,~Example~{\rm\ref{th:17}} on page~{\rm\pageref{page:3}}. In order
to find a description of the clutter~$\bmin\mathcal{F}:=\mathfrak{B}(\mathcal{A})$, let us apply~Algorithm~{\rm\ref{th:6}} to
the vector~$\boldsymbol{\gamma}(\mathcal{F})${\rm:}
{\small
\begin{align*}
\boldsymbol{\phi}\ &\leftarrow\ \boldsymbol{\gamma}(\mathcal{F}):=(0,0,1,0,1,1,1,1)\; ;\\
|\supp(\boldsymbol{\phi})|\ &>\ 0\; ;\\
m\ &\leftarrow\ \underbrace{\min\supp((0,0,1,0,1,1,1,1))}_{3}\; ,\\
M\ \leftarrow\ \underbrace{\hat{0}\;\dot\cup\;\{3\}}_{\{3\}}\; ,\ \ \ \ \
B\ \leftarrow\ \underbrace{\varGamma(3)}_{\{2\}}\; ,\ \ \ \ \
\mathcal{B}\ &\leftarrow\  \underbrace{\emptyset\;\dot\cup\; \{2\}}_{\{\{2\}\}}\; ;\\
\boldsymbol{\beta}\ &\leftarrow\ \underbrace{(0,0,0,0,0,0,0,0)+\boldsymbol{\sigma}(3)}_{
(0,0,1,0,0,0,0,0)
}\; ;\\
\boldsymbol{\phi}\ &\leftarrow\ \underbrace{(0,0,1,0,1,1,1,1)-(0,0,1,0,1,1,1,1)\ast
\widetilde{\boldsymbol{\mathfrak{a}}}(2)}_{
(0,0,0,0,0,1,0,0)
}\; ;\\
|\supp(\boldsymbol{\phi})|\ &>\ 0\; ;\\
m\ &\leftarrow\ \underbrace{\min\supp((0,0,0,0,0,1,0,0))}_{6}\; ,\\
M\ \leftarrow\ \underbrace{\{3\}\;\dot\cup\;\{6\}}_{\{3,6\}}\; ,\ \ \ \ \
B\ \leftarrow\ \underbrace{\varGamma(6)}_{\{1,3\}}\; ,\ \ \ \ \
\mathcal{B}\ &\leftarrow\ \underbrace{\{\{2\}\}\;\dot\cup\; \{1,3\}}_{\{\{2\},\{1,3\}\}}\; ;\\
\boldsymbol{\beta}\ &\leftarrow\ \underbrace{(0,0,1,0,0,0,0,0)+\boldsymbol{\sigma}(6)}_{
(0,0,1,0,0,1,0,0)
}\; ;\\
|\supp(\boldsymbol{\phi})|\ &=\ 1\; ;\\
\text{\em\textsf{Stop.}} & \quad
\end{align*}
}

We see that the set~$\supp(\boldsymbol{\gamma}(\bmin\mathcal{F}))=:M$ of the position numbers of the members of the blocker~$\mathfrak{B}(\mathcal{A})=:\bmin\mathcal{F}$ is the set $\{3,6\}$.

The characteristic vector $\boldsymbol{\gamma}(\bmin\mathcal{F})=:\boldsymbol{\beta}$ of the blocker~$\mathfrak{B}(\mathcal{A})=:\bmin\mathcal{F}$ is the vector~$(0,0,1,0,0,1,0,0)$.

The blocker~$\mathfrak{B}(\mathcal{A})=:\bmin\mathcal{F}=:\mathcal{B}$ of the clutter~$\mathcal{A}:=\{\{1,2\},\{2,3\}\}$ is the clutter~$\{\{2\},\{1,3\}\}$.
\end{example}

\part*{Blocking\,/\,Voting}

\section{Decompositions of the characteristic topes and of the characteristic vectors of families}

\noindent$\bullet$ The vertices $R^i\in\{1,-1\}^t$ of the symmetric cycle~$\boldsymbol{R}$ in the hypercube graph~$\boldsymbol{H}(t,2)$, given in~(\ref{eq:12})(\ref{eq:13}), are just simply defined and useful decomposition components of topes of the oriented matroid~$\mathcal{H}:=(E_t,\{1,-1\}^t)$.

In the context of the combinatorics of finite sets, the vertices $R^i\in\{1,-1\}^{2^t}$ of a distinguished {\em symmetric cycle\/}
\begin{equation*}
\boldsymbol{R}:=(R^0,R^1,\ldots,R^{2\cdot 2^t-1},R^0)
\end{equation*}
in the hypercube graph of topes~$\boldsymbol{H}(2^t,2)$ of the oriented matroid~$\mathcal{H}_{2^t}:=(E_{2^t},\{1,-1\}^{2^t})$, where
\begin{equation}
\label{eq:66}
\begin{split}
R^0:\!&=\mathrm{T}_{2^t}^{(+)}\; ,\\
R^s:\!&={}_{-[s]}R^0\; ,\ \ \ 1\leq s\leq 2^t-1\; ,\\
\end{split}
\end{equation}
and
\begin{equation}
\label{eq:67}
R^{2^t+k}:=-R^k\; ,\ \ \ 0\leq k\leq 2^t-1\; ,
\end{equation}
have an additional meaning:

\begin{remark} Let $\boldsymbol{R}$ be the symmetric cycle in the tope graph of the oriented matroid~$\mathcal{H}_{2^t}:=(E_{2^t},\{1,-1\}^{2^t})$, defined by~{\rm(\ref{eq:66})(\ref{eq:67})}.
\begin{itemize}
\item[\rm(i)] The vertex $R^0:=\mathrm{T}_{2^t}^{(+)}\in\mathrm{V}(\boldsymbol{R})$ is the {\em characteristic tope}
$T_{\emptyset}$ of the {\em empty family} $\emptyset$ on the ground set~$E_t$.

\noindent{}The vertex $R^{2^t}:=-\mathrm{T}_{2^t}^{(+)}\in\mathrm{V}(\boldsymbol{R})$ is the {\em characteristic tope}~$T_{\mathbf{2}^{[t]}}$ of the {\em power set\/} $\mathbf{2}^{[t]}$ of the set~$E_t$.

\item[\rm(ii)] If $1\leq i\leq 2^t-1$, then the vertex $R^i\in\mathrm{V}(\boldsymbol{R})$ is the {\em characteristic tope} $T_{\mathcal{F}}$ of a {\em decreasing family} $\mathcal{F}$ of subsets of the ground set~$E_t$. In other words, the family $\mathcal{F}$ is a particular abstract simplicial complex, when $1<i\leq 2^t-1$.

    Either the subfamily $\bmax\mathcal{F}$ is an $s$-uniform clutter, where $s:=|\varGamma(\max (T_{\mathcal{F}})^{-})|$, or we have $\{|F|\colon F\in\bmax\mathcal{F}\}=\{s,s-1\}$. Indeed, we have
    \begin{equation*}
    \bmax\mathcal{F}=\underbrace{(\mathcal{F}\cap\tbinom{E_t}{s})}_{(\bmax\mathcal{F})\,\cap\,\binom{E_t}{s}}\; \dot\cup\ \;
    \underbrace{\Bigl(\tbinom{E_t}{s-1}
    -(\mathcal{F}\cap\tbinom{E_t}{s})^{\vartriangle}\Bigr)}_{(\bmax\mathcal{F})\,\cap\,\binom{E_t}{s-1}}\; .
    \end{equation*}

\item[\rm(iii)] If $2^t+1\leq i\leq 2\cdot 2^t-1$, then the vertex $R^i\in\mathrm{V}(\boldsymbol{R})$ is the {\em characteristic tope} $T_{\mathcal{F}}$ of an {\em increasing family} $\mathcal{F}$ of subsets of the ground set~$E_t$.

    Either the subfamily $\bmin\mathcal{F}$ is an $s$-uniform clutter, where $s:=|\varGamma(\min (T_{\mathcal{F}})^{-})|$, or we have $\{|F|\colon F\in\bmin\mathcal{F}\}=\{s,s+1\}$. We have
    \begin{equation*}
    \bmin\mathcal{F}=\underbrace{(\mathcal{F}\cap\tbinom{E_t}{s})}_{(\bmin\mathcal{F})\,\cap\,\binom{E_t}{s}}\; \dot\cup\ \;
    \underbrace{\Bigl(\tbinom{E_t}{s+1}
    -(\mathcal{F}\cap\tbinom{E_t}{s})^{\triangledown}\Bigr)}_{(\bmin\mathcal{F})\,\cap\,\binom{E_t}{s+1}}\; .
    \end{equation*}

    If $i=3\cdot 2^{t-1}$, then the clutter~$\bmin\mathcal{F}$ is {\em self-dual}.
\end{itemize}
\end{remark}

\noindent$\bullet$ A distinguished symmetric cycle~$\widetilde{\boldsymbol{R}}:=(\widetilde{R}^0,\widetilde{R}^1,\ldots,\widetilde{R}^{2\cdot 2^t-1},\widetilde{R}^0)$ in the hypercube graph~$\widetilde{\boldsymbol{H}}(2^t,2)$ on the vertex set $\{0,1\}^{2^t}$ is defined\footnote{$\quad$Here $\boldsymbol{\sigma}(e)$ is the $e$th standard unit vector of the space~$\mathbb{R}^{2^t}$.} as follows:
\begin{equation*}
\begin{split}
\widetilde{R}^0:\!&=(0,\ldots,0)\; ,\\
\widetilde{R}^s:\!&=\sum\nolimits_{e\in[s]}\boldsymbol{\sigma}(e)\; ,\ \ \ 1\leq s\leq 2^t-1\; ,
\end{split}
\end{equation*}
and
\begin{equation*}
\widetilde{R}^{2^t+k}:=\mathrm{T}_{2^t}^{(+)}-\widetilde{R}^k\; ,\ \ \ 0\leq k\leq 2^t-1\; .
\end{equation*}

We let~$\mathrm{V}(\widetilde{\boldsymbol{R}}):=(\widetilde{R}^0,\widetilde{R}^1,\ldots,\widetilde{R}^{2\cdot 2^t-1})$ denote the vertex sequence of the cycle~$\widetilde{\boldsymbol{R}}$.

\noindent$\bullet$ Let $\mathcal{F}\subset\mathbf{2}^{[t]}$ be a family of subsets of the ground set~$E_t$, $\emptyset\neq\mathcal{F}\not\ni\hat{0}$.
As earlier, we associate with the family $\mathcal{F}$ its characteristic tope $T_{\mathcal{F}}\in\{1,-1\}^{2^t}$, defined by~(\ref{eq:68}).

Recall that there exists a unique inclusion-minimal subset
\begin{equation*}
\boldsymbol{Q}(T_{\mathcal{F}},\boldsymbol{R})\subset\mathrm{V}(\boldsymbol{R}):=(R^0,R^1,\ldots,R^{2\cdot 2^t-1})
\end{equation*}
of the vertex sequence $\mathrm{V}(\boldsymbol{R})$ of the cycle~$\boldsymbol{R}$, defined by~(\ref{eq:66})(\ref{eq:67}), such that
\begin{equation*}
T_{\mathcal{F}}=\sum_{Q\in\boldsymbol{Q}(T_{\mathcal{F}},\boldsymbol{R})} Q\; .
\end{equation*}
In other words, there exists a unique row vector~$\boldsymbol{x}:=\boldsymbol{x}(T_{\mathcal{F}}):=\boldsymbol{x}(T_{\mathcal{F}},\boldsymbol{R}):=(x_1,\ldots,x_{2^t})
\in\{-1,0,1\}^{2^t}$, such that
\begin{equation}
\label{eq:70}
T_{\mathcal{F}}=\sum_{i\in[2^t]}x_i\cdot R^{i-1}=\boldsymbol{x}\mathbf{M}\; ,
\end{equation}
where
\begin{equation}
\label{eq:73}
\mathbf{M}:=\mathbf{M}(\boldsymbol{R}):=
\left(
\begin{smallmatrix}
R^0\\
R^1 \vspace{-1.5mm}\\
\vdots \vspace{0.5mm}\\
R^{2^t-1}
\end{smallmatrix}
\right)\; .
\end{equation}
One can see this matrix (in the case $t:=3$) in Example~\ref{th:18} on page~\pageref{page:4}. Thus, we have
\begin{equation*}
\boldsymbol{x}=T_{\mathcal{F}}\cdot\mathbf{M}^{-1}\; ,
\end{equation*}
and
\begin{equation*}
\boldsymbol{Q}(T_{\mathcal{F}},\boldsymbol{R}):=\{x_i\cdot R^{i-1}\colon x_i\neq 0\}\; .
\end{equation*}
We use the notation $\mathfrak{q}(T_{\mathcal{F}}):=\mathfrak{q}(T_{\mathcal{F}},\boldsymbol{R}):=|\boldsymbol{Q}(T_{\mathcal{F}},\boldsymbol{R})|$ to denote the cardinality of the set~$\boldsymbol{Q}(T_{\mathcal{F}},\boldsymbol{R})$.

\noindent$\bullet$ Let us consider the subset
\begin{equation*}
\widetilde{\boldsymbol{Q}}(\underbrace{
\boldsymbol{\gamma}(\mathcal{F})
}_{\frac{1}{2}(\mathrm{T}_{2^t}^{(+)}-T_{\mathcal{F}})},
\widetilde{\boldsymbol{R}}):=
\{\tfrac{1}{2}(\mathrm{T}_{2^t}^{(+)}-Q)\colon \; Q\in\boldsymbol{Q}(T_{\mathcal{F}},\boldsymbol{R})\}\subset\mathrm{V}(\widetilde{\boldsymbol{R}})
\; ,
\end{equation*}
and let us use the notation $\mathfrak{q}(\boldsymbol{\gamma}(\mathcal{F})):=\mathfrak{q}(\boldsymbol{\gamma}(\mathcal{F}),\widetilde{\boldsymbol{R}})
:=|\widetilde{\boldsymbol{Q}}(\boldsymbol{\gamma}(\mathcal{F}),
\widetilde{\boldsymbol{R}})|=\mathfrak{q}(T_{\mathcal{F}})$ to denote its cardinality.

In analogy with~(\ref{eq:69}), we have
\begin{equation}
\label{eq:31}
\boldsymbol{\gamma}(\mathcal{F})=-\tfrac{1}{2}\bigl(\mathfrak{q}(\boldsymbol{\gamma}(\mathcal{F}))
-1\bigr)\cdot\mathrm{T}_{2^t}^{(+)} +
\sum_{
\substack{\widetilde{Q}\in\widetilde{\boldsymbol{Q}}(\boldsymbol{\gamma}(\mathcal{F}),\widetilde{\boldsymbol{R}})\colon\\
\widetilde{Q}\neq
(0,\ldots,0)=:\widetilde{R}^0
}
}
\widetilde{Q}\; .
\end{equation}

\noindent$\bullet$ Let $\mathcal{A}\subset\mathbf{2}^{[t]}$ be a nontrivial clutter on the ground set~$E_t$, and let $\mathcal{B}:=\mathfrak{B}(\mathcal{A})$ be its blocker.
We associate with the families $\mathcal{A}^{\triangledown}$, $\mathcal{B}^{\triangledown}$, $\mathcal{A}$ and $\mathcal{B}$ their characteristic topes $T_{\mathcal{A}^{\triangledown}}$, $T_{\mathcal{B}^{\triangledown}}$, $T_{\mathcal{A}}$, $T_{\mathcal{B}}\in\{1,-1\}^{2^t}$, and their characteristic vectors $\boldsymbol{\gamma}(\mathcal{A}^{\triangledown})$, $\boldsymbol{\gamma}(\mathcal{B}^{\triangledown})$, $\boldsymbol{\gamma}(\mathcal{A})$, $\boldsymbol{\gamma}(\mathcal{B})\in\{0,1\}^{2^t}$. See~(\ref{eq:74})--(\ref{eq:81}) in~Example~\ref{th:9}.

\begin{example} \label{th:9}
Suppose $t:=3$, and $E_t=\{1,2,3\}$. Let $\boldsymbol{R}$ by the symmetric cycle in the hypercube graph~$\boldsymbol{H}(2^t,2)$ on the vertex set $\{1,-1\}^{2^t}$, defined by~{\rm(\ref{eq:66})(\ref{eq:67})}.

We are given the {\em blocking pair}\! of {\em clutters\/} $\mathcal{A}:=\{\{1,2\},\{2,3\}\}$ and  $\mathcal{B}:=\mathfrak{B}(\mathcal{A})=\{\{1,3\},\{2\}\}$ on the ground set~$E_t$.

The families $\mathcal{A}^{\triangledown}$, $\mathcal{B}^{\triangledown}$, $\mathcal{A}$ and $\mathcal{B}$ are described by their
characteristic topes
{\small
\begin{align}
\label{eq:74}
T_{\mathcal{A}^{\triangledown}}:\!&=(1,\phantom{+}1,\phantom{+}1,\phantom{+}1,-1,\phantom{+}1,-1,-1)\in\{1,-1\}^{2^t}\; ,\\
\label{eq:75}
T_{\mathcal{B}^{\triangledown}}:\!&=(1,\phantom{+}1,-1,\phantom{+}1,-1,-1,-1,-1)\; ,\\
\label{eq:76}
T_{\mathcal{A}}:\!&=(1,\phantom{+}1,\phantom{+}1,\phantom{+}1,-1,\phantom{+}1,-1,\phantom{+}1)\; ,\\
\label{eq:77}
T_{\mathcal{B}}:\!&=(1,\phantom{+}1,-1,\phantom{+}1,\phantom{+}1,-1,\phantom{+}1,\phantom{+}1)\; ,\\
\intertext{\normalsize{and by their characteristic vectors}}
\label{eq:78}
\boldsymbol{\gamma}(\mathcal{A}^{\triangledown}):\!&=(0,\phantom{+}0,\phantom{+}0,\phantom{+}0,\phantom{+}1,\phantom{+}0,\phantom{+}1,\phantom{+}1)\in\{0,1\}^{2^t}\; ,\\
\label{eq:79}
\boldsymbol{\gamma}(\mathcal{B}^{\triangledown})
:\!&=(0,\phantom{+}0,\phantom{+}1,\phantom{+}0,\phantom{+}1,\phantom{+}1,\phantom{+}1,\phantom{+}1)\; ,\\
\label{eq:80}
\boldsymbol{\gamma}(\mathcal{A}):\!&=(0,\phantom{+}0,\phantom{+}0,\phantom{+}0,\phantom{+}1,\phantom{+}0,\phantom{+}1,\phantom{+}0)\; ,\\
\label{eq:81}
\boldsymbol{\gamma}(\mathcal{B}):\!&=(0,\phantom{+}0,\phantom{+}1,\phantom{+}0,\phantom{+}0,\phantom{+}1,\phantom{+}0,\phantom{+}0)\; .
\end{align}
}

Turning to decompositions of the form~{\rm(\ref{eq:70})}, we see that
{\small
\begin{align}
\label{eq:82}
\boldsymbol{x}(T_{\mathcal{A}^{\triangledown}}
)&=(0,\phantom{-}0,\phantom{-}0,\phantom{-}0,-1,\phantom{-}1,-1,\phantom{-}0)\in\{-1,0,1\}^{2^t}\; ,\\
\label{eq:83}
\boldsymbol{x}(T_{\mathcal{B}^{\triangledown}}
)&=(0,\phantom{-}0,-1,\phantom{-}1,-1,\phantom{-}0,\phantom{-}0,\phantom{-}0)
\; ,\\
\label{eq:84}
\boldsymbol{x}(T_{\mathcal{A}}
)&=(1,\phantom{-}0,\phantom{-}0,\phantom{-}0,-1,\phantom{-}1,-1,\phantom{-}1)
\; ,\\
\label{eq:85}
\boldsymbol{x}(T_{\mathcal{B}}
)&=(1,\phantom{-}0,-1,\phantom{-}1,\phantom{-}0,-1,\phantom{-}1,\phantom{-}0)
\; .
\end{align}
}
Thus, we have the decompositions:
{\small
\begin{equation*}
\begin{split}
T_{\mathcal{A}^{\triangledown}}:=\phantom{-}&(\phantom{+}1,\phantom{+}1,\phantom{+}1,\phantom{+}1,-1,\phantom{+}1,-1,-1)
=-\underbrace{R^4}_{-R^{12}}+\;R^5-\underbrace{R^6}_{-R^{14}}\\
=-\,&(-1,-1,-1,-1,\phantom{+}1,\phantom{+}1,\phantom{+}1,\phantom{+}1)\\
\phantom{=}+\,&(-1,-1,-1,-1,-1,\phantom{+}1,\phantom{+}1,\phantom{+}1)\\
\phantom{=}-\,&(-1,-1,-1,-1,-1,-1,\phantom{+}1,\phantom{+}1)=R^5+R^{12}+R^{14}\\
=\phantom{=}&(-1,-1,-1,-1,-1,\phantom{+}1,\phantom{+}1,\phantom{+}1)\\
\phantom{=}+\,&(\phantom{+}1,\phantom{+}1,\phantom{+}1,\phantom{+}1,-1,-1,-1,-1)\\
\phantom{=}+\,&(\phantom{+}1,\phantom{+}1,\phantom{+}1,\phantom{+}1,\phantom{+}1,\phantom{+}1,-1,-1)\; ,
\end{split}
\end{equation*}
\begin{equation*}
\begin{split}
T_{\mathcal{B}^{\triangledown}}:=\phantom{-}&(\phantom{+}1,\phantom{+}1,-1,\phantom{+}1,-1,-1,-1,-1)
=-\underbrace{R^2}_{-R^{10}}+\;R^3-\underbrace{R^4}_{-R^{12}}\\
=-\,&(-1,-1,\phantom{+}1,\phantom{+}1,\phantom{+}1,\phantom{+}1,\phantom{+}1,\phantom{+}1)\\
\phantom{=}+\,&(-1,-1,-1,\phantom{+}1,\phantom{+}1,\phantom{+}1,\phantom{+}1,\phantom{+}1)\\
\phantom{=}-\,&(-1,-1,-1,-1,\phantom{+}1,\phantom{+}1,\phantom{+}1,\phantom{+}1)=R^3+R^{10}+R^{12}\\
=\phantom{=}&(-1,-1,-1,\phantom{+}1,\phantom{+}1,\phantom{+}1,\phantom{+}1,\phantom{+}1)\\
\phantom{=}+\,&(\phantom{+}1,\phantom{+}1,-1,-1,-1,-1,-1,-1)\\
\phantom{=}+\,&(\phantom{+}1,\phantom{+}1,\phantom{+}1,\phantom{+}1,-1,-1,-1,-1)\; ,
\end{split}
\end{equation*}
\begin{equation*}
\begin{split}
T_{\mathcal{A}}:=\phantom{-}\,&(\phantom{+}1,\phantom{+}1,\phantom{+}1,\phantom{+}1,-1,\phantom{+}1,-1,\phantom{+}1)
=\underbrace{R^0}_{\mathrm{T}_{2^t}^{(+)}}-\underbrace{R^4}_{-R^{12}}+\;R^5-\underbrace{R^6}_{-R^{14}}+\;R^7\\
=\phantom{+}&(\phantom{+}1,\phantom{+}1,\phantom{+}1,\phantom{+}1,\phantom{+}1,\phantom{+}1,\phantom{+}1,\phantom{+}1)\\
\phantom{=}-\,&(-1,-1,-1,-1,\phantom{+}1,\phantom{+}1,\phantom{+}1,\phantom{+}1)\\
\phantom{=}+\,&(-1,-1,-1,-1,-1,\phantom{+}1,\phantom{+}1,\phantom{+}1)\\
\phantom{=}-\,&(-1,-1,-1,-1,-1,-1,\phantom{+}1,\phantom{+}1)\\
\phantom{=}+\,&(-1,-1,-1,-1,-1,-1,-1,\phantom{+}1)=\underbrace{R^0}_{\mathrm{T}_{2^t}^{(+)}}+R^5+R^7+R^{12}+R^{14}\\
=\phantom{+}&(\phantom{+}1,\phantom{+}1,\phantom{+}1,\phantom{+}1,\phantom{+}1,\phantom{+}1,\phantom{+}1,\phantom{+}1)\\
\phantom{=}+\,&(-1,-1,-1,-1,-1,\phantom{+}1,\phantom{+}1,\phantom{+}1)\\
\phantom{=}+\,&(-1,-1,-1,-1,-1,-1,-1,\phantom{+}1)\\
\phantom{=}+\,&(\phantom{+}1,\phantom{+}1,\phantom{+}1,\phantom{+}1,-1,-1,-1,-1)\\
\phantom{=}+\,&(\phantom{+}1,\phantom{+}1,\phantom{+}1,\phantom{+}1,\phantom{+}1,\phantom{+}1,-1,-1)\; ,
\end{split}
\end{equation*} }
\noindent{}and{\small
\begin{equation*}
\begin{split}
T_{\mathcal{B}}:=\phantom{-}&(\phantom{+}1,\phantom{+}1,-1,\phantom{+}1,\phantom{+}1,-1,\phantom{+}1,\phantom{+}1)
=\underbrace{R^0}_{\mathrm{T}_{2^t}^{(+)}}-\underbrace{R^2}_{-R^{10}}+\;R^3-\underbrace{R^5}_{-R^{13}}+\;R^6\\
=\phantom{+}&(\phantom{+}1,\phantom{+}1,\phantom{+}1,\phantom{+}1,\phantom{+}1,\phantom{+}1,\phantom{+}1,\phantom{+}1)\\
\phantom{=}-\,&(-1,-1,\phantom{+}1,\phantom{+}1,\phantom{+}1,\phantom{+}1,\phantom{+}1,\phantom{+}1)\\
\phantom{=}+\,&(-1,-1,-1,\phantom{+}1,\phantom{+}1,\phantom{+}1,\phantom{+}1,\phantom{+}1)\\
\phantom{=}-\,&(-1,-1,-1,-1,-1,\phantom{+}1,\phantom{+}1,\phantom{+}1)\\
\phantom{=}+\,&(-1,-1,-1,-1,-1,-1,\phantom{+}1,\phantom{+}1)=\underbrace{R^0}_{\mathrm{T}_{2^t}^{(+)}}+R^3+R^6+R^{10}+R^{13}\\
=\phantom{+}&(\phantom{+}1,\phantom{+}1,\phantom{+}1,\phantom{+}1,\phantom{+}1,\phantom{+}1,\phantom{+}1,\phantom{+}1)\\
\phantom{=}+\,&(-1,-1,-1,\phantom{+}1,\phantom{+}1,\phantom{+}1,\phantom{+}1,\phantom{+}1)\\
\phantom{=}+\,&(-1,-1,-1,-1,-1,-1,\phantom{+}1,\phantom{+}1)\\
\phantom{=}+\,&(\phantom{+}1,\phantom{+}1,-1,-1,-1,-1,-1,-1)\\
\phantom{=}+\,&(\phantom{+}1,\phantom{+}1,\phantom{+}1,\phantom{+}1,\phantom{+}1,-1,-1,-1)\; .
\end{split}
\end{equation*}
}

Relations of the form~{\rm(\ref{eq:31})} imply that
{\small
\begin{equation*}
\begin{split}
\boldsymbol{\gamma}(\mathcal{A}^{\triangledown}):=\phantom{+}
&(\phantom{+}0,\phantom{+}0,\phantom{+}0,\phantom{+}0,\phantom{+}1,\phantom{+}0,\phantom{+}1,\phantom{+}1)
=-\mathrm{T}_{2^t}^{(+)}+\widetilde{R}^{5}+\widetilde{R}^{12}+\widetilde{R}^{14}\\
=\phantom{+}\,&(-1,-1,-1,-1,-1,-1,-1,-1)\\
\phantom{=}+\,&(\phantom{+}1,\phantom{+}1,\phantom{+}1,\phantom{+}1,\phantom{+}1,\phantom{+}0,\phantom{+}0,\phantom{+}0)\\
\phantom{=}+\,&(\phantom{+}0,\phantom{+}0,\phantom{+}0,\phantom{+}0,\phantom{+}1,\phantom{+}1,\phantom{+}1,\phantom{+}1)\\
\phantom{=}+\,&(\phantom{+}0,\phantom{+}0,\phantom{+}0,\phantom{+}0,\phantom{+}0,\phantom{+}0,\phantom{+}1,\phantom{+}1)\; ,
\end{split}
\end{equation*}
\begin{equation*}
\begin{split}
\boldsymbol{\gamma}(\mathcal{B}^{\triangledown})
:=\phantom{-}&(\phantom{+}0,\phantom{+}0,\phantom{+}1,\phantom{+}0,\phantom{+}1,\phantom{+}1,\phantom{+}1,\phantom{+}1)
=-\mathrm{T}_{2^t}^{(+)}+\widetilde{R}^3+\widetilde{R}^{10}+\widetilde{R}^{12}\\
=\phantom{+}&(-1,-1,-1,-1,-1,-1,-1,-1)\\
\phantom{=}+\,&(\phantom{+}1,\phantom{+}1,\phantom{+}1,\phantom{+}0,\phantom{+}0,\phantom{+}0,\phantom{+}0,\phantom{+}0)\\
\phantom{=}+\,&(\phantom{+}0,\phantom{+}0,\phantom{+}1,\phantom{+}1,\phantom{+}1,\phantom{+}1,\phantom{+}1,\phantom{+}1)\\
\phantom{=}+\,&(\phantom{+}0,\phantom{+}0,\phantom{+}0,\phantom{+}0,\phantom{+}1,\phantom{+}1,\phantom{+}1,\phantom{+}1)\; ,
\end{split}
\end{equation*}
\begin{equation*}
\begin{split}
\boldsymbol{\gamma}(\mathcal{A}):=\phantom{-}&(\phantom{+}0,\phantom{+}0,\phantom{+}0,\phantom{+}0,\phantom{+}1,\phantom{+}0,\phantom{+}1,\phantom{+}0)
=-2\mathrm{T}_{2^t}^{(+)}+\underbrace{\widetilde{R}^0}_{(0,\ldots,0)}+\widetilde{R}^5+\widetilde{R}^7+\widetilde{R}^{12}+\widetilde{R}^{14}\\
=-\,&2\mathrm{T}_{2^t}^{(+)}+\widetilde{R}^5+\widetilde{R}^7+\widetilde{R}^{12}+\widetilde{R}^{14}\\
=\phantom{-\,}&(-2,-2,-2,-2,-2,-2,-2,-2)\\
\phantom{=}+\,&(\phantom{+}1,\phantom{+}1,\phantom{+}1,\phantom{+}1,\phantom{+}1,\phantom{+}0,\phantom{+}0,\phantom{+}0)\\
\phantom{=}+\,&(\phantom{+}1,\phantom{+}1,\phantom{+}1,\phantom{+}1,\phantom{+}1,\phantom{+}1,\phantom{+}1,\phantom{+}0)\\
\phantom{=}+\,&(\phantom{+}0,\phantom{+}0,\phantom{+}0,\phantom{+}0,\phantom{+}1,\phantom{+}1,\phantom{+}1,\phantom{+}1)\\
\phantom{=}+\,&(\phantom{+}0,\phantom{+}0,\phantom{+}0,\phantom{+}0,\phantom{+}0,\phantom{+}0,\phantom{+}1,\phantom{+}1)\; ,
\end{split}
\end{equation*} }
\noindent{}and{\small
\begin{equation*}
\begin{split}
\boldsymbol{\gamma}(\mathcal{B}):=\phantom{-}&(\phantom{+}0,\phantom{+}0,\phantom{+}1,\phantom{+}0,\phantom{+}0,\phantom{+}1,\phantom{+}0,\phantom{+}0)
=-2\mathrm{T}_{2^t}^{(+)}+\underbrace{\widetilde{R}^0}_{(0,\ldots,0)}+\widetilde{R}^3+\widetilde{R}^6+\widetilde{R}^{10}+\widetilde{R}^{13}\\
=-\,&2\mathrm{T}_{2^t}^{(+)}+\widetilde{R}^3+\widetilde{R}^6+\widetilde{R}^{10}+\widetilde{R}^{13}\\
=\phantom{+\,}&(-2,-2,-2,-2,-2,-2,-2,-2)\\
\phantom{=}+\,&(\phantom{+}1,\phantom{+}1,\phantom{+}1,\phantom{+}0,\phantom{+}0,\phantom{+}0,\phantom{+}0,\phantom{+}0)\\
\phantom{=}+\,&(\phantom{+}1,\phantom{+}1,\phantom{+}1,\phantom{+}1,\phantom{+}1,\phantom{+}1,\phantom{+}0,\phantom{+}0)\\
\phantom{=}+\,&(\phantom{+}0,\phantom{+}0,\phantom{+}1,\phantom{+}1,\phantom{+}1,\phantom{+}1,\phantom{+}1,\phantom{+}1)\\
\phantom{=}+\,&(\phantom{+}0,\phantom{+}0,\phantom{+}0,\phantom{+}0,\phantom{+}0,\phantom{+}1,\phantom{+}1,\phantom{+}1)\; .
\end{split}
\end{equation*}
}
\end{example}

\noindent$\bullet$ Corollary~\ref{th:11}(i) and Proposition~\ref{th:10}(iv), restated in dimensionality~$2^t$, suggest the following:

\begin{theorem} Let $\boldsymbol{R}$ be the symmetric cycle in the hypercube graph $\boldsymbol{H}(2^t,2)$ on the vertex set~$\{1,-1\}^{2^t}$, defined by~{\rm(\ref{eq:66})(\ref{eq:67})}.

Let $\mathcal{A}\subset\mathbf{2}^{[t]}$ be a nontrivial clutter on the ground set~$E_t$, and let\/ $\mathcal{B}:=\mathfrak{B}(\mathcal{A})$ be its blocker. Since the characteristic topes of the increasing families~$\mathcal{A}^{\triangledown}$ and $\mathcal{B}^{\triangledown}$ obey the relation
\begin{equation*}
T_{\mathcal{B}^{\triangledown}}=\ro(T_{\mathcal{A}^{\triangledown}})\; ,
\end{equation*}
we have{\rm:}
\begin{itemize}
\item[\rm(i)]
\begin{equation*}
\mathfrak{q}(T_{\mathcal{B}^{\triangledown}}):=|\boldsymbol{Q}(T_{\mathcal{B}^{\triangledown}},\boldsymbol{R})|=
|\boldsymbol{Q}(T_{\mathcal{A}^{\triangledown}},\boldsymbol{R})|=:\mathfrak{q}(T_{\mathcal{A}^{\triangledown}})\; ,
\end{equation*}
and
\begin{equation*}
\boldsymbol{x}(T_{\mathcal{B}^{\triangledown}}
)=\boldsymbol{x}(T_{\mathcal{A}^{\triangledown}}
)
\cdot\overline{\mathbf{U}}(2^t)\cdot\overline{\mathbf{T}}(2^t)\; .
\end{equation*}

\item[\rm(ii)]
Suppose that the subset $(T_{\mathcal{A}^{\triangledown}})^-=\supp(\boldsymbol{\gamma}(\mathcal{A}^{\triangledown}))\subset E_{2^t}$ is a disjoint union
\begin{equation*}
[i_1,j_1]\;\dot\cup\;[i_2,j_2]\;\dot\cup\;\cdots\;\dot\cup\;[i_{\varrho-1},j_{\varrho-1}]\;\dot\cup\;[i_{\varrho},j_{\varrho}]
\end{equation*}
of intervals such that
\begin{equation*}
j_1+2\leq i_2,\ \ j_2+2\leq i_3,\ \ \ldots,\ \
j_{\varrho-2}+2\leq i_{\varrho-1},\ \ j_{\varrho-1}+2\leq i_{\varrho}\; ,
\end{equation*}
for some $\varrho$. We have
\begin{equation*}
\mathfrak{q}(T_{\mathcal{B}^{\triangledown}})=\mathfrak{q}(T_{\mathcal{A}^{\triangledown}})=2\varrho-1\; ;\\
\end{equation*}
\begin{align*}
\boldsymbol{x}(T_{\mathcal{A}^{\triangledown}})&=\sum_{1\leq k\leq\varrho-1}\boldsymbol{\sigma}(j_k+1)-\sum_{1\leq \ell\leq\varrho}\boldsymbol{\sigma}(i_{\ell})\; ,\\
\intertext{and}
\boldsymbol{x}(T_{\mathcal{B}^{\triangledown}})&=
\sum_{1\leq k\leq\varrho-1}\boldsymbol{\sigma}(2^t-j_k+1)-\sum_{1\leq \ell\leq\varrho}\boldsymbol{\sigma}(2^t-i_{\ell}+2)\; .
\end{align*}
\end{itemize}
\end{theorem}

See~expressions~(\ref{eq:74})(\ref{eq:75}) and~(\ref{eq:82})(\ref{eq:83}) in~Example~\ref{th:9}.

\subsection{A clutter $\{\{a\}\}$} $\quad$

As earlier (in Section~\ref{section:1}), let $\{\{a\}\}$ be a clutter on the ground set~$E_t$, whose only member is a {\em one-element\/} subset $\{a\}\subset E_t$.

\noindent$\bullet$
Let us associate with the characteristic tope $\boldsymbol{\mathfrak{a}}(a):=T_{\{\{a\}\}^{\triangledown}}$ of the principal increasing family~$\{\{a\}\}^{\triangledown}$ the row vector $\boldsymbol{x}(\boldsymbol{\mathfrak{a}}(a)):=\boldsymbol{x}(\boldsymbol{\mathfrak{a}}(a),\boldsymbol{R})$ \mbox{$\in\{-1,0,1\}^{2^t}$,} described in~(\ref{eq:70}), where $\boldsymbol{R}$ is the symmetric cycle in the hypercube graph $\boldsymbol{H}(2^t,2)$, defined by~(\ref{eq:66})(\ref{eq:67}). Recall that
\begin{equation}
\label{eq:72}
\boldsymbol{x}(\boldsymbol{\mathfrak{a}}(a))=\boldsymbol{\mathfrak{a}}(a)\cdot\mathbf{M}^{-1}\; ,
\end{equation}
where the matrix~$\mathbf{M}$ is defined by~$(\ref{eq:73})$, and
\begin{equation*}
\boldsymbol{Q}(\boldsymbol{\mathfrak{a}}(a),\boldsymbol{R}):=\{x_i\cdot R^{i-1}\colon
x_i\neq 0\}\; ,\ \ \ \text{and}\ \ \ \boldsymbol{\mathfrak{a}}(a)=\sum_{Q\in\boldsymbol{Q}(\boldsymbol{\mathfrak{a}}(a),\boldsymbol{R})}Q\; ,
\end{equation*}
see~Example~\ref{th:18}.

\begin{example}
\label{th:18}
Suppose $t:=3$, and~$E_t=\{1,2,3\}$.

The characteristic topes associated with the principal increasing families that are generated by the clutters~$\{\{a\}\}$, for the elements $a\in E_t$ of the ground set, are as follows{\rm:}
{\small
\begin{align*}
\boldsymbol{\mathfrak{a}}(1):=\boldsymbol{\mathfrak{a}}(1;2^t):=T_{\{\{1\}\}^{\triangledown}}
&=(1,-1,\phantom{-}1,\phantom{-}1,-1,-1,\phantom{-}1,-1)\in\{1,-1\}^{2^t}\; ,\\
\boldsymbol{\mathfrak{a}}(2):=\boldsymbol{\mathfrak{a}}(2;2^t):=T_{\{\{2\}\}^{\triangledown}}
&=(1,\phantom{-}1,-1,\phantom{-}1,-1,\phantom{-}1,-1,-1)\; ,\\
\boldsymbol{\mathfrak{a}}(3):=\boldsymbol{\mathfrak{a}}(3;2^t):=T_{\{\{3\}\}^{\triangledown}}
&=(1,\phantom{-}1,\phantom{-}1,-1,\phantom{-}1,-1,-1,-1)\; .\\
\intertext{\normalsize{The corresponding ``$\boldsymbol{x}$-vectors'', given in~{\rm(\ref{eq:72})} for the symmetric cycle~$\boldsymbol{R}$ in the hypercube graph~$\boldsymbol{H}(2^t,2)$ defined by~{\rm(\ref{eq:66})(\ref{eq:67})}, are{\rm:}}}
\boldsymbol{x}(\boldsymbol{\mathfrak{a}}(1)):=\boldsymbol{x}(\boldsymbol{\mathfrak{a}}(1),\boldsymbol{R})
&=(0,-1,\phantom{-}1,\phantom{-}0,-1,\phantom{-}0,\phantom{-}1,-1)\in\{-1,0,1\}^{2^t}\; ,\\
\boldsymbol{x}(\boldsymbol{\mathfrak{a}}(2)):=\boldsymbol{x}(\boldsymbol{\mathfrak{a}}(2),\boldsymbol{R})
&=(0,\phantom{-}0,-1,\phantom{-}1,-1,\phantom{-}1,-1,\phantom{-}0)\; ,\\
\boldsymbol{x}(\boldsymbol{\mathfrak{a}}(3)):=\boldsymbol{x}(\boldsymbol{\mathfrak{a}}(3),\boldsymbol{R})
&=(0,\phantom{-}0,\phantom{-}0,-1,\phantom{-}1,-1,\phantom{-}0,\phantom{-}0)\; .
\end{align*}
}
Indeed, we see that
\begin{equation*}
\begin{split}
\boldsymbol{x}(\boldsymbol{\mathfrak{a}}(2))\cdot\mathbf{M}
&=
\left(
\begin{smallmatrix}
0&\phantom{-}0&-1&\phantom{-}1&-1&\phantom{-}1&-1&\phantom{-}0
\end{smallmatrix}
\right)
\cdot
\left(
\begin{smallmatrix}
\phantom{-}1& \phantom{-}1& \phantom{-}1& \phantom{-}1& \phantom{-}1& \phantom{-}1& \phantom{-}1& \phantom{-}1\\
-1& \phantom{-}1& \phantom{-}1& \phantom{-}1& \phantom{-}1& \phantom{-}1& \phantom{-}1& \phantom{-}1\\
-1& -1& \phantom{-}1& \phantom{-}1& \phantom{-}1& \phantom{-}1& \phantom{-}1& \phantom{-}1\\
-1& -1& -1& \phantom{-}1& \phantom{-}1& \phantom{-}1& \phantom{-}1& \phantom{-}1\\
-1& -1& -1& -1& \phantom{-}1& \phantom{-}1& \phantom{-}1& \phantom{-}1\\
-1& -1& -1& -1& -1& \phantom{-}1& \phantom{-}1& \phantom{-}1\\
-1& -1& -1& -1& -1& -1& \phantom{-}1& \phantom{-}1\\
-1& -1& -1& -1& -1& -1& -1& \phantom{-}1
\end{smallmatrix}
\right)\\
&=
\left(
\begin{smallmatrix}
1&\phantom{-}1&-1&\phantom{-}1&-1&\phantom{-}1&-1&-1
\end{smallmatrix}
\right)=:\boldsymbol{\mathfrak{a}}(2)\; .
\end{split}
\end{equation*}
\label{page:4}
\end{example}

\noindent$\bullet$
For the row vector $\boldsymbol{y}(1+a):=\boldsymbol{y}(1+a;2^t)\in\{-1,0,1\}^{2^t}$, defined by
\begin{equation*}
\boldsymbol{y}(1+a):=\boldsymbol{x}({}_{-\{1+a\}}\mathrm{T}_{2^t}^{(+)}
)=:\boldsymbol{x}(T_{\{\{a\}\}}
)\; ,
\end{equation*}
we have (see~\cite[Sect.~2]{M-SC-II}):
\begin{equation*}
\boldsymbol{y}(1+a)=\boldsymbol{\sigma}(1)-\boldsymbol{\sigma}(1+a)+\boldsymbol{\sigma}(2+a)\; .
\end{equation*}
In other words,
\begin{equation*}
\boldsymbol{Q}(T_{\{\{a\}\}},\boldsymbol{R})=\{\underbrace{R^0}_{\mathrm{T}_{2^t}^{(+)}}\, ,R^{1+a},R^{2^t+a}\}\; ,
\end{equation*}
and
\begin{equation*}
T_{\mathfrak{B}(\{\{a\}\})}=T_{\{\{a\}\}}=
\mathrm{T}_{2^t}^{(+)}+R^{1+a}+R^{2^t+a}\; .
\end{equation*}
Equivalently,
\begin{equation*}
\widetilde{\boldsymbol{Q}}(\boldsymbol{\gamma}(\{\{a\}\}),\widetilde{\boldsymbol{R}})=\{\underbrace{\widetilde{R}^0}_{(0,\ldots,0)},
\widetilde{R}^{1+a},\widetilde{R}^{2^t+a}\}\; ,
\end{equation*}
and
\begin{equation*}
\begin{split}
\boldsymbol{\gamma}(\mathfrak{B}(\{\{a\}\}))=\boldsymbol{\gamma}(\{\{a\}\})&=-\tfrac{1}{2}(3-1)\cdot\mathrm{T}_{2^t}^{(+)}
+\widetilde{R}^{1+a}+\widetilde{R}^{2^t+a}\\ &=-\mathrm{T}_{2^t}^{(+)}
+\widetilde{R}^{1+a}+\widetilde{R}^{2^t+a}\; .
\end{split}
\end{equation*}

\subsection{A clutter $\{A\}$} $\quad$

As in Section~\ref{section:2}, let $\{A\}$ be a clutter on the ground set~$E_t$, whose only member is a {\em nonempty subset} $A\subseteq E_t$.

\noindent$\bullet$ Dealing with the symmetric cycle~$\boldsymbol{R}$ in the hypercube graph~$\boldsymbol{H}(2^t,2)$, defined by~(\ref{eq:66})(\ref{eq:67}), with the matrix~$\mathbf{M}$ given in~(\ref{eq:73}), and with ``$\boldsymbol{x}$-vectors'' described in~(\ref{eq:70}), for the row vector
\begin{equation}
\label{eq:71}
\boldsymbol{y}(\varGamma^{-1}(A)):=\boldsymbol{x}({}_{-\{\varGamma^{-1}(A)\}}\mathrm{T}_{2^t}^{(+)})
=:\boldsymbol{x}(T_{\{A\}})=T_{\{A\}}\cdot\mathbf{M}^{-1}\in\{-1,0,1\}^{2^t}\; ,
\end{equation}
we have (see~\cite[Sect.~2]{M-SC-II}):
\begin{equation*}
\boldsymbol{y}(\varGamma^{-1}(A))=
\begin{cases}
\phantom{-}\boldsymbol{\sigma}(1)-\boldsymbol{\sigma}(\varGamma^{-1}(A))+\boldsymbol{\sigma}(1+\varGamma^{-1}(A))\; , & \text{if $A\neq E_t$},\\
\quad\vspace{-3mm}\\
-\boldsymbol{\sigma}(2^t)\; , & \text{if $A=E_t$}.
\end{cases}
\end{equation*}
In other words,
\begin{equation*}
\boldsymbol{Q}(T_{\{A\}},\boldsymbol{R})=
\begin{cases}
\{\,\underbrace{R^0}_{\mathrm{T}_{2^t}^{(+)}}\, ,R^{\varGamma^{-1}(A)},R^{2^t+\varGamma^{-1}(A)-1}\}\; , & \text{if $A\neq E_t$},\\
\quad\vspace{-3mm}\\
\{R^{2\cdot 2^t -1}\}\; , & \text{if $A=E_t$},
\end{cases}
\end{equation*}
and
\begin{equation*}
T_{\{A\}}=
\begin{cases}
\mathrm{T}_{2^t}^{(+)}+R^{\varGamma^{-1}(A)}+R^{2^t+\varGamma^{-1}(A)-1}\; , & \text{if $A\neq E_t$},\\
\quad\vspace{-3mm}\\
R^{2\cdot 2^t -1}\; , & \text{if $A=E_t$}.
\end{cases}
\end{equation*}
Equivalently,
\begin{equation*}
\widetilde{\boldsymbol{Q}}(\boldsymbol{\gamma}(\{A\}),\widetilde{\boldsymbol{R}})=
\begin{cases}
\{\underbrace{\widetilde{R}^0}_{(0,\ldots,0)},\widetilde{R}^{\varGamma^{-1}(A)},\widetilde{R}^{2^t+\varGamma^{-1}(A)-1}\}\; , & \text{if $A\neq E_t$},\\
\quad\vspace{-3mm}\\
\{\widetilde{R}^{2\cdot 2^t -1}\}\; , & \text{if $A=E_t$},
\end{cases}
\end{equation*}
and
\begin{equation*}
\boldsymbol{\gamma}(\{A\})=
\begin{cases}
\phantom{=}-\tfrac{1}{2}(3-1)\cdot\mathrm{T}_{2^t}^{(+)}+\widetilde{R}^{\varGamma^{-1}(A)}+\widetilde{R}^{2^t+\varGamma^{-1}(A)-1}\\
=-\mathrm{T}_{2^t}^{(+)}+\widetilde{R}^{\varGamma^{-1}(A)}+\widetilde{R}^{2^t+\varGamma^{-1}(A)-1}\; , & \text{if $A\neq E_t$},\\
\quad\vspace{-3mm}\\
\phantom{=-}\widetilde{R}^{2\cdot 2^t -1}\; , & \text{if $A=E_t$}.
\end{cases}
\end{equation*}

\noindent$\bullet$
Recall that $\boldsymbol{\gamma}(\{A\}^{\triangledown})=\sideset{}{^\ast}\prod_{a\in A}\widetilde{\boldsymbol{\mathfrak{a}}}(a)$, and $\boldsymbol{\gamma}(\mathfrak{B}(\{A\})^{\triangledown})=\rn(\boldsymbol{\gamma}(\{A\}^{\triangledown}))$.

\begin{remark}[cf.~Remark~\ref{th:8}]
For a nonempty subset~$A\subseteq E_t$, we have
\begin{itemize}
\item[\rm(i)]
\begin{equation*}
\boldsymbol{\gamma}(\{A\}^{\triangledown})=\sideset{}{^\ast}\prod_{a\in A}\bigl(\tfrac{1}{2}\bigr(\mathrm{T}_{2^t}^{(+)}-\boldsymbol{x}(\boldsymbol{\mathfrak{a}}(a))
\cdot\mathbf{M}\bigr)\bigr)\; .
\end{equation*}

\item[\rm(ii)]
\begin{equation*}
\boldsymbol{\gamma}(\mathfrak{B}(\{A\})^{\triangledown})
=\mathrm{T}_{2^t}^{(+)}
-\Bigl(\sideset{}{^\ast}\prod_{a\in A}\bigl(\tfrac{1}{2}\bigl(\mathrm{T}_{2^t}^{(+)}-\boldsymbol{x}(\boldsymbol{\mathfrak{a}}(a))\cdot\mathbf{M}\bigr)\bigr)\Bigr)
\cdot\overline{\mathbf{U}}(2^t)\; .
\end{equation*}
\end{itemize}
\end{remark}

\noindent$\bullet$ Since the {\em blocker\/} of the clutter $\{A\}$ is the clutter $\mathfrak{B}(\{A\})=\{\{a\}\colon a\in A\}$, and
$\boldsymbol{\gamma}(\mathfrak{B}(\{A\}))=\sum_{a\in A}\boldsymbol{\gamma}(\{\{a\}\})$, we have
\begin{equation*}
\boldsymbol{\gamma}(\mathfrak{B}(\{A\}))=\sum_{a\in A}(-\mathrm{T}_{2^t}^{(+)}+\widetilde{R}^{1+a}+\widetilde{R}^{2^t+a})\; ,
\end{equation*}
that is,
\begin{equation*}
\boldsymbol{\gamma}(\mathfrak{B}(\{A\}))=-|A|\cdot\mathrm{T}_{2^t}^{(+)}+\sum_{a\in A}(\widetilde{R}^{1+a}+\widetilde{R}^{2^t+a})\; .
\end{equation*}

\subsection{A clutter $\mathcal{A}:=\{A_1,\ldots,A_{\alpha}\}$} $\quad$

As in Section~\ref{section:3}, let $\mathcal{A}:=\{A_1,\ldots,A_{\alpha}\}$ be a nontrivial clutter on the ground set $E_t$.

\noindent$\bullet$ In analogy with~\cite[Rem.~2.2]{M-SC-II}, dealing with the symmetric cycle~$\boldsymbol{R}$ in the hypercube graph~$\boldsymbol{H}(2^t,2)$, defined by~(\ref{eq:66})(\ref{eq:67}), with the matrix~$\mathbf{M}$ given in~(\ref{eq:73}), with ``$\boldsymbol{x}$-vectors'' described in~(\ref{eq:70}), and with ``$\boldsymbol{y}$-vectors'' given in~(\ref{eq:71}), we have
\begin{equation*}
\boldsymbol{x}(T_{\mathcal{A}})=(1-\#\mathcal{A})\cdot\boldsymbol{\sigma}(1)+\sum_{A\in\mathcal{A}}\boldsymbol{y}(\varGamma^{-1}(A))\; ,
\end{equation*}
that is,
\begin{equation*}
\boldsymbol{x}(T_{\mathcal{A}})=
\begin{cases}
\phantom{-}\boldsymbol{\sigma}(1)+\sum_{A\in\mathcal{A}}(-\boldsymbol{\sigma}(\varGamma^{-1}(A))+\boldsymbol{\sigma}(1+\varGamma^{-1}(A)))\; , & \text{if $\mathcal{A}\neq\{E_t\}$},\\
\quad\vspace{-3mm}\\
-\boldsymbol{\sigma}(2^t)\; , & \text{if $\mathcal{A}=\{E_t\}$},
\end{cases}
\end{equation*}
or
\begin{equation*}
T_{\mathcal{A}}=
\begin{cases}
\phantom{-}\mathrm{T}_{2^t}^{(+)}+\sum_{A\in\mathcal{A}}(R^{\varGamma^{-1}(A)}+R^{2^t+\varGamma^{-1}(A)-1})\; , & \text{if $\mathcal{A}\neq \{E_t\}$},\\
\quad\vspace{-3mm}\\
\phantom{-}R^{2\cdot 2^t -1}\; , & \text{if $\mathcal{A}= \{E_t\}$}.
\end{cases}
\end{equation*}

We also have
\begin{equation*}
\boldsymbol{\gamma}(\mathcal{A})=
\begin{cases}
-(\#\mathcal{A})\cdot\mathrm{T}_{2^t}^{(+)}+\sum_{A\in\mathcal{A}}(\widetilde{R}^{\varGamma^{-1}(A)}+\widetilde{R}^{2^t+\varGamma^{-1}(A)-1})\; , & \text{if $\mathcal{A}\neq \{E_t\}$},\\
\quad\vspace{-3mm}\\
\phantom{-}\widetilde{R}^{2\cdot 2^t -1}\; , & \text{if $\mathcal{A}= \{E_t\}$}.
\end{cases}
\end{equation*}

\noindent$\bullet$ Theorem~\ref{th:5} can be accompanied with the following statement:
\begin{corollary}
If $\mathcal{A}:=\{A_1,\ldots,A_{\alpha}\}$ is a nontrivial clutter on the ground set~$E_t$, then we have{\rm:}
\begin{itemize}
\item[\rm(i)]
\begin{equation*}
\boldsymbol{\gamma}(\mathcal{A}^{\triangledown})=
\mathrm{T}_{2^t}^{(+)}-
\Bigl(\sideset{}{^\ast}\prod_{i\in[\alpha]}\Bigl(\mathrm{T}^{(+)}_{2^t}
-\sideset{}{^\ast}\prod_{a^i\in A_i}\bigl(\tfrac{1}{2}\bigl(\mathrm{T}_{2^t}^{(+)}-\boldsymbol{x}(\boldsymbol{\mathfrak{a}}(a^i))
\cdot\mathbf{M}\bigr)\bigr)\Bigr)\Bigr)\; .
\end{equation*}
\item[\rm(ii)]
\begin{equation*}
\boldsymbol{\gamma}(\mathfrak{B}(\mathcal{A})^{\triangledown})
=\Bigl(\sideset{}{^\ast}\prod_{i\in[\alpha]}\Bigl(\mathrm{T}^{(+)}_{2^t}
-\sideset{}{^\ast}\prod_{a^i\in A_i}\bigl(\tfrac{1}{2}\bigl(\mathrm{T}_{2^t}^{(+)}-\boldsymbol{x}(\boldsymbol{\mathfrak{a}}(a^i))
\cdot\mathbf{M}\bigr)\bigr)\Bigr)\Bigr)\cdot\overline{\mathbf{U}}(2^t)
\; .
\end{equation*}
\end{itemize}
\end{corollary}

\vspace{5mm}
\end{document}